\documentclass[a4paper,oneside,11pt,reqno]{amsart}
\usepackage[utf8]{inputenc}

\usepackage{amsfonts}
\usepackage{amsmath}
\usepackage{amssymb}
\usepackage{amsthm}
\usepackage{geometry}
 \usepackage{lineno}
\usepackage{mathrsfs} 
\usepackage[dvipsnames]{xcolor}

\usepackage{bbm}
\usepackage{mathtools}

\usepackage{dsfont}

\usepackage[pagebackref=true,colorlinks=true, linkcolor=Blue, citecolor=Blue]{hyperref}
\renewcommand*\backref[1]{\ifx#1\relax \else (Cited on #1) \fi}

\usepackage{appendix}
\usepackage{enumerate}
\usepackage{float}

\usepackage{tikz}
\usepackage{tikz-cd} 
\usetikzlibrary{arrows, arrows.meta}

\usepackage{zref-clever}
\zcsetup{cap}
\zcRefTypeSetup{equation}{
Name-sg = {} ,
name-sg = {},
Name-pl = {},
name-pl = {}
}
\zcRefTypeSetup{assumptions}{
    Name-sg = Assumptions,
    name-sg = assumptions,
    Name-pl = Assumptions,
    name-pl = assumptions,
}
\zcsetup{nameinlink=true}

\AddToHook{env/remark/begin}{%
\zcsetup{countertype={definition=remark}}}
\AddToHook{env/conjecture/begin}{%
\zcsetup{countertype={definition=conjecture}}}
\AddToHook{env/corollary/begin}{%
\zcsetup{countertype={definition=corollary}}}
\AddToHook{env/theorem/begin}{%
\zcsetup{countertype={definition=theorem}}}
\AddToHook{env/lemma/begin}{%
\zcsetup{countertype={definition=lemma}}}
\AddToHook{env/proposition/begin}{%
\zcsetup{countertype={definition=proposition}}}
\AddToHook{env/example/begin}{%
\zcsetup{countertype={definition=example}}}
\AddToHook{env/assumptions/begin}{%
\zcsetup{countertype={definition=assumptions}}}

\usepackage{scalerel}[2016/12/29]

\theoremstyle{plain}
\newtheorem{definition}{Definition}
\newtheorem{proposition}[definition]{Proposition}
\newtheorem{lemma}[definition]{Lemma}
\newtheorem{corollary}[definition]{Corollary}
\newtheorem{theorem}[definition]{Theorem}
\newtheorem{remark}[definition]{Remark}

\newtheorem{assumptions}[definition]{Assumptions}

\numberwithin{definition}{section}
\numberwithin{equation}{section}

\renewcommand*{\P}{\mathbb{P}}

\newcommand*{\R}{\mathbb{R}}

\newcommand*{\N}{\mathbb{N}}
\newcommand*{\F}{\mathscr{F}}

\newcommand*{\Mcal}{\mathcal{M}}

\newcommand*{\Pcal}{\mathcal{P}}
\newcommand*{\Gcal}{\mathcal{G}}

\newcommand*{\Lcal}{\mathcal{L}}

\newcommand*{\Ncal}{\mathcal{N}}

\newcommand{\norm}[1]{\left \lVert  #1 \right \rVert}
\newcommand{\abs}[1]{\left\lvert #1 \right\rvert}
\newcommand*{\Tl}{T\Ssup{\Lambda}}
\newcommand*{\Xl}{X\Ssup{\Lambda}}

\newcommand*{\GG}{\mathscr{G}}

\newcommand*{\RR}{\mathscr{R}}
\newcommand*{\LL}{\mathscr{L}}

\newcommand*{\1}{\mathds{1}}

\renewcommand*{\L}{\Lambda}

\renewcommand*{\d}{\mathrm{d}}
\newcommand*{\e}{\mathrm{e}}

\newcommand*{\SpecEnt}{\mathscr{I}}
\newcommand*{\RelEnt}{I}

\newcommand{\Ssup}[1]{^{\scriptscriptstyle{({#1}})}}

\newcommand{\birthsinL}[1]{B_\L(#1)}
\newcommand{\deathsinL}[1]{D_\L(#1)}

\newcommand*{\sepset}{\colon}

\title[Reversible birth-and-death dynamics in continuum]{Reversible birth-and-death dynamics in continuum: \\ free-energy dissipation and attractor properties}
\author[B. Jahnel]{Benedikt Jahnel}
\author[J. K\"oppl]{Jonas K\"oppl}
\author[Y. Steenbeck]{Yannic Steenbeck}
\author[A. Zass]{Alexander Zass}

\address[Benedikt Jahnel]{TU Braunschweig, Institut für Mathematische Stochastik, 
Germany, and Weierstrass Institute, Berlin, Germany.}
\email{jahnel@wias-berlin.de}
\address[Jonas K\"oppl]{Universität Münster,
Germany.}
\email{jonas.koeppl@uni-muenster.de}
\address[Yannic Steenbeck]{TU Braunschweig, Institut für Mathematische Stochastik, 
Germany.}
\email{yannic.steenbeck@tu-braunschweig.de}
\address[Alexander Zass]{Weierstrass Institute, Berlin,
Germany.}
\email{zass@wias-berlin.de}
\date{\today}


\date{\today}

\definecolor{amethyst}{rgb}{0.6, 0.4, 0.8}

\usepackage[textwidth=20mm,textsize=small]{todonotes}

\usepackage{soul}

\begin{document}

\begin{abstract}
We consider continuous-time birth-and-death dynamics in $\R^d$ admitting infinite-volume Gibbs point processes as reversible measures. The birth rates are translation invariant, have finite range, and are uniformly bounded from above and away from zero. For all translation-invariant initial measures, we show that the specific relative entropy is non-increasing along trajectories and that, if the entropy is initially finite, then every long-time weak limit point is a Gibbs point process for the same interaction. Our proofs rest on a lower bound for entropy dissipation in terms of the reduced Palm measure of the evolved law. 
\end{abstract}

\maketitle



\section{Introduction}

The theory of Gibbs point processes, describing interacting systems of particles in continuous space, is by now a classical topic in probability and statistical mechanics and, for various questions, satisfactory results have been obtained: to list a few non-exhaustive examples, existence~\cite{ruelle_1969,Dereudre_Vasseur_2020,Roelly_Zass_2020,Betsch2023a}, uniqueness~\cite{ruelle_1969,Fernndez2007,Houdebert2022,Betsch2023}, phase transitions \cite{Ruelle1971WRPT,Giacomin1995Agreement, CCK95,Georgii_1996,Dereudre2025,DRC26,bowen_2026}, the equivalence of various ensembles \cite{georgii_equivalence_1995}, and the Gibbs variational principle~\cite{georgii_equivalence_1995,dereudre2009variational,JKSZ24}.
On the lattice, these results on the equilibrium behavior of large systems of interacting particles are complemented by a detailed study of associated dynamics, see for example the classical reference \cite{liggett_interacting_2005}.
In contrast, very little is known about dynamical aspects of infinite systems in the continuum, in particular out-of-equilibrium dynamics and convergence to equilibrium. 

\subsection{Well-defined dynamics} 
The first issue lies in the fact that the unbounded nature of point processes in the continuum means that already establishing a well-defined dynamics proves to be a major challenge, and cannot be recovered by applying Liggett's classical results, cf.~\cite{penrose2008spatial}. Inspired by biologically motivated models (e.g., locally regulated populations), the study of birth-and-death processes in the continuum, with birth and death rates depending on the configuration of the system, was initiated in~\cite{Preston1975} for bounded systems. Building on this, the seminal work~\cite{HS78} shows the well-posedness of the dynamics for a one-dimensional system, as well as a characterization of its reversible measures, using a martingale-problem approach. Extending their results to more general cases is then an intrinsically non-trivial matter, and a few approaches have been explored in the literature: the first uses Dirichlet forms to construct equilibrium dynamics and characterize stationary measures as Gibbs states via their correlation functions, see e.g.~\cite{Kondratiev2005spectralgap,Kond2008,Finkelshtein2014Dynamical}; the second characterizes the Markov process as a solution of a martingale problem, see e.g. ~\cite{Kurtz1980Representations,Garcia1995Birth,garcia2006spatial}; finally, it is worth noting that the lookdown-construction approach presented in~\cite{Etheridge2019Genealogical} covers a broad class of spatial birth-and-death processes, with the aim of capturing the genealogy of a population. 

\subsection{Out-of-equilibrium dynamics}\label{sec:finite-state-space-example}
While the above works are able to construct well-defined dynamics and investigate key properties of the processes such as characterizations of the stationary measures, the long-time behavior of the dynamics, and in particular their weak limit points, remains much less understood than on the lattice setting. In the continuum, non-equilibrium dynamics have been constructed by considering time-evolved correlation functions~\cite{Kondratiev2006nonequilibrium,Kondratiev2008nonequilibrium}. Relaxation to equilibrium has been previously investigated in the finite-volume setting using spectral methods: in~\cite{Bertini2002spectralgap,Wu2004spectralgap}, the authors establish positive spectral-gap bounds for finite-volume Glauber dynamics and explore its connection to exponential convergence to equilibrium in the uniqueness regime. In the infinite-volume setting, exponential convergence to the unique equilibrium was established in~\cite{Kondratiev2010nonequilibrium}. However, these results do not address the characterization of the weak limit points of the dynamics in regimes of phase coexistence, as considered in this work.

We investigate the long-time behavior of Markov processes out of equilibrium with an approach based on relative-entropy dissipation, motivated by the following simple observation for continuous-time Markov chains on finite state spaces. For concreteness, let us denote the state space by $E$ and consider an irreducible generator $\mathscr{L}$. In this situation, it is well known that there exists a unique measure $\mu$ that is time-stationary with respect to the Markov semigroup $(P_t)_{t \geq 0}$ generated by $\mathscr{L}$. Because of the irreducibility, $\mu$ puts strictly positive mass on every state $x \in E$ and we can define the relative entropy of another probability measure $\nu$ with respect to $\mu$ by 
\begin{align*}
    \RelEnt(\nu \lvert \mu) = \sum_{x \in E}\nu(x) \log \frac{\nu(x)}{\mu(x)},
\end{align*}
where we use the convention that $0\log 0 = 0$. For an initial distribution $\nu$ let us denote the distribution at time $t\geq 0$ by $\nu_t$. By Jensen's inequality one can then easily check that for any $\nu$ the map $t \mapsto \RelEnt(\nu_t \lvert \mu)$ is non-increasing and is constant if and only if $\nu = \mu$. This tells us that the functional $\nu \mapsto \RelEnt(\nu \lvert \mu)$ is a strict Lyapunov function for the unique fixed point $\mu$ of the measure-valued ODE 
\begin{align*}
    \partial \nu_t = \nu_t \mathscr{L}.
\end{align*}
This directly implies that $(\nu_t)_{t \geq 0}$ converges to $\mu$ as $t$ tends to infinity for any initial distribution $\nu$. Of course, this is nothing but the classical ergodic theorem for finite-state Markov processes, but this proof shows that the convergence to the unique time-stationary measure also fits precisely into the physical picture of convergence to equilibrium via relative-entropy dissipation. 
Moreover, it can also be adapted to more complex situations. 
In the context of classical interacting particle systems on the integer lattice, this strategy has been successfully applied, first to reversible dynamics \cite{holley_free_1971,higuchi_results_1975} and later also to non-reversible systems \cite{Sullivan1976,kunsch_time_1984,handa1996entropy,jahnel_attractor_2019, jahnel_dynamical_2023,JK25}. 
Finally, it is worth noting that, in a related direction, gradient-flow approaches to infinite-volume dynamics have recently been developed in~\cite{Erbar2025,Herry2025}, with specific relative entropy and free energy, respectively, as the underlying functionals.

\subsection{Contributions of this work}
Our goal is to investigate the long-time behavior of infinite-volume birth-and-death dynamics in continuum space through relative-entropy dissipation and characterize the long-time weak limit points. In particular, our main contribution is two-fold. First, we show, in \zcref{Theorem:Main_Theorem_Decrease_in_Entropy} a free-energy dissipation result which, to the best of our knowledge, so far existed only in the non-interacting case, e.g.~\cite{Dello2024Wasserstein,huesmann2025}.
Second, in \zcref{Theorem:Main_Theorem_Limit_Points}, we apply this free-energy dissipation result to establish the attractor property of the dynamics for starting measures of finite specific relative entropy, in the sense that any of its long-time weak limit points is a Gibbs point process w.r.t.\ the considered interaction.

On the way, we prove various technical results for birth-and-death dynamics in continuous space, which could be of independent interest. Most importantly, we show that there is a finite speed of propagation, i.e., with very high probability no information about the state of the process at a given position or its boundary condition can travel faster than linearly, see \zcref{Lemma:finite_speed_of_propagation_version_1}. We also provide, in \zcref{sec:construction_dynamics} and \zcref{sec:appendix:infiniteconstruction}, a reformulation of the construction of the infinite-volume dynamics of~\cite{garcia2006spatial} in terms of measurable solution maps.

While we focus on a birth-and-death process with translation-invariant and finite-range birth rates that are uniformly bounded from above and away from zero, we believe the foundations laid here could serve as a stepping stone for future work, see \zcref{sec:outlook}. 

\subsection{Outline of the paper}

The rest of the paper is structured as follows: in \zcref{sec:setting_results}, after introducing the notation and setting, we present our main results in \zcref{sec:results} and discuss possible future directions in \zcref{sec:outlook}.
We lay out the proof strategy and prove our two main theorems in~\zcref{section:proof-strategy}, with the proofs split up into various auxiliary results provided in 
\zcref[S,range]{sec:construction_dynamics,sec:time-stationary} and the appendices. We start by providing the construction of the infinite-volume dynamics in \zcref{sec:construction_dynamics}, which we then approximate by a process that only acts locally in \zcref{section:finite_volume_dynamics}, and estimate the thereby introduced error. In \zcref{sec:reversible-Gibbs}, we show a useful equivalence between reversible measures and Gibbs measures. \zcref{sec:relative-entropy-decreasing} provides the entropy estimates and finite-volume comparison used to prove monotonicity of the specific relative entropy. In \zcref{sec:attractor-property}, we establish a finite-volume entropy-dissipation identity and bounds on Fisher information. Finally, in \zcref{sec:time-stationary}, we show that the set of translation-invariant Gibbs measures coincides with the set of translation-invariant time-stationary measures.

\section{Setting and main result}\label{sec:setting_results}

We consider Markov processes on the space $\Omega$ of locally finite counting measures (point configurations) in $\R^d$, $d\ge 1$, where several points may sit at the same location. 

The space $\R^d$ is endowed with the Euclidean norm $\abs{\cdot}$ and the associated Borel $\sigma$-algebra $\mathfrak{B}(\R^d)$. We write $\L\Subset\R^d$, or $\L\in\mathfrak{B}_b(\mathbb{R}^d)$, if $\L$ is a bounded Borel subset of $\R^d$, set $\L_n := [-n/2,n/2)^d$, for $n\in\N:=\{1,2,\ldots\}$, and call \emph{box} any set of the form \([a_1, b_1) \times \dots \times [a_d, b_d)\Subset\R^d\), $a_i,b_i\in\R$. We denote by $B_R(x)$ the ball with radius $R>0$ centered at $x\in \R^d$. Let $o$ denote the origin of $\R^d$.

We will often use the short-hand notation $f\lesssim g$ if there exists a finite constant $c>0$, independent of $f$ and $g$, such that $f\le c \cdot g$. We write $f\lesssim_A g$ if $c=c(A)$ depends on some parameter $A$, for example a subset $\L\subset\R^d$.

\subsection{Point process formalism}

An element $\eta\in\Omega$ is of the form
\begin{equation*}
    \eta = \sum_{i\ge 1} \delta_{x_i}
    \equiv \{x_1,x_2,\dots\},\quad x_i\in\R^d,
\end{equation*}
where the set notation is to be intended as a multiset in order to allow for multiplicities. We denote by $\Omega_{\mathrm{s}}$ the subspace of simple configurations, in which multiplicities are not allowed, and by $\Omega_{\mathrm{f}}$ the space of finite point configurations. 

A function $f\colon\Omega \to \R$ is said to be \emph{local} if there exists $\Lambda\Subset\R^d$ such that $f(\eta) = f(\eta_\Lambda)$ for all $\eta\in\Omega$; in this case, we also say that $f$ is \emph{$\Lambda$-local}. We write $\eta_\L := \eta\cap\L$, with  $\L\subset\R^d$, and correspondingly denote by $\Omega_\L$ the space of point configurations in $\L$. For $\L\Subset\R^d$, let $\eta\mapsto N_\L(\eta)=|\eta_\L|$ denote the {\em number of points} of $\eta$ in $\L$. Given two point configurations $\eta = \sum_{i\geq1}\delta_{x_i},\zeta = \sum_{j\geq 1}\delta_{y_j}\in\Omega$, we denote by $\eta\zeta \equiv \eta\cup\zeta:=\sum_{i\geq 1}\delta_{x_i}+\sum_{j\geq1}\delta_{y_j}$ their concatenation. We write $\theta_z$ for the translation by a vector $z\in\R^d$, acting either on points by $x\mapsto x-z$ or on point configurations by $\sum_{i\geq 1}\delta_{x_i}\mapsto \sum_{i\geq 1}\delta_{\theta_z x_i}$.
We call \emph{point process} any probability measure on the set of point configurations, and denote the set of all such measures by $\Pcal(\Omega)$. The restriction of $\mu\in\Pcal(\Omega)$ to $\L\Subset\R^d$ is denoted by $\mu_\L\in\Pcal(\Omega_\L)$. Let $\Pcal_\theta(\Omega)\subset\Pcal(\Omega)$ denote the subset of \emph{translation-invariant} probability measures, i.e.\ elements $\mu\in\Pcal(\Omega)$ such that $\mu\circ\theta_z^{-1} = \mu$ for all $z\in\R^d$. We denote by \(\Pcal_1(\Omega)\) the set of all \(\mu \in \Pcal(\Omega)\) with a locally finite intensity measure, i.e.\ with \(\mu[N_\Lambda] < \infty\) for all \(\Lambda \Subset \R^d\), and set \(\Pcal_{1, \theta}(\Omega) = \Pcal_1(\Omega) \cap \Pcal_\theta(\Omega)\). 

Let $C_c(\R^d)$ denote the space of continuous, compactly supported functions from $\R^d$ to $\R$. We endow $\Omega$ with the vague topology, i.e.\ the coarsest topology on $\Omega$ that makes $\eta\mapsto \int g\, \d\eta$ continuous for $g\in C_c(\R^d)$. Similarly, we endow $\Omega_{\mathrm{s}}$ with the induced topology, so that both spaces are Polish. The corresponding Borel \(\sigma\)-algebra is generated by the local counting variables \(\{N_\Lambda \colon \Lambda \Subset \R^d\}\). Additionally, let \(\widehat{\Omega}\) be the space of counting measures \(\widehat{\eta}\) on \(\R^d \times [0, \infty)\) with \(\widehat{\eta}(\cdot \times [0, \infty)) \in \Omega\). It is again equipped with the induced subspace topology w.r.t.\ the locally finite counting measures on \(\R^d \times [0, \infty)\) equipped with the vague topology, and its Borel sigma algebra is generated by the local counting variables \(\{N_{\Lambda \times F} \colon \Lambda \Subset \R^d, F \Subset [0, \infty)\}\). See also~\cite[Chapter 9]{DaleyVereJones2008} and~\cite[Chapters 1 and 4]{Kal17} for the relevant topological and measure-theoretic background. For a function $f:\R^d\times\Omega \to \R$, we write \(\Vert f\Vert_\infty := \sup_{x \in \R^d, \eta \in \Omega} \vert f(x, \eta)\vert\).

Let \(E := D([0, \infty); \Omega)\) denote the Skorokhod space of \(\Omega\)-valued c\`adl\`ag paths and \(E_{\mathrm{s}} := D([0, \infty); \Omega_s)\) the Skorokhod space of \(\Omega_s\)-valued c\`adl\`ag paths, both equipped with the corresponding Skorokhod \(J_1\) topology, making them Polish spaces, cf.\ \cite[Chapter 3, Section 5]{EK86}. Similarly, let \(\widehat{E} := D([0, \infty); \widehat{\Omega})\) be the Skorokhod space of \(\widehat{\Omega}\)-valued c\`adl\`ag paths, also equipped with the Skorokhod \(J_1\) topology.

\subsection{Gibbs point processes}

Given an \emph{energy functional} $H$, i.e.\ a translation-invariant measurable function $H:\Omega_{\mathrm{f}}\to \R\cup\{+ \infty\}$,
we are interested in the translation-invariant {\em Gibbs point processes} associated to $H$, which are defined as follows. Let $\pi$ denote the homogeneous intensity-one {\em Poisson point process} on $\Omega$, and define the {\em finite-volume Gibbs point process} in $\L\Subset\R^d$ with {\em boundary condition} $\omega\in \Omega$ as the following probability measure on $\Omega_\L$:
\begin{equation*}
    \nu\Ssup{\omega}_\L(\d\eta_\L):=Z_\L(\omega)^{-1}\e^{-H_\L(\eta_\L\omega_{\L^c})}\pi_\L(\d \eta_\L),
\end{equation*}
where the \emph{partition function} $Z_\L(\omega)$ is the usual normalization constant, and 
\begin{equation*}
    H_\L(\eta):=\lim_{n\uparrow\infty}\big(H(\eta_{\L_n})-H(\eta_{\L_n\setminus\L})\big)
\end{equation*} 
is the \emph{(conditional) energy} in $\Lambda$, where the limit is understood pointwise for every $\eta\in\Omega$.
Under the finite-range assumption of \zcref{assumptions}(iii) below, this limit exists since the term in the right-hand side eventually becomes independent of $n$.

We say that a probability measure $\nu$ on $\Omega$ is an {\em infinite-volume Gibbs point process} associated to the energy functional $H$ if it satisfies the {\em DLR equations} (named after Dobrushin, Lanford, and Ruelle,~\cite{dobrushin_1968,lanford_ruelle_1969}): 
\begin{equation*}
    \int\nu(\d \omega)f(\omega)=\int \nu(\d \omega) \int\nu_\L\Ssup{\omega}(\d \eta_\L)\, f(\eta_\L\omega_{\L^c}),\quad \Lambda \Subset \R^d, f\ge0 \text{ and measurable}.
\end{equation*}
Our \zcref{assumptions} below will guarantee that the set of translation-invariant Gibbs point processes $\GG_\theta$ is non-empty. Moreover, the question of \emph{phase transitions}, or non-uniqueness of the Gibbs measures is a central topic in the field, and one of the features of our main result is that it works in both the uniqueness as well as the phase-transition regime, i.e.\ when the set of translation-invariant infinite-volume Gibbs point processes $\GG_\theta$ contains more than one element. 

\subsection{Birth-and-death processes}
Based on the equilibrium setting just described, we wish to consider Markovian dynamics that are reversible with respect to elements of $\GG_\theta$. For this purpose, let 
\begin{equation*}
    h(x,\eta):= \lim_{n\uparrow\infty}\big(H(\eta_{\L_n}+\delta_x) - H(\eta_{\L_n})\big)
\end{equation*} 
be the {\em conditional energy} of a point $x\in\R^d$ in a configuration $\eta\in\Omega$, and define the {\em birth rate}, based on $h$, at $x$ given $\eta$ as
\begin{equation*}
    b(x,\eta) := \e^{-h(x,\eta)}.
\end{equation*} 
This is the \emph{Papangelou intensity} associated to the energy functional $H$, as the infinite-volume Gibbs point processes are characterized by the \emph{GNZ equations} (named after Georgii, Nguyen, and Zessin,~\cite{Georgii1976Canonical,xanh1979integral}):
\begin{equation*}
    \int \nu(\d\eta) \int \eta(\d x) f(x,\eta-\delta_x) = \int \nu(\d\eta) \int_{\R^d} \d x\, f(x,\eta)\, \e^{-h(x,\eta)},\quad f:\R^d\times\Omega\to [0, \infty) \text{ measurable}.
\end{equation*}
We will always work under the following
\begin{assumptions}\label{assumptions}
    The birth rate is a measurable function $b:\R^d\times\Omega\to [0, \infty)$ that
\begin{enumerate}[(i)]
    \item is \emph{uniformly elliptic}, i.e.\ $\norm{b}_{\infty},\norm{1/b}_{\infty}<+\infty$;
    \item is \emph{translation invariant}, i.e.\ $b(x,\eta) = b(o,\theta_x\eta)$ for all $x\in\R^d$ and $\eta\in\Omega$;
    \item has \emph{finite range}, i.e.\ there exists $R>0$ such that $b(x, \eta)= b(x, \eta_{B_R(x)})$ for all \(x \in \R^d\) and \(\eta \in \Omega\).
\end{enumerate}
\end{assumptions}
This means, for example, that $\GG_\theta$ is the set of translation-invariant Gibbs point processes associated to an energy functional whose birth rate satisfies \zcref{assumptions}\emph{(i)--(iii)}.
See~\zcref{sec:outlook} for a discussion on their optimality. We note that they guarantee that the set $\GG_\theta$ is non-empty by standard results, as they imply that $H$ is non-degenerate, hereditary, finite range, and stable (after the harmless normalization $H(\emptyset)=0$), and hence fits into the setting of~\cite[Theorem 1]{dereudre_2019}.

\begin{remark}
    The above assumptions are satisfied, for example, by the following two models, both known to exhibit a non-uniqueness phase transition.
    \begin{enumerate}[(i)]
        \item The \emph{area-interaction functional}
        \begin{equation*}
            H(\eta) = \alpha|B_{R/2}(\eta)| := \alpha\Big\lvert \bigcup_{x\in \eta}B_{R/2}(x) \Big\rvert, \quad \alpha\in\R\setminus\{0\}.  
        \end{equation*}
        The resulting birth rate has finite range $R$ and is uniformly elliptic, since $\abs{h(x,\eta)} = \abs{\alpha} |B_{R/2}(x)\setminus B_{R/2}(\eta)|\leq \abs{\alpha}\abs{B_R(o)}$. Moreover, translation invariance follows from the translation invariance of the Lebesgue measure.
        Note that the interaction is attractive (resp.\ repulsive),  if $\alpha$ is positive (resp.\ negative), and the birth rate is higher (resp.\ lower) closer to points of $\eta$. 
        It is a classical result see e.g.~\cite{Ruelle1971WRPT,Giacomin1995Agreement,CCK95}, that in dimensions \(d \geq 2\), the area interaction model at radius $1/2$ exhibits a non-uniqueness phase transition when the intensity of the underlying Poisson point process and the inverse-temperature parameter are equal and sufficiently large. In our parametrization, taking $\alpha=1$, the change of variables $x\mapsto x/R$ yields a model with radius $1/2$, Poisson intensity $R^d$, and inverse temperature $R^d$. Hence, non-uniqueness holds for sufficiently large $R$.
        \item The \emph{diluted pairwise interaction} introduced in~\cite{DRC26}, defined as follows: let $R>0$ and $\phi:\R_+\mapsto\R\cup\{+\infty\}$ be a compactly supported radial pair potential with $r^{d-1}\phi(r)\in L^1$, and consider
        \begin{equation*}
            H(\eta) = \int_{B_{R}(\eta)} \d x\int_{B_{R}(\eta)} \d y\, \phi(\abs{x-y}).
        \end{equation*}
        The resulting birth rate is clearly finite range and translation invariant. Moreover, it is uniformly elliptic, since
        \begin{equation*}
        \begin{split}
            \abs{h(x,\eta)} &=
            \abs{\int_{B_R(x)\setminus B_R(\eta)}\d y \Big(\int_{B_R(\eta)} \d z \, \phi(\abs{y-z}) + \int_{B_R(x) \cup B_R(\eta)}\d z \, \phi(\abs{y-z})\Big)}\\
            &\leq 2 d\abs{B_1(o)}^2 R^d \int_0^{+\infty} r^{d-1}\abs{\phi(r)}\, \d r.
        \end{split}
        \end{equation*}
        Note that in order to recover their high-intensity phase transition in our context, we need to rescale the potential, taking $\phi_R(r) = R^{-2d}\phi(r/R)$.
    \end{enumerate}
\end{remark}

\medbreak

Now, for a given configuration $\eta$, we would like to consider the infinite-volume \emph{birth-and-death process} $(X\Ssup\eta_t)_{t\ge 0}$ associated to the (formal) generator
\begin{equation*}\label{eq:genBD}
    (\LL f)(\gamma) := \int_{\R^d} \d x \, b(x,\gamma)\big( f(\gamma+\delta_x) - f(\gamma) \big) + \int \gamma(\d x) \big( f(\gamma-\delta_x) - f(\gamma)\big)
\end{equation*}
and with starting configuration $X_0\Ssup{\eta}=\eta$. 

\subsubsection{Graphical representation}

As the generator is based on the conditional energy $h$, one indeed expects the infinite-volume Gibbs measures associated to the energy functional $H$ to be reversible with respect to the dynamics.
A priori, however, because of the unbounded local density of particles, it is not clear that such a Markov process even exists. Following the classical theory of Markov processes, cf.~\cite[Theorem~4.4.2]{EK86}, this follows if we can show that the associated martingale problem is well posed. This is indeed the case: in~\cite{garcia2006spatial}, the authors show that any solution to the martingale problem for $\LL$ is also a solution to a given SDE, and vice versa. Moreover, this SDE, representing a thinning of a space-time Poisson random measure, has a unique solution, so that the martingale problem is well posed. 

\begin{figure}
    \centering
    \includegraphics[width=0.7\linewidth]{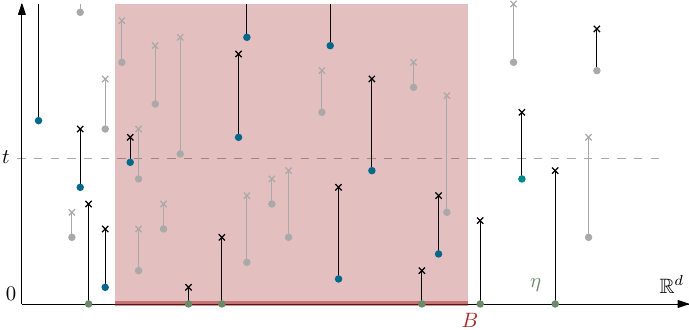}
    \caption{Illustration of the graphical representation of \zcref{Definition:Markov_Process}. At time $t=0$, an initial configuration $\eta$ (in green) is given; for positive times, the driving space-time Poisson measure $\Ncal$ proposes new points to be born at a given location and time (dots), with assigned lifetimes (length of the segment, where crosses represent death times); points are accepted (blue) or rejected (gray) depending on the current configuration. For any $B\Subset\R^d$ (red line) and $t>0$, $X_t\Ssup{\eta}(B)$ counts the number of initial and accepted points in $B$ that have survived until time $t$ (i.e.~all accepted points in the red shaded region whose lifeline intersects the dashed line).}
    \label{fig:graphical}
\end{figure}

We provide the details of the martingale-problem approach and the construction of the dynamics in~\zcref{sec:construction_dynamics}. Here, we only state the following proposition, which gives a graphical representation of the birth-and-death process as the unique (weak) solution of an SDE; note that such a process will take values in the space $E$ of c\`adl\`ag functions on counting measures on $\R^d$; see \zcref{fig:graphical} for an illustration.

\begin{proposition}[Graphical representation]\label{Definition:Markov_Process}
    There exist a probability space $(\Xi,\Gcal,\P)$, a filtration $(\Gcal_t)_{t\geq 0}$, and a Poisson random measure $\Ncal$ on $\R^d\times[0,\infty)^3$ with intensity measure $\d x\otimes \d s \otimes\e^{-r}\d r\otimes\d u$ such that, for all initial configurations $\eta = \sum_{i\ge 1}\delta_{x_i} \in \Omega$, there exists a pathwise-unique $(\Gcal_t)_{t\geq0}$-adapted process $X\Ssup{\eta}$ with c\`adl\`ag paths which $\P$-a.s.\ solves, for all $B \Subset \R^d$ and $t\geq 0$, the following SDE:
    \begin{equation}\label{eq:sdeBD}
    \begin{split}
        X\Ssup{\eta}_t(B) = \int_{B\times (0,t] \times [0,\infty)^2}\Ncal(\d x,\d s,\d r,\d u)\ &\1_{[0,b(x,X\Ssup{\eta}_{s-})]}(u)\1_{(t-s,\infty)}(r) \\
        &+ \int_{B\times [0,\infty)}\hat{\eta}(\d x,\d r)\ \1_{(t,\infty)}(r),
    \end{split}
    \end{equation}
    with $X\Ssup{\eta}_0(B) = \eta(B)$, where $\hat\eta = \sum_{i\ge 1}\delta_{(x_i,\tau_i)}$ is the point process on $\R^d\times[0,\infty)$ obtained by associating to each $x_i$ a mark $\tau_i$, and \((\tau_i)_i\) are i.i.d.\ unit-exponential random variables independent of \(\Ncal\).
\end{proposition}

\medskip 
We denote the {\em transition semigroup} of the Markov process $X\Ssup{\eta}$ by
\begin{align*}
    (T_t f)(\eta):= \mathbb{E}[f(X_t\Ssup{\eta})],\quad t\geq 0,
\end{align*}
for all $f$ for which the right-hand side is well-defined, 
and note that, as shown in \zcref{lemma:infinite_volume_process_forward_equation_local_bounded_measurable} below,
\begin{align*}
    \frac{\d}{\d s}T_s f\Big\vert_{s = t} = T_t \LL f, \quad f\in \Lcal,
\end{align*} 
where $\Lcal$ denotes the set of all bounded, local, measurable functions $f\colon\Omega \to \R$.
Moreover, let us mention that, due to the finite speed of propagation exhibited in \zcref{Lemma:finite_speed_of_propagation_version_1}, the existence of the infinite-volume dynamics can also be proven by other means. 

\subsubsection{Finite-volume dynamics}\label{section:local-dynamics}

Since the projection of the dynamics to a finite volume $\Lambda$ is not Markovian due to interactions with $\Lambda^c$, we introduce modified dynamics \(( \Tl_t )_{t \geq 0}\), in which the birth rate incorporates random boundary conditions sampled from a translation-invariant Gibbs measure \(\nu\in\GG_\theta\). We later compare their behavior to the original dynamics, see~\zcref{sec:relative-entropy-decreasing}, and use them to show the monotonicity of the specific relative entropy, see~\zcref{sec:monotonicity_sketch}. Furthermore, these dynamics are reversible w.r.t.\ the marginal \(\nu_\Lambda\) and have a strongly continuous semigroup on \(L^1(\nu_\Lambda)\), thus allowing us to use semigroup theory to compute the integral identities needed for entropy dissipation in~\zcref{sec:attractor-property}.

We construct a Markov process \((X^{(\Lambda, \eta)}_t)_{t \geq 0}\) with transition semigroup \((\Tl_t)_{t \geq 0}\), which will be given in terms of an SDE similar to the one above, associated to the formal generator
\begin{align*}
    (\mathscr{L}_\Lambda f)(\eta_\Lambda)
    = \int_{\Lambda} \d x\, b^\nu_\Lambda(x, \eta_\Lambda) \big[f(\eta_\Lambda+\delta_x) - f(\eta_\Lambda)\big] + \sum_{x\in \eta_\Lambda}\, \big[f(\eta_\Lambda -\delta_x) - f(\eta_\Lambda)\big],
\end{align*} 
where the \emph{localized} birth rate \(b_\Lambda^\nu\) is defined by
\begin{align*}
    b_\Lambda^\nu(x, \eta_\Lambda)
    := \int \nu_{\Lambda \,\vert\, \Lambda^c}(\d\zeta_{\Lambda^c} \vert \eta_\Lambda)\, b(x, \eta_\Lambda \zeta_{\Lambda^c}), 
    \quad x \in \Lambda, \eta_\Lambda \in \Omega_\Lambda,
\end{align*}
with \(\nu_{\Lambda \,\vert\, \Lambda^c}\) a regular version of the conditional distribution under \(\nu\in \GG_\theta\) of the particles in \(\Lambda^c\) given the configuration inside \(\Lambda\). See~\zcref{section:finite_volume_dynamics} for more details.


\subsection{Reversible measures}\label{sec:intro_reversible}

Our main objective is to understand the birth-and-death dynamics when started out of equilibrium, from another translation-invariant measure.

In order to identify all weak limit points of the dynamics as Gibbs measures, we start out by providing a more convenient characterization of reversible measures and Gibbs measures. For this, we introduce the following concept of reversibility with respect to the semigroup $(T_t)_{t\geq 0}$:
\begin{definition}
    We say that \(\mu \in \Pcal(\Omega)\) is {\em reversible} if, for every \(t \geq 0\) and every \(f, g \in \Lcal\), 
    \begin{equation*}
        \mu[(T_t f) \, g] = \mu[f \, (T_t g)].  
    \end{equation*} 
    The set of translation-invariant reversible probability measures is denoted by \(\RR_\theta\).
\end{definition}
The expected equivalence between (translation-invariant) reversible measures and infinite-volume Gibbs point processes holds, and is proven, in~\zcref{sec:reversible-Gibbs}, with the help of the GNZ equations and Palm measures, which we now introduce.

If \(\mu \in \Pcal_\theta(\Omega)\), there exists a unique measure $\mu_o$ on \(\Omega\), its {\em (unnormalized) Palm measure}, such that
\begin{align*}
    \int\mu(\d \eta)\int \eta(\d x) f(x,\theta_x\eta)
    = \int \d x\ \mu_o[f(x,\cdot)],\qquad f\colon \R^d\times\Omega\to[0,\infty) \text{ measurable},
\end{align*} 
describing the point process from a typical point, shifted to the origin. Removing such a point leads to the {\em (unnormalized) reduced Palm measure} $\mu^!_o$, characterized by
\begin{align*}
    \int\mu(\d \eta)\int \eta(\d x) f(x,\theta_x\eta-\delta_o)
    = \int \d x\ \mu_o^![f(x,\cdot)],\qquad f\colon \R^d\times\Omega\to[0,\infty) \text{ measurable}.
\end{align*}
See for example~\cite{Kallenberg1983Random}. We note that the unnormalized reduced Palm measure \(\mu_o^!\) is always \(\sigma\)-finite; moreover, in the case that \(\mu \in \Pcal_{1, \theta}(\Omega)\), i.e.\ if \(\mu\) additionally has finite intensity, it is actually a finite measure with \(\mu_o^![1] = \mu[N_{[0, 1)^d}]\).

\begin{proposition}\label{Prop:Reversible_Measures_equal_Gibbs_Measures}
    The following statements are equivalent for \(\mu \in \Pcal_\theta(\Omega)\).
    \begin{enumerate}[(i)]
        \item \(\mu \in \RR_\theta\).
        \item \(\mu\) is \emph{infinitesimally reversible}, i.e.\
            \begin{equation}\label{equation:definition_infinitesimally_reversible}
                \mu\big[(\mathscr{L} f) \, g\big] = \mu\big[f \, (\mathscr{L} g)\big]
            \end{equation} for all \(f, g \in \Lcal\) such that \(\mathscr{L} f, \mathscr{L} g \in \Lcal\).
        \item \(\mu\) is such that
        \begin{align}
            b(o, \cdot) \mu =
            \mu_o^!.
        \end{align} 
        \item \(\mu \in \GG_\theta\).
    \end{enumerate}
    In particular, \(\RR_\theta = \GG_\theta\).
\end{proposition}

\subsection{Entropy and Fisher information}

Before stating our main results, we need to introduce two key functionals: the specific relative entropy and the modified Fisher information. 

\begin{definition}
    Let $\mu,\nu$ be two finite measures on $\Omega$. The {\em local relative entropy} of $\mu$ w.r.t.\ $\nu$ on $\L\Subset\R^d$ is defined as
\begin{align*}
    \RelEnt_\L(\mu|\nu):=\RelEnt(\mu_\L|\nu_\L) := \int\mu(\d\eta)\log\frac{\d \mu_\L}{\d\nu_\L}(\eta_\L),
\end{align*} 
if the density $\d \mu_\L/\d\nu_\L$ exists, and $+\infty$ otherwise. Similarly,
\begin{equation*}
    \RelEnt(\mu | \nu) := \int\mu(\d\eta)\log\frac{\d \mu}{\d\nu}(\eta),
\end{equation*} 
if the density $\d \mu/\d\nu$ exists, and $+\infty$ otherwise.

The {\em specific relative entropy} of $\mu$ w.r.t.\ $\nu$ is defined as
\begin{align*}
    \SpecEnt(\mu|\nu):=\liminf_{n\uparrow\infty}\frac{\RelEnt_{\L_n}(\mu|\nu)}{\abs{\L_n}}.
\end{align*}
\end{definition}
If $\nu\in\GG_\theta$ and $\mu\in\Pcal(\Omega)$, the specific relative entropy $\SpecEnt(\mu|\nu)$ is also known as \emph{excess free energy}, and it exists as an actual limit if $\mu$ is translation-invariant, see \zcref{prop:existence-relative-entropy} below. 

\begin{definition}
   Let $\nu \in \GG_\theta$ and $\mu_\L\in\Pcal(\Omega_\Lambda)$, $\L\Subset\R^d$, such that $\mu_\L \ll \nu_\L$. The \emph{modified Fisher information} (cf.~\cite{Dello2024Wasserstein,huesmann2025} for the non-interacting case) is defined as
\begin{align}\label{def:fisher-information}
    \mathcal{J}_\L(\mu_\L\lvert \nu) := \int_\Lambda \d x\int \nu_\L(\d\eta)\  b^{\nu}_\L(x,\eta_\L)\ D_x\frac{\d\mu_\L}{\d\nu_\L}(\eta_\L) \ D_x\log\frac{\d\mu_\L}{\d\nu_\L}(\eta_\L),
\end{align}
where
\begin{equation*}
    D_x F(\eta) := F(\eta+\delta_x)-F(\eta)
\end{equation*}
is the {\em add-one-cost operator}. In the case that \(\mu_\Lambda \not\ll \nu_\Lambda\), we simply set \(\mathcal{J}_\L(\mu_\L\lvert \nu) := +\infty\).
\end{definition} 
We use here the convention that, for \(a, b \geq 0\), \((b-a)(\log b - \log a)=0\) if \(a=b=0\), while \((b-a)(\log b - \log a)=+\infty\) if exactly one of \(a, b\) is \(0\). 

Finally, for \(\mu \in \Pcal_{1, \theta}(\Omega)\), let
\begin{align}\label{eqn:definition-sigma}
    \sigma(\mu)
    := \RelEnt\big(b(o, \cdot)\mu \,\big\vert\, \mu_o^! \big) +  \RelEnt\big(\mu_o^! \,\big\vert\, b(o, \cdot)\mu \big)
\end{align}
and let $\sigma(\mu):=+\infty$ if $\mu\in \Pcal_{\theta}(\Omega)$ has infinite intensity. Note that $\sigma$ is always non-negative, see~\zcref{rmk:sigmanonnegative}.

\subsection{Main results}\label{sec:results}
In order to analyze the long-time behavior of the dynamics for arbitrary initial distributions, we rely on a Lyapunov-type approach based on the specific relative entropy. The advantage of this strategy is that it also works in the phase-transition regime and for initial distributions that are not absolutely continuous with respect to the reversible measures. 

We can now present our first main result, which establishes the entropy dissipation, connecting entropy production as time progresses and Fisher information. 

\begin{theorem}\label{Theorem:Main_Theorem_Decrease_in_Entropy}
Assume the birth rate $b$ satisfies \zcref{assumptions} and let \(\nu\in \GG_\theta\). 
\begin{enumerate}[(i)]
    \item For any starting measure \(\mu \in \Pcal_\theta(\Omega)\), the map \(t \mapsto \SpecEnt(\mu T_t | \nu)\) is non-increasing.
    \item If \(\mu \in \Pcal_\theta(\Omega)\) satisfies \(\SpecEnt(\mu | \nu) < \infty\), then
    \begin{align*}
    \SpecEnt(\mu | \nu) - \SpecEnt(\mu T_t | \nu) &\geq \int_0^t \d s\ \liminf_{n\uparrow\infty} \frac{\mathcal{J}_{\L_n}(\mu T\Ssup{\L_n}_s\vert\nu)}{\abs{\L_n} }\geq \int_0^t\d s\ \sigma(\mu T_s) \ge 0.
\end{align*}
\item For \(\mu \in \Pcal_\theta(\Omega)\) we have \(\sigma(\mu) = 0\) if and only if \(\mu \in \GG_\theta\).
\end{enumerate}
\end{theorem}
Our main take away from \zcref{Theorem:Main_Theorem_Decrease_in_Entropy} is the following attractor property. We work here with the weak topology on $\Pcal(\Omega)$ induced by the vague topology on $\Omega$, i.e.\ the coarsest topology on $\Pcal(\Omega)$ that makes $\mu\mapsto \mu[f]$ continuous for every bounded continuous function $f:\Omega\to\R$.
\begin{theorem}\label{Theorem:Main_Theorem_Limit_Points}
    Assume the birth rate $b$ satisfies \zcref{assumptions}, let $\nu\in\GG_\theta$ and \(\mu \in \Pcal_\theta(\Omega)\) with \(\SpecEnt(\mu | \nu) < \infty\). Then, for any $\mu^*$ such that there exists an increasing sequence of times $(t_k)_{k\ge 0}$, with $t_k\uparrow\infty$ and \(\lim_{k\uparrow\infty}\mu T_{t_k}=\mu^*\) in the weak topology, we have that $\mu^*\in\GG_\theta$.
    In particular, \(\mu^*\) is reversible.
\end{theorem}

As a corollary of the attractor property, we obtain the following characterization of the time-stationary measures for the dynamics generated by $\mathscr{L}$. In our context, a measure \(\mu \in \Pcal(\Omega)\) is called \emph{time-stationary} if \(\mu T_t[f] = \mu[f]\) holds for all \(f \in \Lcal\) and all \(t \geq 0\).
Let
\begin{align*}
    \mathcal{S}_\theta
    := \{\mu \in \Pcal_\theta \,\colon\, \mu \text{ is time-stationary}\}.
\end{align*}
In this notation, the attractor property implies the following. 
\begin{corollary}\label{corollary:time-stationary-iff-gibbs}
    Assume that the birth rate $b$ satisfies \zcref{assumptions}. We then have \(\mathcal{S}_\theta = \GG_\theta\).
\end{corollary}
Since time-stationarity is weaker than reversibility, the inclusion \(\GG_\theta = \RR_\theta \subseteq \mathcal{S}_\theta\) already follows from \zcref{Prop:Reversible_Measures_equal_Gibbs_Measures}. 
For the inclusion \(\mathcal{S}_\theta \subseteq \GG_\theta\), we notice that \zcref{Theorem:Main_Theorem_Limit_Points} directly yields that every time-stationary \(\mu \in \Pcal_\theta(\Omega)\) with \(\SpecEnt(\mu \,\vert\, \nu) < \infty\) has to be Gibbsian too, i.e., \(\mu \in \GG_\theta\). Hence, it suffices to show that a time-stationary \(\mu \in \Pcal_\theta(\Omega)\) automatically has finite specific relative entropy. This is done in~\zcref{sec:time-stationary}. 
\medskip

Our main results confirm the physical intuition that the excess free energy is non-increasing in time and moreover governs the long-term behavior of large systems of interacting particles. The proofs of our two main results are inspired by the method developed in \cite{Sullivan1976} for interacting particle systems on the lattice; we explain the general strategy in~\zcref{section:proof-strategy} before providing all of the proofs in the subsequent sections. 

\subsection{Generalizations and outlook}\label{sec:outlook}

Before we turn toward the proofs of our main results, let us briefly discuss some possible generalizations and comment on future research directions.

\subsubsection{Model assumptions}

It is natural to ask under what weaker assumptions on the interaction functional our main results remain valid. The uniform upper bound on the birth rate and the finite-range assumption are e.g.\ used in the construction of the dynamics and the finite speed of propagation result~\zcref{Lemma:comparison_finite_infinite}, whereas the uniform positive lower bound also enters the factor property of~\zcref{Le:factor_property_of_the_Gibbs_measure}.
One possible generalization would allow bounded birth rates that are not bounded away from zero, as for example arise from purely repulsive pair interactions.
More generally, one could investigate birth rates with some sort of temperedness in local particle counts, using the full framework of~\cite{garcia2006spatial} as a starting point for the construction. An extension of the entropy arguments used here would then require suitable moment estimates and control of the boundary terms in the finite-volume comparison.
Translation invariance, on the other hand, ensures that the dynamics preserve the translation invariance of the starting measure, which is crucial for our arguments based on Palm measures. Extending this approach to spatially inhomogeneous birth rates would require additional ideas.

\subsubsection{Outlook}

Can the entropy-dissipation inequality in~\zcref{Theorem:Main_Theorem_Decrease_in_Entropy} be strengthened to an equality, and can the dissipation be expressed as a thermodynamic limit of finite-volume modified Fisher information, thus yielding a de Bruijn-type identity? In the follow-up work~\cite{JKSZ26}, we answer this question affirmatively for continuous interactions and for a large class of starting measures. The proof follows a different approach, since obtaining such an identity with the strategy of the current paper presents two major hurdles.
First, while one would intuitively believe that the converse upper bound in \zcref{Le:compariability_of_finite_volume_entropies_original_vs_modified} should also hold, the structure of the Donsker--Varadhan variational formula makes it much harder to prove. While for the lower bound we only needed to use a single appropriately chosen test function, the upper bound requires bounding the functional from above for arbitrary test functions. Second, while for the Poisson case the density limit of the Fisher information is shown to exist with a subadditivity argument, see \cite[Lemma 7.8]{huesmann2025}, in the interacting case the situation is not as clear, as the measure is not exactly factorizable and one loses the precise superadditivity. In contrast to the proof of \zcref{Le:quasi_superadditivity_entropy_wrt_gibbs_measure}, the factor property of $\nu$ does not seem to yield an $o(\abs{\Lambda})$ error term.

Some other possible directions that are inspired by analogous results for interacting particle systems on the lattice are as follows. First, one could see how well our results can be extended to the case of non-reversible birth-and-death dynamics as for example in \cite{kunsch_time_1984,jahnel_dynamical_2023}. 
Second, it would be interesting to see if one can also apply a similar strategy based on the change of relative entropy to analyze the long-time behavior without translation invariance. For reversible lattice systems, this can be done in one and two dimensions, see \cite{holley_one_1977} and \cite{JK25}, and it is not clear if these  arguments can be extended to the continuum.

\section{Proof strategy}\label{section:proof-strategy}

For this section, we will take the existence and further properties of the dynamics as granted. We take care of these technically important details in~\zcref{sec:construction_dynamics}. Moreover, from now on, we will always assume that all properties required in \zcref{assumptions} are satisfied but will not explicitly mention it in the statements of all of the lemmas and propositions below. 

\subsection{Monotonicity of the specific relative entropy}\label{sec:monotonicity_sketch}
As already explained in~\zcref{sec:finite-state-space-example}, the monotonicity of the map $t \mapsto \RelEnt(\nu P_t \lvert \mu)$, where $(P_t)_{t\geq 0}$ is a Markov semigroup with invariant measure $\mu$ and $\nu$ is some arbitrary initial distribution, is a direct consequence of Jensen's inequality. However, since the projection of the infinite volume dynamics to a finite volume $\Lambda$ is not Markovian due to interactions with $\Lambda^c$, the finite volume entropy functionals $t \mapsto \RelEnt_\Lambda(\mu T_t \lvert \nu)$ are in general not monotone. For this reason, we use the local dynamics $(T^{(\Lambda)}_t)_{t\geq0}$, introduced in~\zcref{section:local-dynamics} as an auxiliary process.  For fixed $\Lambda \Subset \R^d$, we can then use Jensen's inequality to see that $t\mapsto \RelEnt_\Lambda(\mu T^{(\Lambda)}_t\lvert \nu)$ is non-increasing. To use this finite-volume information to say something about the infinite-volume dynamics, we need to bound the error we made by this local approximation. As a first step, define the specific relative entropy functional associated to the local dynamics via
\begin{align*}
    \widetilde\SpecEnt_{\mu}(t)
        := \liminf_{n \uparrow \infty} n^{-d}\RelEnt_{\L_n}\bigl(\mu T\Ssup{\L_n}_t |\nu\bigr).
\end{align*}
As a pointwise liminf of non-increasing functions, $\widetilde\SpecEnt_{\mu}(\cdot)$ is also non-increasing. 
The main goal of this subsection is to prove the first part of our main result \zcref{Theorem:Main_Theorem_Decrease_in_Entropy}, which is implied by the following comparison result. 
\begin{proposition}\label{Prop:specific_entropy_is_decreasing_in_time}
        Let \(\mu \in \Pcal_\theta(\Omega)\).
        We have that $\SpecEnt(\mu T_t | \nu)\leq \widetilde \SpecEnt_{\mu}(t)$.
        In particular, \(t \mapsto \SpecEnt(\mu T_t | \nu)\) is non-increasing.
\end{proposition}

The proof of this comparison rests on two technical ingredients. First, we show that the finite-volume relative entropy is quasi-superadditive, i.e., superadditive up to sub-volume-order corrections. 

\begin{lemma}[Quasi-superadditivity of the relative entropy w.r.t.\ Gibbs measures]\label{Le:quasi_superadditivity_entropy_wrt_gibbs_measure}
     Let \(\mu \in \Pcal_1(\Omega)\), $\nu \in \mathscr{G}_\theta$, and let \(\Lambda, \Delta \Subset \R^d\) be disjoint. Then,
    \begin{align*}
        \RelEnt_{\Lambda \cup \Delta}(\mu |\nu)
        &\geq \RelEnt_{\Lambda}(\mu |\nu) + \RelEnt_{\Delta}(\mu |\nu) \\
        &\qquad - \big\{ \log(\Vert b\Vert_\infty \Vert 1/b\Vert_\infty) \mu[N_{\partial^R \Lambda}] + (\Vert b\Vert_\infty - 1/\Vert 1/b\Vert_\infty) \abs{\partial^R \Lambda} \big\}. \nonumber
    \end{align*} 
\end{lemma}

The second ingredient is the following approximation property that essentially tells us that in order to approximate the infinite-volume process in a finite observation window $\Lambda \Subset \R^d$ up until time $T$, it suffices to look at the finite-volume dynamics in a linearly enlarged window $\Lambda \oplus B(0,cT)$ for some constant $c>0$. 
To make sense of this mathematically, we make the natural identification between \(\Lambda\)-local measurable functions on \(\Omega\) and measurable functions on \(\Omega_\Lambda\).
\begin{lemma}[Finite speed of propagation]\label{Lemma:finite_speed_of_propagation_version_1}
    For any box $\L\Subset\R^d$ with each side length at least \(1\), let $\widehat{\L} \Subset \R^d$ be some superset of \(\Lambda\), with \(\ell := \mathrm{dist}(\Lambda, \widehat{\Lambda}^{c}) > 0\). Let \(\mu\in\Pcal_\theta(\Omega)\).
    Then, for any \(\Lambda\)-local \(f\in \Lcal\) and every bounded measurable \(g\), we have
    \begin{align}
        \big\vert \mu[g T_t f] - \mu[g T_t\Ssup{\widehat{\Lambda}}f] \big\vert
        \lesssim_t
        \vert \Lambda\vert \Vert f \Vert_\infty \Vert g \Vert_\infty  \e^{-\mathrm{dist}(\Lambda, \widehat{\Lambda}^{c}) \rho(\mathrm{dist}(\Lambda, \widehat{\Lambda}^{c}))}
    \end{align} 
    for some increasing function \(\rho:\R_+ \to \R_+\) (depending on \(t, d, \Vert b \Vert_\infty, R\)) with \(\rho(\ell) \xrightarrow[\ell\uparrow \infty]{} \infty\).
    
    Moreover, if \(\mu \in \Pcal_{1, \theta}(\Omega)\), i.e.\ if additionally \(\mu[N_{[0, 1)^d}] < \infty\), the following holds.
    Consider a family $\{f_\L\sepset  \L\Subset\R^d \text{ box}\}$ of measurable functions such that \(f_\Lambda\) is \(\Lambda\)-local and \(\vert f_\Lambda \vert \lesssim N_{\Lambda}\) (uniformly in $\Lambda$). If \(\Lambda, \widehat{\Lambda} \uparrow \mathbb{R}^d\) in such a way that \(\log(\vert\Lambda\vert)/\mathrm{dist}(\Lambda, \widehat{\Lambda}^{c})\) stays bounded from above, then, as $\L\uparrow\R^d$,
    \begin{align*}
        \frac{1}{\vert \Lambda\vert}\big\lvert \mu T_t[f_\Lambda] - \mu T_t\Ssup{\widehat{\Lambda}}[f_\Lambda] \big\rvert
        \xrightarrow[]{} 0.
    \end{align*} 
\end{lemma}

Using this as input, we can now show that the change in relative entropy under the global dynamics can be lower bounded by the change with respect to the local dynamics. 
For this, we have to compare finite-volume entropies in enlarged volumes. Consider a set $\L$, contained in a larger set $\L^\ast$; the lemma below shows that the finite-volume entropy \(\RelEnt_{\Lambda}(\mu T_t | \nu)\) is at most only slightly larger than \(\RelEnt_{\Lambda}\bigl(\mu T_t\Ssup{\Lambda^\ast} \,\big\vert\, \nu \bigr)\) if \(\Lambda\) is sufficiently far from the boundary \(\partial \Lambda^\ast\).

\begin{lemma}\label{Le:compariability_of_finite_volume_entropies_original_vs_modified}
            Let \(\mu \in \Pcal_\theta(\Omega)\) and \(\nu \in \Gcal_\theta\).
            Fix \(t > 0\) and \(\Lambda \Subset \mathbb{R}^d\).
            For every \(r\) with
            \begin{align*}
                \RelEnt_{\Lambda}(\mu T_t | \nu) > r,
            \end{align*} there exists \(\ell > 0\), such that, for every \(\Lambda^* \Subset \R^d\) with \(\Lambda \subseteq \Lambda^*\) and for every translate \(\widetilde{\Lambda} \subseteq \Lambda^*\) of \(\Lambda\) with \(\mathrm{dist}(\widetilde{\Lambda}, (\Lambda^*)^c) \geq \ell\), the following holds:
            \begin{align}
                 \RelEnt_{\widetilde{\Lambda}}\bigl(\mu T_t\Ssup{\Lambda^\ast} \,\big\vert\, \nu \bigr) 
                 > r.
            \end{align}
        \end{lemma}

To control the error term in~\zcref{Le:quasi_superadditivity_entropy_wrt_gibbs_measure} and use the second part of~\zcref{Lemma:finite_speed_of_propagation_version_1}, we need to have  some control over the expected number of points in bounded volumes. This is the content of the following technical lemma that follows more or less directly from the graphical construction of the finite-volume dynamics which is detailed in~\zcref{sec:appendix:finite_volume_dynamics}. 

\begin{lemma}\label{lemma:control-number-of-points}
    Let $\mu \in \mathcal{P}_\theta(\Omega)$. Then for every $t\geq 0$, $\Lambda \Subset \R^d$,  and measurable $B\subset \Lambda$ we have 
    \begin{equation}\label{equation:bounds_particle_numbers_local_dynamics}
        e^{-t}\mu[N_B] \leq (\mu T^{(\Lambda)}_t)[N_B] \leq e^{-t}\mu[N_B] + \norm{b}_\infty(1-e^{-t})\abs{B}. 
    \end{equation}
    In particular, if there is a $k\in \N$ such that $I_{\Lambda_k}(\mu T^{(\Lambda_k)}_t \lvert \nu) < \infty$, then also $\mu[N_{[0,1)^d}]<\infty$. 
\end{lemma}
The lower bound in \eqref{equation:bounds_particle_numbers_local_dynamics} is essentially the contribution from the initial points that survived until time $t$, while the upper bound additionally gets a contribution from the points that were born after time $0$ and haven't yet died. \medskip

With all of these intermediate results at our hands, we can provide the proof of~\zcref{Prop:specific_entropy_is_decreasing_in_time}, and thereby~\zcref{Theorem:Main_Theorem_Decrease_in_Entropy}\emph{(i)}. The finite speed of propagation estimate in \zcref{Lemma:finite_speed_of_propagation_version_1} is proved in \zcref{section:finite_volume_dynamics}, while the proofs of the remaining technical lemmas are deferred to Section~\ref{sec:relative-entropy-decreasing}. 

\begin{proof}[Proof of~\zcref{Prop:specific_entropy_is_decreasing_in_time}]
    For every finite volume \(\Lambda\), the function \(t \mapsto \RelEnt_{\Lambda}\bigl(\mu \Tl_t \big|\nu\bigr)\) is non-increasing by Jensen's inequality. Hence, \(\widetilde{\SpecEnt}_\mu(t)\) is a non-increasing function of \(t\) as a limit of non-increasing functions.
            It follows that we only have to show for fixed \(t > 0\) that
            \begin{align}\label{equation:order_in_modified_entropy}
                \SpecEnt(\mu T_t | \nu)
                \leq \widetilde{\SpecEnt}_{\mu}(t), \quad \mu \in \Pcal_\theta(\Omega),
            \end{align} 
            in order to see that \(t \mapsto \SpecEnt(\mu T_t | \nu)\) is non-increasing, as then by the semigroup property
            \begin{align*}
                \SpecEnt(\mu T_{s+t} | \nu)
                = \SpecEnt((\mu T_{s}) T_t | \nu)
                \leq \widetilde{\SpecEnt}_{\mu T_s}(t)
                \leq \widetilde{\SpecEnt}_{\mu T_s}(0)
                = \SpecEnt(\mu T_{s} | \nu).
            \end{align*}
            If \(\widetilde{\SpecEnt}_{\mu}(t) = \infty\) then \eqref{equation:order_in_modified_entropy} is trivially true and we can therefore assume without loss of generality that there is a $k\in \N$ such that $I_{\Lambda_k}(\mu T_t^{(\Lambda_k)}\lvert \nu) < \infty$. Thus, by~\zcref{lemma:control-number-of-points} we must have $\mu[N_{[0,1)^d}]< \infty$.
            By \eqref{equation:bounds_particle_numbers_local_dynamics} and translation invariance, we get
            \begin{align}\label{equation:bounds_particle_numbers_local_dynamics_2}
                (\mu \Tl_t)[N_B]
                \leq (\Vert b\Vert_\infty + \mu[N_{[0, 1)^d}])\abs{B}.
            \end{align} 
            Take \(r < \RelEnt_{\Lambda_n}(\mu T_t \vert \nu)\) and fix \(n \in \N\).
            We now show \begin{equation}\label{equation:intermediate_lower_bound_relative_entropy_local_dynamics}
            \begin{split}
                &\widetilde{\SpecEnt}_{\mu}(t)
                = \liminf_{k \uparrow \infty} \frac{\RelEnt_{\L_k}\bigl(\mu T\Ssup{\L_k}_t |\nu\bigr)}{\abs{\Lambda_k}} \\
                &\geq \frac{r}{\abs{\Lambda_n}}  
                - \frac{\abs{\partial^R \Lambda_n}}{\abs{\Lambda_n}}\big\{ \log(\Vert b\Vert_\infty \Vert 1/b\Vert_\infty) (\mu[N_{[0,1)^d}] + \Vert b \Vert_\infty) + (\Vert b\Vert_\infty - 1/\Vert 1/b\Vert_\infty)  \big\},
            \end{split}
            \end{equation} 
            so that the proof is finished by first letting \(r \uparrow \RelEnt_{\Lambda_n}(\mu T_t \vert \nu)\) and then \(n \to \infty\). By \zcref{Le:compariability_of_finite_volume_entropies_original_vs_modified}, there exists a buffer distance \(\ell > 0\) such that
            \begin{align}\label{equation:lower_bound_r_relative_entropy_local_dynamics}
                \RelEnt_{\widetilde{\Lambda}}(\mu T_t\Ssup{\Lambda_k} |\nu) > r
            \end{align} for every \(k \geq n\) and every translate \(\widetilde{\Lambda} \subseteq \Lambda_k\) of \(\Lambda_n\) satisfying \(\mathrm{dist}(\widetilde{\Lambda}, (\Lambda_k)^c) \geq \ell\).
             Let \(N_k \in \N\) be the maximal number of disjoint translates \(\widetilde{\Lambda} \subseteq \Lambda_k\) of \(\Lambda_n\) with \(\mathrm{dist}(\widetilde{\Lambda}, (\Lambda_k)^c) \geq \ell\), and let \((\Lambda^{(k)}_i)_{1 \leq i \leq N_k}\) be a family of such disjoint translates. It then holds \(\frac{N_k \abs{\Lambda_n}}{\abs{\Lambda_k}} \xrightarrow[k \to \infty]{} 1\).  By repeatedly using \zcref{Le:quasi_superadditivity_entropy_wrt_gibbs_measure} and then \eqref{equation:lower_bound_r_relative_entropy_local_dynamics}, we get
            \begin{equation}
            \begin{split}
                \RelEnt_{\Lambda_k}(\mu T_t\Ssup{\Lambda_k} |\nu) \geq \sum_{i = 1}^{N_k} \RelEnt_{\Lambda_i^{(k)}}(\mu T_t\Ssup{\Lambda_k} |\nu) 
                - \log(\Vert b\Vert_\infty \Vert 1/b \Vert_\infty) &\sum_{i = 1}^{N_k} (\mu T_t\Ssup{\Lambda_k})[N_{\partial^R \Lambda_i^{(k)}}]
                \\
                - (\Vert b\Vert_\infty - 1/\Vert 1/b\Vert_\infty) &\sum_{i = 1}^{N_k} \abs{\partial^R \Lambda_i^{(k)}}
            \end{split}
            \end{equation} 
            for \(k\) large enough so that \(N_k \geq 1\). 
            Using \eqref{equation:bounds_particle_numbers_local_dynamics_2} and the choice of \(\ell\) we therefore obtain 
            \begin{align*}
                &\RelEnt_{\Lambda_k}(\mu T_t\Ssup{\Lambda_k} |\nu) \\
                &\geq N_k r
                - N_k \abs{\partial^R \Lambda_n}\big\{ \log(\Vert b\Vert_\infty \Vert 1/b\Vert_\infty) (\mu[N_{[0,1)^d}] + \Vert b \Vert_\infty) + (\Vert b\Vert_\infty - 1/\Vert 1/b\Vert_\infty)  \big\} \nonumber
            \end{align*} 
            Dividing by \(\abs{\Lambda_k}\) and taking the \(\liminf_{k \to \infty}\) yields \eqref{equation:intermediate_lower_bound_relative_entropy_local_dynamics} and therefore the statement of this proposition.
\end{proof}

\subsection{The Fisher information and a lower bound}
Now we at least know that the specific relative entropy is non-increasing along the measure-valued trajectories $t \mapsto \mu T_t$ for $\mu \in \mathcal{P}_\theta(\Omega)$, but in order to say something about limit points, we have to obtain a nicer representation of the functional $\mu \mapsto \mathscr{I}(\mu \lvert \nu) - \mathscr{I}(\mu T_t \lvert \nu)$ that measures the dissipation of the relative entropy between time $0$ and time $t$. 
Our first goal for this section is to derive the rather explicit representation 
\begin{align*}
    \RelEnt_{\Lambda}(\mu |\nu) - \RelEnt_{\Lambda}(\mu \Tl_t |\nu)
    = \int_{0}^{t}\d s\, \mathcal{J}_\L(\mu \Tl_s\vert\nu),
\end{align*} 
for the change of the relative entropy under the local dynamics. This identity is the starting point for the proof of the second statement in \zcref{Theorem:Main_Theorem_Decrease_in_Entropy}.  \medskip 

First, we establish the above finite-volume entropy-dissipation identity in terms of the modified Fisher information. We then express the latter as the (spatially averaged) relative entropy of some particular measures on a double layer system, i.e., with configuration space $\Omega_\Lambda \times \Omega_\Lambda$.
These measures are related to the  addition and removal of a point and are defined as follows. The required background on local Palm distributions can be found in \cite{DaleyVereJones2008}, see in particular ~\cite[Equations 13.1.4 and 13.1.5]{DaleyVereJones2008}. 
    \begin{definition}\label{def:Palm}
        Let \(\mu_\Lambda \in \Pcal(\Omega_\Lambda)\) with finite intensity and intensity measure absolutely continuous to Lebesgue measure on \(\Lambda\).
        Let \((\mu_\Lambda^x)_{x \in \Lambda}\) be the (unnormalized) local Palm measures for \(\mu_\Lambda\), i.e.,
        \begin{align*}
            \int \mu_{\Lambda}(\d\eta_\Lambda)\int \eta_{\Lambda}(\d x) f(x, \eta_\Lambda)
            = \int_{\Lambda}\d x \int \mu_\Lambda^x(\d\eta_\Lambda) \, f(x, \eta_\Lambda),
        \end{align*} 
        for suitable functions \(f\) (non-negative or integrable with respect to the Campbell measure).
        
        We introduce the measures \(\mu_\Lambda \ast G\Ssup{\Lambda}_x\) and \(\mu_\Lambda \circ G\Ssup{\Lambda}_x\) on \(\Omega_\Lambda \times \Omega_\Lambda\) defined, respectively, via
        \begin{align*}
            &\int (\mu_\Lambda \ast G\Ssup{\Lambda}_x)(\d\eta_\Lambda, \d\zeta_\Lambda) \, f(\eta_\Lambda, \zeta_\Lambda)\\
            &= \int \mu_\Lambda(\d \eta_\Lambda) \, b^\nu_\Lambda(x, \eta_\Lambda) \, f(\eta_\Lambda, \eta_{\Lambda} +\delta_{x})
            + \int \mu_\Lambda^{x}(\d \eta_{\Lambda}) \, f(\eta_\Lambda, \eta_\Lambda -\delta_x) 
        \end{align*}
        and
        \begin{align*}
            &\int (\mu_\Lambda \circ G\Ssup{\Lambda}_x)(\d\eta_\Lambda, \d\zeta_\Lambda) \, f(\eta_\Lambda, \zeta_\Lambda) \\
            &= \int \mu_\Lambda(\d \eta_\Lambda) \, b^\nu_\Lambda(x, \eta_\Lambda) \, f(\eta_{\Lambda} +\delta_{x}, \eta_\Lambda)
            + \int \mu_\Lambda^{x}(\d \eta_{\Lambda}) \, f( \eta_\Lambda -\delta_x, \eta_\Lambda) ,
        \end{align*}
        where \(b^\nu_\Lambda(x, \eta_\Lambda)\) is the birth rate of the modified dynamics introduced in \zcref{equation:definition_local_birth_rate_gibbs_exterior}.
        The corresponding infinite-volume objects are, for \(\mu \in \Pcal_{1, \theta}(\Omega)\), given by measures \(\mu \ast G_o\) and \(\mu \circ G_o\) on \(\Omega \times \Omega\), defined respectively via
        \begin{align*}
            &\int (\mu \ast G_o)(\d\eta, \d\zeta) \, f(\eta, \zeta)
            = \int \mu(\d \eta) \, b(o, \eta) \, f(\eta, \eta +\delta_{o})
            + \int \mu_o(\d \eta) \, f(\eta, \eta - \delta_o) 
        \end{align*}
        and
        \begin{align*}
            &\int (\mu \circ G_o)(\d\eta, \d\zeta) \, f(\eta, \zeta)
            = \int \mu(\d \eta) \, b(o, \eta) \, f(\eta +\delta_{o}, \eta)
            + \int \mu_o(\d \eta) \, f(\eta - \delta_o, \eta).
        \end{align*}
        We then set for \(x \in \R^d\),
        \begin{align*}
             \mu \ast G_x := (\mu \ast G_o) \circ (\theta_{-x}, \theta_{-x})^{-1},
             \quad \mu \circ G_x := (\mu \circ G_o) \circ (\theta_{-x}, \theta_{-x})^{-1}.
        \end{align*}
    \end{definition}
    \begin{remark}\label{rmk:sigmanonnegative}
    By a direct computation, we have that, for $\mu\in\Pcal_{1,\theta}(\Omega)$,
    \begin{align}\label{equation:relative_entropy_birth_death_flux_identity}
        \RelEnt\big(b(o, \cdot)\mu \,\big\vert\, \mu_o^! \big) +  \RelEnt\big(\mu_o^! \,\big\vert\, b(o, \cdot)\mu \big)
        = \RelEnt(\mu \ast G_o \,\vert\, \mu \circ G_o).
    \end{align}
    Notice that \(\RelEnt(\mu \ast G_o \,\vert\, \mu \circ G_o) \geq 0\), because \(\mu \ast G_o[1] = \mu \circ G_o[1]\), while the two summands \(\RelEnt\big(b(o, \cdot)\mu \,\big\vert\, \mu_o^! \big)\), \(\RelEnt\big(\mu_o^! \,\big\vert\, b(o, \cdot)\mu \big)\) can be negative on their own. So we have that $\sigma(\mu)\geq 0$ for any $\mu\in \Pcal_\theta(\Omega)$.
    \end{remark}
   With this notation at hand, we can now formulate the first desired representation of the change of relative entropy under the local dynamics. 
    \begin{proposition}[Finite-volume entropy-dissipation identity]\label{Prop:representation_of_finite_volume_derivative_of_entropy_for_modified_semigroup}
        For any  \(\Lambda \Subset \R^d\) and \(\mu \in \Pcal(\Omega)\) with \(\RelEnt_{\Lambda}(\mu \vert \nu) < \infty\), we have 
        \begin{align*}
            \RelEnt_{\Lambda}(\mu |\nu) - \RelEnt_{\Lambda}(\mu \Tl_t |\nu)
            = \int_0^{t} \d s \, \mathcal{J}_\L(\mu \Tl_s\vert\nu).
        \end{align*} 
    \end{proposition}

    We now express the modified Fisher information in terms of the double-layer measures introduced above.

    \begin{proposition}[Double-layer representation]\label{proposition:double-layer-to-modified-fisher}
    For all \(\Lambda \Subset \R^d\) and \(\mu_\Lambda \in \Pcal(\Omega_\Lambda)\) with finite intensity measure absolutely continuous to Lebesgue measure on \(\Lambda\), we have the representation
    \begin{equation}\label{equation:finite_volume_fisher_information_birth_death_flux_per_site_representation}
    \begin{split}
        \mathcal{J}_\L(\mu_\Lambda \vert\nu)
        = \int_{\Lambda} \d x\, \RelEnt\bigl( \mu_\Lambda \ast G\Ssup{\Lambda}_x \,\big\vert\, \mu_\Lambda  \circ G\Ssup{\Lambda}_x \bigr).
    \end{split}
    \end{equation}
\end{proposition}

To conclude the proof of \zcref{Theorem:Main_Theorem_Decrease_in_Entropy} it remains to show that we can lower bound the density limit of the Fisher information by the more explicit functional $\mu \mapsto \sigma(\mu T_s)$ as defined in \eqref{eqn:definition-sigma}. 
    \begin{proposition}\label{proposition:lower_bound_xi_mu_fisher_information_finite_volume_dynamics}
            Let \(\mu \in \mathcal{P}_{1, \theta}(\Omega)\). Then, for all $s \geq 0$ we have 
            \begin{equation}\label{equation:lower_bound_xi_mu_fisher_information_finite_volume_dynamics}
                \liminf_{n \uparrow \infty}  n^{-d} \mathcal{J}_{\Lambda_n}(\mu T\Ssup{\L_n}_s\vert\nu)\geq \sigma(\mu T_s).
            \end{equation}
    \end{proposition}
    For the proof we make use of the identity \eqref{equation:relative_entropy_birth_death_flux_identity} for \(\sigma(\mu)\) and the double-layer representation derived in \zcref{proposition:double-layer-to-modified-fisher}, see~\zcref{Sec:Proof_of_Comparison_of_xi_Lambda_and_xi} for the details. 
    \medskip 

    This allows us to conclude the proof of our first main result \zcref{Theorem:Main_Theorem_Decrease_in_Entropy}. Recall that \zcref{Theorem:Main_Theorem_Decrease_in_Entropy}\emph{(i)} is a direct consequence of~\zcref{Prop:specific_entropy_is_decreasing_in_time}. The proofs of the technical lemmas are given in Section~\ref{sec:attractor-property}. 
    
    \begin{proof}[Proof of \zcref{Theorem:Main_Theorem_Decrease_in_Entropy}(ii) and (iii).]
        Regarding \textit{(ii)}, we have for any \(\mu \in \Pcal_\theta(\Omega)\) with \(\SpecEnt(\mu \vert \nu) < \infty\), by the inequality \(\SpecEnt(\mu T_t | \nu)  \leq \widetilde \SpecEnt_{\mu}(t)\) from \zcref{Prop:specific_entropy_is_decreasing_in_time} and the existence of \(\SpecEnt(\mu | \nu)\) as limit from \zcref{prop:existence-relative-entropy}, by \zcref{Prop:representation_of_finite_volume_derivative_of_entropy_for_modified_semigroup} and Fatou's lemma and by \zcref{proposition:lower_bound_xi_mu_fisher_information_finite_volume_dynamics},
        \begin{align*}
            \SpecEnt(\mu | \nu) - \SpecEnt(\mu T_t | \nu) 
            &\geq \SpecEnt(\mu | \nu) - \widetilde{\SpecEnt}_{\mu}(t) 
            = \limsup_{n \uparrow \infty} \frac{1}{\abs{\Lambda_n}}\Big(\RelEnt_{\L_n}\bigl(\mu  |\nu\bigr) - \RelEnt_{\L_n}\bigl(\mu T\Ssup{\L_n}_t |\nu\bigr)\Big) \\
            &\geq \liminf_{n \uparrow \infty} \frac{1}{\abs{\Lambda_n}}\Big(\RelEnt_{\L_n}\bigl(\mu  |\nu\bigr) - \RelEnt_{\L_n}\bigl(\mu T\Ssup{\L_n}_t |\nu\bigr)\Big)
            = \liminf_{n \uparrow \infty}  \int_0^{t} \d s \, \frac{\mathcal{J}_{\L_n}(\mu T\Ssup{\L_n}_s\vert\nu)}{\abs{\L_n} } \\
            &\geq \int_0^t \d s\ \liminf_{n\uparrow\infty} \frac{\mathcal{J}_{\L_n}(\mu T\Ssup{\L_n}_s\vert\nu)}{\abs{\L_n} } 
            \geq \int_0^t\d s\ \sigma(\mu T_s).
        \end{align*}
        Regarding \textit{(iii)}, to see that $\sigma(\mu)$ vanishes if and only if $\mu \in \mathscr{G}_\theta$ first note that \(\sigma(\mu) = 0\) if and only if \(b(o, \cdot) \mu = \mu_o^!\). But this is equivalent to \(\mu \in \GG_\theta\) by \zcref{Prop:Reversible_Measures_equal_Gibbs_Measures}.
    \end{proof}

\subsection{The attractor property}

The main idea behind the proof of the attractor property is quite simple and already appears in \cite{holley_free_1971}. First note that if $\rho \in \mathcal{P}(\Omega)$ is a limit point of $(\mu T_t)_{t\geq 0}$ in the weak topology for some translation-invariant initial condition $\mu$, then $\rho$ must also be translation invariant. We now want to argue as follows. If $\rho$ were not a Gibbs measure, then by~\zcref{Prop:Reversible_Measures_equal_Gibbs_Measures} there is some $\delta >0$ such that the entropy production functional $\sigma(\cdot)$ evaluated in $\rho$ is bounded away from zero, i.e., 
\begin{equation*}
    \sigma(\rho) = I(\rho * G_0 \lvert \rho \circ G_0) \geq \delta. 
\end{equation*}
So if we additionally knew that $\mathscr{I}(\rho \lvert \nu) < \infty$, and that, whenever $(\mu T_t)_{t\geq 0}$ is close to $\rho$, then it spends at least a fixed amount of time in a neighborhood of $\rho$ where $\sigma(\hat{\rho}) \geq \delta/4$, then we could use a telescoping argument to get a contradiction to the assumption that $\rho \notin \mathscr{G}_\theta$. 

The following proposition does the first step of this program and tells us that $\rho$ must also have finite specific relative entropy with respect to the Gibbs measure $\nu$. 

\begin{proposition}\label{prop:existence-relative-entropy}
            For any \(\nu \in \Gcal_\theta\), the specific relative entropy $\mathscr{I}(\mu |\nu)$ with respect to $\nu$ exists for all \(\mu \in \Pcal_\theta(\Omega)\) as a limit,
            \begin{align}\label{equation:existence-relative-entropy-limit}
                \SpecEnt(\mu|\nu)
                = \lim_{n\to\infty} \frac{\RelEnt_{\L_n}(\mu|\nu)}{\abs{\L_n}}.
            \end{align}
            Furthermore, the map \(\mu \mapsto \mathscr{I}(\mu |\nu)\) is lower semicontinuous w.r.t.\ weak convergence on \(\Pcal_\theta(\Omega)\).
\end{proposition}

We can use this and the monotonicity from~\zcref{Prop:specific_entropy_is_decreasing_in_time} to obtain 
\begin{equation}\label{eqn:assumption-finite-entropy}
    \mathscr{I}(\rho \lvert \nu) \leq \liminf_{k\to \infty}\mathscr{I}(\mu T_{t_k}\lvert \nu) \leq \sup_{t\geq 0}\mathscr{I}(\mu T_t \lvert \nu) = \mathscr{I}(\mu \lvert \nu) < \infty. 
\end{equation}
So if we now additionally knew that there is some $\varepsilon>0$ and a $K \in \N$ such that 
\begin{equation}\label{eqn:main-difficulty}
    \forall \varepsilon' \in [0, \varepsilon] \ \forall k \geq K: \quad \sigma(\mu T_{t_k + \varepsilon'}) \geq \delta/4, 
\end{equation}
then we could conclude that 
\begin{equation*}
    \mathscr{I}(\mu \lvert \nu) \geq \liminf_{t\to \infty} [\mathscr{I}(\mu \lvert \nu) -\mathscr{I}(\mu T_t \lvert \nu)] \geq \int_0^\infty \d s \sigma(\mu T_s) 
    \geq \sum_{j=1}^\infty \int_0^\epsilon \d\varepsilon' \sigma(\mu T_{t_{k_j}+\varepsilon'}) \geq \sum_{j=1}^\infty \frac{\varepsilon\delta}{4} = \infty,
\end{equation*}
which contradicts \eqref{eqn:assumption-finite-entropy}, thus we must have $\rho \in \mathscr{G}_\theta$. \medskip 

Establishing \eqref{eqn:main-difficulty} is the technical part of this  program and requires the following technical lemma. 

\begin{lemma}\label{lemma:tau_L_convergence_from_weak_convergence_and_bounded_specific_entropy}
            Fix \(\nu \in \GG_\theta\). Let \(\rho \in \Pcal_\theta(\Omega)\) and let \((\rho_k)_{k \in \N}\) be a sequence in \(\Pcal_\theta(\Omega)\) such that \(\rho_k \Rightarrow \rho\) as \(k \to \infty\) and such that
            \begin{align}\label{equation:assumption_bounded_specific_entropy}
                \sup_{k \in \N} \SpecEnt(\rho_k \vert \nu) 
                < \infty.
            \end{align}
            Then, for all local, measurable \(f\) with \(\abs{f} \lesssim (1+N_\Lambda)\) for some \(\Lambda \Subset \R^d\), we have that
            \begin{align}\label{equation:convergence_of_expectations_bounded_local_measurable}
                \rho_k[f] \xrightarrow[k \to \infty]{} \rho[f].
            \end{align}
\end{lemma}

With all of these ingredients at hand, we can finally provide the proof of the attractor property. 

\begin{proof}[Proof of \zcref{Theorem:Main_Theorem_Limit_Points}]
        Suppose \(\rho \in \Pcal(\Omega)\) is some limit point in the weak topology of \((\mu T_t)_{t \geq 0}\) along the sequence $(t_k)_{k\ge 0}$, but $\rho\not\in\GG_\theta$.
        Necessarily, \(\rho \in \Pcal(\Omega)_{\theta}\), as translation invariance passes to the limit. 
        Lower semicontinuity of the specific relative entropy \zcref{prop:existence-relative-entropy} as well as the monotonicity from \zcref{Theorem:Main_Theorem_Decrease_in_Entropy} imply
        \begin{align}\label{equation:boundedness_specific_entropies_along_semigroup}
            \SpecEnt(\rho \vert \nu)
            \leq \liminf_{k \to \infty} \SpecEnt(\mu T_{t_k} \vert \nu)
            \leq \sup_{t \geq 0} \SpecEnt(\mu T_{t} \vert \nu)
            \leq \SpecEnt(\mu \vert \nu) < \infty.
        \end{align} 
        In particular, \(\rho \in \Pcal_{1, \theta}(\Omega)\) by \eqref{equation:first_moment_bound_by_finite_entropy}. 
        
        Then, by \zcref{Prop:Reversible_Measures_equal_Gibbs_Measures}, there is some \(\delta > 0\) such that
        \begin{align*}
            \RelEnt(\rho \ast G_0 \,\vert\, \rho \circ G_0)
            \geq \delta.
        \end{align*}
        We will show that there is some \(\epsilon > 0\) and \(K \in \N\) such that
        \begin{align}\label{equation:lower_bound_entropy_dissipation_goal}
            \sigma(\mu T_{t_k + \epsilon'})
            = \RelEnt(\mu T_{t_k + \epsilon'} \ast G_0 \,\vert\, \mu T_{t_k + \epsilon'} \circ G_0)
            \geq \delta/4, \quad \epsilon' \in [0, \epsilon], k \geq K.
        \end{align} 
        Hence, by \zcref{Theorem:Main_Theorem_Decrease_in_Entropy}, choosing some arbitrary subsequence \((t_{k_j})_{j}\) of \((t_k)_k\) with \(k_1 \geq K\) and such that \(t_{k_{j+1}} \geq t_{k_j} + \epsilon\) for all \(j\),
        \begin{align}
                \SpecEnt(\mu | \nu)
                &\geq
                \liminf_{t \uparrow \infty} \SpecEnt(\mu | \nu) - \SpecEnt(\mu T_{t} | \nu) \geq \int_{0}^{\infty}\d s\, \sigma(\mu T_s) 
                \geq \sum_{j} \int_{0}^{\epsilon}\d \epsilon'\, \sigma(\mu T_{t_{k_j} + \epsilon'}) 
                \geq \sum_{j} \frac{\epsilon \delta}{4}
                = \infty,
        \end{align} which is in contradiction to the assumption \(\SpecEnt(\mu | \nu) < \infty\).

        Now, by the Donsker--Varadhan formula \eqref{equation:donsker_varadhan_variational_formula_finite_measures} and approximation, there is a bounded, continuous, \(\Delta\)-local function \(F\), for some bounded cube \(\Delta \supseteq B_{R}(o)\), such that
        \begin{align}\label{equation:lower_bound_entropy_dissipation_functional_limit}
            (\rho \ast G_0)[1] + (\rho \ast G_0)[F] - (\rho \circ G_0)[\e^{F}]
            \geq \delta/2.
        \end{align}
        We have, again by the Donsker--Varadhan formula \eqref{equation:donsker_varadhan_variational_formula_finite_measures},
        \begin{align}\label{equation:lower_bound_entropy_dissipation_functional_mu_T_t}
            \RelEnt(\mu T_{t} \ast G_0 \,\vert\, \mu T_{t} \circ G_0)
            \geq (\mu T_t \ast G_0)[1] + (\mu T_t \ast G_0)[F] - (\mu T_t \circ G_0)[\e^{F}]
        \end{align} for all \(t \geq 0\), so that it suffices to compare the individual terms.
        We will only show how to obtain an \(\epsilon > 0\) and \(K \in \N\) such that
        \begin{align}\label{equation:entropy_dissipation_individual_term_small_difference}
            \big\vert (\rho \circ G_0)[\e^{F}] - (\mu T_{t_k + \epsilon'} \circ G_0)[\e^{F}] \big\vert
            \leq \delta/12, \quad \epsilon' \in [0, \epsilon], k \geq K,
        \end{align} as the other terms can be handled analogously, so that \eqref{equation:lower_bound_entropy_dissipation_goal} follows with bounds like \eqref{equation:entropy_dissipation_individual_term_small_difference} and \eqref{equation:lower_bound_entropy_dissipation_functional_mu_T_t}, \eqref{equation:lower_bound_entropy_dissipation_functional_limit}.

        We will first show how to obtain the bound \eqref{equation:entropy_dissipation_individual_term_small_difference} from that for \(\epsilon' = 0\) and then extend it to \(\epsilon' \in (0, \epsilon]\) by continuity arguments. We have
        \begin{align}\label{equation:integral_test_function_birth_death_term}
            &(\mu T_{t_k} \circ G_0)[\e^{F}] \\
            &= \int (\mu T_{t_k})(\d\eta) b(o, \eta) \e^{F(\eta + \delta_o, \eta)}  + \int (\mu T_{t_k})(\d\eta) \int_{[0, 1)^d} \eta (\d x)\,  \e^{F(\theta_x \eta - \delta_o, \theta_x\eta)}. \nonumber
        \end{align} 
        The first integrand in \eqref{equation:integral_test_function_birth_death_term}, 
        \begin{align*}
            \eta \mapsto b(o, \eta) \e^{F(\eta + \delta_o, \eta)},
        \end{align*} 
        is measurable and local by locality of \(F\) and finite range of \(b\). It is also bounded by \(\Vert b\Vert_\infty \e^{\Vert F\Vert_\infty}\). Similarly, the second integrand in \eqref{equation:integral_test_function_birth_death_term}, 
        \begin{align*}
            \eta \mapsto  \int_{[0, 1)^d} \eta (\d x)\,  \e^{F(\theta_x \eta - \delta_o, \theta_x\eta)},
        \end{align*} 
        is measurable, \((\Delta \oplus [0, 1)^d)\)-measurable, and bounded by the random variable \(\e^{\Vert F\Vert_\infty} N_{\Delta \oplus [0, 1)^d}\).
        Consequently, using \eqref{equation:boundedness_specific_entropies_along_semigroup}, we can apply \zcref{lemma:tau_L_convergence_from_weak_convergence_and_bounded_specific_entropy} to obtain convergence of both summands in \eqref{equation:integral_test_function_birth_death_term} and therefore we can choose \(K \in \N\) large enough such that, for every \(k \geq K\), \begin{align}\label{equation:entropy_dissipation_individual_term_small_difference_epsilon_prime_0}
            \big\vert (\rho \circ G_0)[\e^{F}] - (\mu T_{t_k} \circ G_0)[\e^{F}] \big\vert
            \leq \delta/24.
        \end{align}
        Let us finally see that there is an \(\epsilon > 0\) such that
        \begin{align}\label{equation:entropy_dissipation_individual_term_small_difference_epsilon_prime_greater_than_0}
            \big\vert (\mu T_{t_k + \epsilon'} \circ G_0)[\e^{F}] - (\mu T_{t_k } \circ G_0)[\e^{F}] \big\vert
            \leq \delta/24, \quad \epsilon' \in [0, \epsilon], k \in \N.
        \end{align}
        Indeed, letting \(G(\eta) := \int_{[0, 1)^d} \eta (\d x)\,  \e^{F(\theta_x \eta - \delta_o, \theta_x \eta)}\), by the Markov / semigroup property,
        \begin{align*}
            &\Big\vert \int (\mu T_{t_k + \epsilon'})(\d\eta) \, G(\eta) - \int (\mu T_{t_k})(\d\eta) \, G(\eta)  \Big\vert \\
            &\leq \mu T_{t_k}\Big[\big\vert G - T_{\epsilon'} G\big\vert \mathbf{1}_{N_{\Delta \oplus [0,1)^d} \leq M}\Big] + \mu T_{t_k}\Big[\big\vert G - T_{\epsilon'} G\big\vert \mathbf{1}_{N_{\Delta \oplus [0,1)^d} > M}\Big]  \nonumber
        \end{align*} 
        The triangle inequality and stochastic domination in the graphical construction yield for \(M \geq \Vert b \Vert_\infty \abs{\Delta \oplus [0,1)^d}\) that
        \begin{align*}
             &\mu T_{t_k}\Big[\big\vert G - T_{\epsilon'} G\big\vert \mathbf{1}_{N_{\Delta \oplus [0,1)^d} > M}\Big] \\
             &\leq \mu T_{t_k}\Big[\e^{\Vert F\Vert_\infty} \Big(N_{\Delta \oplus [0,1)^d} + \e^{-\epsilon'} N_{\Delta \oplus [0,1)^d} + \Vert b\Vert_\infty \vert\Delta \oplus [0,1)^d\vert  (1 - \e^{-\epsilon'})  \Big) \mathbf{1}_{N_{\Delta \oplus [0,1)^d} > M}\Big] \\
             &\leq 2\e^{\Vert F\Vert_\infty} \mu T_{t_k}\Big[N_{\Delta \oplus [0,1)^d} \mathbf{1}_{N_{\Delta \oplus [0,1)^d} > M}\Big]. 
        \end{align*}
        Here, we can first choose, by \eqref{equation:boundedness_specific_entropies_along_semigroup} and the uniform integrability established in the proof of \zcref{lemma:tau_L_convergence_from_weak_convergence_and_bounded_specific_entropy}, \(M\) large enough, so that
        \begin{align*}
            \sup_{k \in \N} 2\e^{\Vert F\Vert_\infty} \mu T_{t_k}\Big[N_{\Delta \oplus [0,1)^d} \mathbf{1}_{N_{\Delta \oplus [0,1)^d} > M}\Big]
            \leq \delta / 96.
        \end{align*} 
        Now, we choose \(\epsilon > 0\) small enough, so that
        \begin{align*}
            &\sup_{k \in \N} \mu T_{t_k}\Big[\big\vert G - T_{\epsilon'} G\big\vert \mathbf{1}_{N_{\Delta \oplus [0,1)^d} \leq M}\Big] 
            \leq \delta / 96
        \end{align*} 
        for all \(\epsilon' \in [0, \epsilon]\).
        This is possible, because
        \begin{align*}
            \big\vert G(\eta) - T_{\epsilon'} G(\eta)\big\vert
            \xrightarrow[\epsilon' \downarrow 0]{} 0,
        \end{align*} 
        uniformly for \(\eta\) with \(N_{\Delta \oplus [0,1)^d}(\eta) \leq M\), \(M\) fixed, by decomposing \(T_{\epsilon'} G\) into \(G\) times the probability that nothing happens for the dynamics in \(\Delta \oplus [0,1)^d\) during \([0, \epsilon]\) and into the other case. By \((\Delta \oplus [0,1)^d)\)-locality and the bound \(G \leq \e^{\Vert F\Vert_\infty} N_{\Delta \oplus [0,1)^d}\), one has indeed only to control the probability that nothing happens and the expected number of particles inside \((\Delta \oplus [0,1)^d)\) times the indicator that a birth or death happens inside \((\Delta \oplus [0,1)^d)\) during \([0, \epsilon]\). The bounded birth rates \(b\) and the bounded total death rate of the initial condition \(\eta\) in \(\Delta \oplus [0,1)^d\) by \(N_{\Delta \oplus [0,1)^d \leq M}\) easily ensure this control.
        
        By an analogous argument we get that, for \(\epsilon > 0\) small enough and for all \(\epsilon' \in [0, \epsilon]\),
        \begin{align*}
            &\sup_{k \in \N} \big\vert (\mu T_{t_k + \epsilon'} \circ G_0)[\e^{F}] - (\mu T_{t_k } \circ G_0)[\e^{F}] \big\vert \\
            &\leq \sup_{k \in \N} \Big\vert \int (\mu T_{t_k + \epsilon'})(\d\eta) b(o, \eta) \e^{F(\eta + \delta_o, \eta)} - \int (\mu T_{t_k})(\d\eta) b(o, \eta) \e^{F(\eta + \delta_o, \eta)}  \Big\vert \nonumber \\
            &\qquad + \sup_{k \in \N}  \Big\vert \int (\mu T_{t_k + \epsilon'})(\d\eta) G(\eta) - \int (\mu T_{t_k})(\d\eta) G(\eta) \big\vert \nonumber \\
            &\leq \delta/48 + \delta/48 = \delta/24 \nonumber
        \end{align*}  
        i.e.\ \eqref{equation:entropy_dissipation_individual_term_small_difference_epsilon_prime_greater_than_0} holds.
        Finally, \eqref{equation:entropy_dissipation_individual_term_small_difference_epsilon_prime_0} and \eqref{equation:entropy_dissipation_individual_term_small_difference_epsilon_prime_greater_than_0} together yield \eqref{equation:entropy_dissipation_individual_term_small_difference}.
    \end{proof}

\section{Construction of the infinite-volume birth-and-death dynamics}\label{sec:construction_dynamics}

We present here a re-elaboration of the results of~\cite{garcia2006spatial}, which hold for any finite-range birth rate that is uniformly bounded from above meaning that, in particular, they apply to our setting and allow us to construct the Markov process described in the above \zcref{Definition:Markov_Process}.

To avoid confusion with our standing assumptions on $b$, we call such a function $q$. In other words, let \(q \colon \R^d \times \Omega \to [0, \infty)\) be a measurable and bounded function, \(\Vert q\Vert_\infty < \infty\), with finite range, i.e.\ such that there exists an \(R > 0\) such that
\begin{align*}
    q(x, \eta)
    = q(x, \eta_{B_R(x)})
\end{align*} 
for all \(x \in \R^d\) and \(\eta \in \Omega\).
We want to construct, for every initial condition \(\eta \in \Omega\), the infinite-volume birth-and-death process \(X^\eta = (X^{\eta}_t)_{t \geq 0}\) with values in \(E\) and formal generator
\begin{align}\label{definition:formal_generator}
    (\LL f)(\eta) 
    := \int_{\R^d} \d x \ q(x,\eta) [f(\eta+\delta_x) - f(\eta)] + \sum_{x\in\eta} [f(\eta-\delta_x) - f(\eta)].
\end{align}

The notion of a process with formal generator \zcref{definition:formal_generator} is made precise via the following martingale problem, cf.~\cite[Section 2.1]{garcia2006spatial}.
\begin{definition}[Martingale problem for \((\mathscr{L}, \mathcal{D}, \mu)\)]\label{definition:martingale_problem}
    Fix \(\mu \in \Pcal(\Omega)\) and denote by \((X_t)_{t \geq 0}\) the coordinate process on \(E = D([0, \infty); \Omega)\).
    Consider the class of test functions
    \begin{align}
        \mathcal{D}
        := \{\eta \mapsto \e^{- \int \eta(\d x) g(x)} \,\colon\, g \in C_c(\R^d), g \geq 0\}.
    \end{align}
    We say that a probability measure \(P\) on \(E\) \emph{solves the martingale problem for} \((\mathscr{L}, \mathcal{D}, \mu)\), if \(P \circ X_0^{-1} = \mu\) and if, for any \(f \in \mathcal{D}\),
    \begin{align}\label{equation:martingale_problem}
        f(X_t) - f(X_0) - \int_0^t (\mathscr{L} f)(X_s) \d s, \quad t \geq 0,
    \end{align} 
    is a martingale w.r.t.\ the canonical filtration \(\sigma(X_s \colon 0 \leq s \leq t)\), \(t \geq 0\).
\end{definition}
\begin{remark}[Relation to uniqueness of associated Markov processes]\label{remark:martingale_problem_markov_process}
    For the reader not acquainted with martingale problems, we recall here a relevant part of the classical connection to Markov processes.
    Suppose \(X\) is an \(\Omega\)-valued, time-homogeneous Markov process with c\`adl\`ag paths, initial distribution \(\mu\), and transition semigroup \((T_t)_{t \geq 0}\) on the bounded measurable functions \(\Omega \to \R\).
    Assume that the forward equation
    \begin{align}\label{equation:forward_equation_L}
        T_t f
        = f + \int_0^t \, T_s\mathscr{L}f \, \d s, \quad t \geq 0,
    \end{align} 
    holds for all \(f \in \mathcal{D}\). Then the law of \(X\) solves the martingale problem for \((\mathscr{L}, \mathcal{D}, \mu)\).
    
    Indeed, for \(f \in \mathcal{D}\) we have that \(f\) and \(\mathscr{L} f\) are bounded, so that we can compute with the Markov property and \eqref{equation:forward_equation_L},
    \begin{align*}
        &\mathbb{E}[f(X_t) - f(X_s) \,\vert\, \sigma(X_r \colon 0 \leq r \leq s) ]
        = (T_{t-s} f)(X_s) - f(X_s)
        = \int_0^{t-s} \, (T_r \mathscr{L}f)(X_s) \, \d r \\
        &= \int_s^{t} \, (T_{r-s} \mathscr{L}f)(X_s) \, \d r 
        = \int_s^{t} \, \mathbb{E}[(\mathscr{L}f)(X_{s + (r-s)}) \,\vert\, \sigma(X_u \colon 0 \leq u \leq s) ] \, \d r \\
        &= \mathbb{E}\Big[\int_s^t \, (\mathscr{L}f)(X_r) \, \d r \,\Big\vert\, \sigma(X_u \colon 0 \leq u \leq s) \Big],
    \end{align*} 
    meaning that \eqref{equation:martingale_problem} is a martingale.
\end{remark}

Indeed, in the sense of the martingale problem, there is a unique process with given initial distribution and with formal generator \(\mathscr{L}\), and it is Markovian.
\begin{proposition}[Solutions of the martingale problem]\label{proposition:martingale_problem}
    For any \(\mu \in \Pcal(\Omega)\), there is a unique probability measure \(P^\mu\) solving the martingale problem for \((\mathscr{L}, \mathcal{D}, \mu)\) from \zcref{definition:martingale_problem}.
    Moreover, \(P^\mu\) is the law of a time-homogeneous \(\Omega\)-valued Markov process with c\`adl\`ag paths.
\end{proposition}

There is a rather explicit construction of the solution to the martingale problem for \((\mathscr{L}, \mathcal{D}, \mu)\) as solution to the following SDE introduced in \cite{garcia2006spatial}.

\begin{definition}[birth-and-death SDE, graphical construction]\label{definition:SDE}
    Fix \(\mu \in \Pcal(\Omega)\).
    Consider a complete probability space \((\Xi, \mathcal{G}, \mathbb{P})\) equipped with a filtration \((\mathcal{G}_t)_{t \geq 0}\).
    Denote \(\mathfrak{N} := \{H \subseteq \Xi \,\colon\, \exists G \in \mathcal{G} \colon H \subseteq G, \mathbb{P}(G) = 0 \}\) for the family of \(\mathbb{P}\)-null sets. Additionally require that the filtered probability space \((\Xi, \mathcal{G}, (\mathcal{G}_t)_{t \geq 0}, \mathbb{P})\)  also has the property that
    \begin{align*}
        \mathfrak{N} \subseteq \mathcal{G}_0,
    \end{align*}
    and that it supports the following random variables.
    
    There is a \(\mathcal{G}_0\)-measurable lifespan-marked initial configuration \(\widehat{\eta}\), such that \(\eta := \widehat{\eta}(\cdot \times [0, \infty)) \sim \mu\), and \(\widehat{\eta}\) is obtained, conditionally on \(\eta\), by independently assigning a unit-exponential random variable to each occurrence of a point in \(\eta\).
    Furthermore, there is a Poisson random measure \(\Ncal\) on \(\R^d \times [0, \infty)^3\) with intensity measure \(\d x \otimes \d s \otimes \e^{-r} \d r \otimes \d u\) compatible with the filtration in the following sense,
    \begin{align*}
        \sigma(\Ncal_{\R^d \times [0, t] \times [0, \infty)^2}) \subseteq \mathcal{G}_t, \quad
        \sigma(\Ncal_{\R^d \times (t, \infty) \times [0, \infty)^2}) \perp\!\!\!\perp \mathcal{G}_t, \quad t \geq 0.
    \end{align*}
    In particular, \(\mathcal{N}\) is independent from \(\widehat{\eta}\).
    
    Now, for \(B \Subset \R^d\) and \(t \geq 0\), consider the SDE
    \begin{equation}\label{equation:SDE}
    \begin{split}
        X_t(B)
        &= \int_{B \times (0, t] \times [0, \infty)^2} \Ncal(\d x, \d s, \d r, \d u) {\1}_{(t-s, \infty)}(r) {\1}_{[0, q(x, X_{s-})]}(u)
        \\
        &+ \int_{B \times [0, \infty)} \widehat{\eta}(\d x, \d r) {\1}_{(t, \infty)}(r).
    \end{split}
    \end{equation}

    We call \emph{solution} to the SDE \zcref{equation:SDE} any \(\Omega\)-valued, \((\mathcal{G}_t)_{t \geq 0}\)-adapted process \(X = (X_t)_{t \geq 0}\) with  c\`adl\`ag paths such that \zcref{equation:SDE} holds \(\mathbb{P}\)-a.s.\ simultaneously for every \(B \Subset \R^d\) and \(t \geq 0\). We call a solution \(X\) to the SDE \zcref{equation:SDE} \emph{strong} if it is adapted to the filtration 
    \begin{align*}
        \sigma(\widehat{\eta}, \Ncal_{\R^d \times [0, t] \times [0, \infty)^2}) \lor \mathfrak{N}, \quad t \geq 0.
    \end{align*}
    We say that \emph{pathwise uniqueness} holds for \zcref{equation:SDE} if \(\mathbb{P}(X_t = Y_t \text{ for every } t \geq 0) = 1\) whenever \(X, Y\) are two solutions of \zcref{equation:SDE} defined on the same filtered probability space, driven by the same \(\widehat{\eta}\) and \(\mathcal{N}\).
\end{definition}
\begin{remark}
    The solution concept of Garcia--Kurtz requires only that the underlying probability space is complete, not that the \(\mathbb{P}\)-null sets are included in the filtration. We think that it then is not always correct that there is an adapted solution to the SDE with everywhere c\`adl\`ag paths. One can then only require almost surely c\`adl\`ag paths on an arbitrary such filtered probability space. Including \(\mathbb{P}\)-null sets in the filtration allows for modifications on null-sets while preserving adaptedness.
\end{remark}

\begin{proposition}[Solutions of the birth-and-death SDE]\label{proposition:SDE}
    Consider the situation of \zcref{definition:SDE}.
    Then, 
    \begin{enumerate}[(i)]
        \item The SDE \zcref{equation:SDE} admits a strong solution \(X\) and pathwise uniqueness holds.

        \item In particular, there is uniqueness in law: for a fixed initial distribution \(\mu \in \Pcal(\Omega)\), independently of the underlying filtered probability space, all solutions to the SDE \zcref{equation:SDE} have the same law as \(E\)-valued random variables.
        
        \item The law of any solution \(X\) of \zcref{equation:SDE} with initial condition \(\mu \in \Pcal(\Omega)\), as \(E\)-valued random variable, is the unique solution of the martingale problem for \((\mathscr{L}, \mathcal{D}, \mu)\) from \zcref{proposition:martingale_problem}.

        \item Any solution \(X = (X_t)_{t \geq 0}\) of \zcref{definition:SDE} is time-homogeneous Markov w.r.t.\ its natural filtration. Moreover, writing \(X^\eta\) for a solution with fixed initial condition \(\eta \in \Omega\), the laws of \(X^\eta = (X^{\eta}_t)_{t \geq 0}\), \(\eta \in \Omega\), constitute a time-homogeneous Markov family on \(\Omega\).
        
        \item Moreover, if \(\mu \in \Pcal(\Omega_{\mathrm{s}})\), the strong solution \(X\) can be chosen as an \(E_{\mathrm{s}}\)-valued random variable.
        
        \item Suppose additionally that the birth rate \(q\) is translation invariant,
        \begin{align}\label{equation:translation_covariance_birth_rate_q}
            q(\theta_z x, \theta_z \eta) = q(x, \eta)
        \end{align} 
        for all \(x, z \in \R^d\) and \(\eta \in \Omega\).
        Then, writing \(P^\mu\) for the law of the solution \(X\) with initial distribution \(\mu \in \Pcal(\Omega)\), we have
        \begin{align*}
            P^{\mu \circ \theta_z^{-1}} 
            = P^{\mu} \circ \theta_z^{-1}.
        \end{align*}
    \end{enumerate}  
\end{proposition}

\begin{proposition}[Transition semigroup \((T_t)_{t \geq 0}\)]\label{proposition:transition_semigroup}
    For each \(\eta \in \Omega\), let \(X^\eta = (X_t^\eta)_{t \geq 0}\) be any \(\Omega\)-valued process with c\`adl\`ag paths solving the SDE \eqref{equation:SDE} with fixed initial condition \(\eta\).
    Define
    \begin{align}
        (T_t f)(\eta)
        := \mathbb{E}[f(X_t^\eta)], \quad t \geq 0, \, \eta \in \Omega,
    \end{align} for \(f \colon \Omega \to \R\) bounded, measurable. 
    Then, 
    \begin{enumerate}[(i)]
        \item  \((T_t)_{t \geq 0}\) is a Markov semigroup on the bounded, measurable functions \(\Omega \to \R\). In particular, \(T_t f \colon \Omega \to \R\) is measurable.
        \item  In the case that \(q\) is translation invariant, \((T_t)_{t \geq 0}\) is translation covariant, i.e.\
        \begin{align*}
            (T_t f) \circ \theta_z
            = T_t (f \circ \theta_z)
        \end{align*} 
        for any bounded measurable \(f\), \(z \in \R^d\) and \(t \geq 0\).

        \item For \(f \in \Lcal\), the forward equation
        \begin{align}
            T_t f(\eta)
            = f(\eta) + \int_0^t \, T_s\mathscr{L} f(\eta) \, \d s, \quad t \geq 0,
        \end{align} 
        holds pointwise for all \(\eta \in \Omega\).
    \end{enumerate}
\end{proposition}

The proof of these results can be found in \zcref{sec:appendix:infiniteconstruction}. We reformulate the existence result of \cite{garcia2006spatial} and its proof via \emph{solution maps}: deterministic maps that take the lifespan-marked initial condition \(\widehat{\eta}\) and a random counting measure \(\Ncal\) as input and return a birth-and-death path starting at \(\eta\). The solution maps are then applied to \(\widehat{\eta}\), obtained by independently associating to each particle in \(\eta\) a unit-exponential lifespan, and to the Poisson random measure \(\Ncal\).
This approach allows for a somewhat transparent validation of the measurability in the initial condition and the Markov property. Furthermore, it also allows us to directly define couplings of the infinite-volume birth-and-death process with finite-volume versions of it.
The idea of such solution maps for SDEs is of course classical and underlies, e.g., the celebrated Yamada--Watanabe theorem.

\section{Comparison with finite-volume dynamics: finite speed of propagation}\label{section:finite_volume_dynamics}

Fix \(\nu \in \GG_\theta\) and let \(\Lambda \Subset \R^d\).
Choose a regular version \(\nu_{\Lambda \,\vert\, \Lambda^c}\) of the conditional distribution under \(\nu\) of the particles in \(\Lambda^c\) given the configuration inside \(\Lambda\), i.e.\ such that \(\nu_{\Lambda \,\vert\, \Lambda^c}\) is a kernel from \(\Omega_\Lambda\) to \(\Omega_{\Lambda^c}\) and
\begin{align*}
    \int \nu_\Lambda(\d\eta_\Lambda) \int \nu_{\Lambda \,\vert\, \Lambda^c}(\d\zeta_{\Lambda^c} \vert \eta_\Lambda) \, f(\eta_\Lambda \zeta_{\Lambda^c})
    = \int \nu(\d\eta) \, f(\eta)
\end{align*}  
holds for all measurable \(f \colon \Omega \to [0, \infty)\).
Define the birth rate \(b_\Lambda^\nu\) by
\begin{align}\label{equation:definition_local_birth_rate_gibbs_exterior}
    b_\Lambda^\nu(x, \eta_\Lambda)
    := \int \nu_{\Lambda \,\vert\, \Lambda^c}(\d\zeta_{\Lambda^c} \vert \eta_\Lambda)\, b(x, \eta_\Lambda \zeta_{\Lambda^c}), 
    \quad x \in \Lambda, \eta_\Lambda \in \Omega_\Lambda,
\end{align} 
which is jointly measurable, bounded by \(\Vert b \Vert_\infty\) and agrees with the infinite-volume birth rate \(b\) for \(x\) with \(\mathrm{dist}(x, \Lambda^c) > R\), and consider the formal generator
\begin{align*}
    (\mathscr{L}_\Lambda f)(\eta_\Lambda)
    = \int_{\Lambda} \d x\, b^\nu_\Lambda(x, \eta_\Lambda) \big[f(\eta_\Lambda+\delta_x) - f(\eta_\Lambda)\big] + \sum_{x\in\eta_\Lambda} \big[f(\eta_\Lambda -\delta_x) - f(\eta_\Lambda)\big]
\end{align*} 
for measurable \(f \colon \Omega_\Lambda \to \R\) such that the right-hand side is well-defined for all \(\eta_\Lambda \in \Omega_\Lambda\).
We will see that it is associated with the Markov process \((X^{(\Lambda, \eta)}_t)_{t \geq 0}\) which is, for \(\eta \in \Omega_\Lambda\), the solution to the SDE
\begin{equation}\label{equation:SDE_local}
\begin{split}
        X\Ssup{\L,\eta}_t(B) 
        = \int_{B \times (0,t] \times [0,\infty)^2} \Ncal(\d x,\d s,\d r, \d u)\, &\1_{[0,b^\nu_\L(x,X\Ssup{\L,\eta}_{s-})]}(u) \1_{(t-s,\infty)}(r) \\
        &+ \int_{B\times [0,\infty)}\hat{\eta}(\d x, \d r)\, \1_{(t,\infty)}(r),
\end{split}
\end{equation} 
for \(B \subseteq \Lambda\) and \(t \geq 0\), analogous to \zcref{definition:SDE}, with \(\widehat{\eta}\) a marking of \(\eta\) with i.i.d.\ unit-exponential random variables and an independent Poisson random measure \(\Ncal\) on \(\R^d \times [0, \infty)^3\) with intensity \(\d x \otimes \d s \otimes \e^{-r} \d r \otimes \d u\).

Indeed, this SDE is well-posed as in the infinite-volume case:
\begin{proposition}[SDE for local dynamics]\label{proposition:SDE_local}
    The following hold.
    \begin{enumerate}[(i)]
        \item \emph{Unique solution of the SDE.} For every \(\eta \in \Omega_\Lambda\) and on every filtered probability space as in \zcref{definition:SDE} with initial distribution \(\delta_\eta\), the SDE \eqref{equation:SDE_local} admits an \(\Omega_\Lambda\)-valued strong solution \(X\Ssup{\L,\eta}\), with c\`adl\`ag paths and locally only finitely many jumps, and pathwise uniqueness holds.

        \item \emph{Markov property.} For each fixed \(\eta \in \Omega_\Lambda\), all solutions \(X\Ssup{\L,\eta}\) of \zcref{equation:SDE_local}, in the sense of \zcref{definition:SDE}, with initial configuration \(\eta\) have the same law, independently of the underlying filtered probability space, and  \(X\Ssup{\L,\eta}\) is a time-homogeneous Markov process w.r.t.\ the natural filtration.
        Furthermore, the laws of  \(X\Ssup{\L,\eta} = (X\Ssup{\L,\eta}_t)_{t \geq 0}\), \(\eta \in \Omega_\Lambda\), constitute a time-homogeneous Markov family on \(\Omega_\Lambda\).

        \item \emph{Martingale problem.} Let \(\mathcal{D}_\Lambda = \{\eta_\Lambda \mapsto \e^{-\int_\Lambda \eta_\Lambda(\d x) g(x)} \,\colon\, g \in C_c(\R^d), g \geq 0 \}\). Then, for every \(\mu_\Lambda \in \Pcal(\Omega_\Lambda)\), the unique solution to the martingale problem for \((\mathscr{L}_\Lambda, \mathcal{D}_\Lambda, \mu_\Lambda)\), defined analogously to \zcref{definition:martingale_problem}, is given by \(P^{\Lambda, \mu} := \int \mu_\Lambda(\d\eta) \, \mathrm{Law}(X^{(\Lambda, \eta)})\).
    \end{enumerate}
\end{proposition} 
The proof is provided in \zcref{section:finite_volume_dynamics_construction}.

We now let the transition semigroup \((\Tl_t)_{t \geq 0}\) of our local dynamics \(\Xl\) be defined via
\begin{align}\label{equation:definition_semigroup_local_dynamics}
    (\Tl_t f)(\eta)
    := \mathbb{E}[f(X_t^{(\Lambda, \eta)})], \quad t \geq 0, \, \eta \in \Omega_\Lambda,
\end{align} 
for all bounded, measurable \(f \colon \Omega_\Lambda \to \R\).
We will also denote by \(\Tl_t f\) the corresponding function for all \(f \colon \Omega_\Lambda \to \R \cup \{-\infty, \infty\}\) such that the right-hand side of \eqref{equation:definition_semigroup_local_dynamics} is well-defined, e.g.\ for non-negative measurable \(f\).

\begin{proposition}[Transition semigroup for local dynamics]\label{proposition:finite_volume_dynamics_semigroup}
    The following hold.
    \begin{enumerate}[(i)]
        \item  \((\Tl_t)_{t \geq 0}\) is a Markov semigroup on the bounded, measurable functions \(\Omega_\Lambda \to \R\). In particular, \(\Tl_t f \colon \Omega_\Lambda \to \R\) is measurable for all \(t \geq 0\).

        \item The operators \((\Tl_t)_{t \geq 0}\) extend to a \(C_0\)-semigroup on \(L^1(\nu_\Lambda)\). The action of its generator agrees \(\nu_\Lambda\)-a.e.\ with the formal generator \(\mathscr{L}_\Lambda\) and the domain of its generator contains the bounded measurable functions on \(\Omega_\Lambda\).
    \end{enumerate}
\end{proposition}

Moreover, the GNZ equations imply reversibility of \(\nu_\Lambda\) w.r.t.\ the local dynamics \((\Xl_t)_{t \geq 0}\).

\begin{proposition}[Finite-volume reversibility]\label{proposition:finite_volume_reversibility}
    Let \(f, g \colon \Omega_\Lambda \to [0, \infty]\) be measurable.
    Then
    \begin{align*}
        \nu_\Lambda[f(\Tl_t g)]
        = \nu_\Lambda[(\Tl_t f) g]
    \end{align*} 
    for all \(t \geq 0\).
\end{proposition}
The proofs of  \zcref{proposition:finite_volume_dynamics_semigroup} and \zcref{proposition:finite_volume_reversibility} are given in \zcref{section:finite_volume_dynamics_semigroup}.

For any initial condition \(\eta\in\Omega\), we have a canonical coupling between the infinite-volume Markov process \((X_t\Ssup{\eta})_{t \geq 0}\) and the Markov process \((X_t\Ssup{\Lambda, \eta})_{t \geq 0}\) in the finite volume \(\Lambda\) with stochastic boundary conditions, using a common Poisson random measure \(\mathcal{N}\) for the birth proposals and common i.i.d.\ exponentially distributed lifespans for shared particles of the initial configuration. Existence of such a coupling is shown in \zcref{section:finite_volume_dynamics_construction}. Herein, \((X_t\Ssup{\Lambda, \eta})_{t \geq 0}\) means for \(\eta \in \Omega\) that \(X_t\Ssup{\Lambda, \eta}\) is started in \(\eta_\Lambda \in \Omega_\Lambda\).
We use this coupling to compare the corresponding evolutions, obtaining the following ``finite-speed of propagation'' result.
\begin{lemma}[Finite speed of propagation]\label{Lemma:comparison_finite_infinite}
    Let \(T > 0\) be fixed. For any box $\L\Subset\R^d$ with each side length at least \(1\), let $\widehat{\L} \Subset \R^d$ be some superset of $\L$, with \(\ell := \mathrm{dist}(\Lambda, \widehat{\Lambda}^{c}) > 0\). We have
    \begin{align*}
        \mathbb{P}\Big(\exists t \in [0, T] \colon X_{t, \Lambda}\Ssup{\widehat{\Lambda}, \eta} \neq X_{t, \Lambda}^{(\eta)}  \Big)
        \lesssim_{d, R} \abs{\Lambda} 
        \e^{- \ell \rho(\ell)},
    \end{align*}
    for some increasing \(\rho\) (depending on \(T, d, \Vert b \Vert_\infty, R\)) with \(\rho(\ell) \xrightarrow[\ell\uparrow \infty]{} \infty\).
\end{lemma}
\begin{proof}
    First note that \(X_{t, \Lambda}\Ssup{\widehat{\Lambda}, \eta} \neq X_{t, \Lambda}^{(\eta)}\) for some \(t \in [0, T]\) is possible only if the following event happens
    \begin{align*}
        E_{\Lambda}(T,\ell)
        &:= \Big\{\exists n \in \mathbb{N} ,\, 0 < t_1 < \dots < t_n \leq T,\, x_1, \dots, x_n \in \mathbb{R}^d \text{ such that }\\
        &\qquad\quad \mathcal{N}\big(\{x_i\}  \times \{t_i\} \times \mathbb{R}_+ \times [0, \Vert b\Vert_\infty] \big) = 1 \text{ for } i = 1, \dots, n, \\
        &\qquad\quad\mathrm{dist}(x_1, \widehat{\Lambda}^c) \leq R,\, x_n \in \Lambda, \ \vert x_i -  x_{i+1} \vert \leq R \text{ for } i = 1, \dots, n-1
        \Big\},
    \end{align*}
    because every disagreement \(X_{t, \Lambda}\Ssup{\widehat{\Lambda}, \eta} \neq X_{t, \Lambda}^{(\eta)}\) has to originate in the boundary layer of \(\widehat{\Lambda}\), i.e.\ at a point \(x \in \widehat{\Lambda}\) with \(\mathrm{dist}(x, \widehat{\Lambda}^c) \leq R\), and propagate through successive births.
    In the definition of the event \(E_{\Lambda}(T,\ell)\), there are necessarily \(n \geq \Big\lceil \ell/R \Big\rceil\) successive birth proposals as
    \begin{align*}
        \ell = \mathrm{dist}(\Lambda, \widehat{\Lambda}^c) 
        \leq \mathrm{dist}(x_n, \widehat{\Lambda}^c)
        \leq \mathrm{dist}(x_1, \widehat{\Lambda}^c) + \sum_{i = 1}^{n-1} \abs{x_{i+1} - x_i}
        \leq R n
    \end{align*} for \(x_1, \dots, x_n\) from the definition of \(E_{\Lambda}(T,\ell)\).
    
    Tile \(\mathbb{R}^d\) into boxes \(Q_j\) of side length \(2R\), and call two boxes \emph{neighboring} if they share any boundary points. In particular they are allowed to coincide.
    We then have \(E_\Lambda(T,\ell) \subseteq \widetilde{E}_\Lambda(T,\ell)\), where
    \begin{align*}
        \widetilde{E}_\Lambda(T,\ell)
        &:= 
        \Big\{\exists n \geq \Big\lceil \ell/R \Big\rceil,\, 0 < t_1 < \dots < t_n \leq T,\text{ and boxes } Q_{j_1}, \dots, Q_{j_n} \text{ such that } \\
        &\qquad\quad  \mathcal{N}\big(Q_{j_i}  \times \{t_i\} \times \mathbb{R}_+ \times [0, \Vert b\Vert_\infty]\big) = 1 \text{ for } i = 1, \dots, n,\\
        &\qquad\quad Q_{j_n} \cap \Lambda \neq \emptyset \text{ and the boxes } Q_{j_i}, Q_{j_{i+1}} \text{ are neighboring for } i = 1, \dots, n-1
        \Big\}.
    \end{align*}
    For a fixed, ordered sequence \(Q_{j_1}, \dots, Q_{j_n}\) of (consecutively neighboring) boxes and the waiting times
    \begin{align*}
        &\tau_0 := 0, \\
        &\tau_i = \inf\{t > \tau_{i-1} \,\colon\, \mathcal{N}\big(Q_{j_i}  \times \{t\} \times \mathbb{R}_+ \times [0, \Vert b\Vert_\infty]\big) = 1 \},  
        \quad 0 < i \leq n,
    \end{align*} 
    we have that \(\tau_n = \sum_{i = 1}^{n} (\tau_i - \tau_{i-1})\) is a sum of i.i.d.\ \(\mathrm{Exp}(\Vert b\Vert_\infty \abs{Q_1})\)-distributed random variables, by the strong Markov property of the Poisson random measure \(\Ncal\).
    Hence, with \(Z \sim \mathrm{Poisson}(T\Vert b\Vert_\infty \abs{Q_1})\), we get
    \begin{equation}
    \begin{split}
        &\mathbb{P}\Big(\exists 0 < t_1 < \dots < t_n \leq T  \text{ s.t. } \mathcal{N}\big(Q_{j_i}   \times \{t_i\} \times \mathbb{R}_+ \times  [0, \Vert b\Vert_\infty]\big) = 1 \text{ for } i = 1, \dots, n\Big)  \\
        &= \mathbb{P}\Big(\tau_n \leq T\Big) 
        =  \mathbb{P}\Big(Z \geq n \Big)
        \leq \exp\Big(-n \log\Big(\frac{n}{eT \Vert b\Vert_\infty \abs{Q_1}}\Big) \Big).
    \end{split}
    \end{equation}
    We can now bound the probability of \(\widetilde{E}_\Lambda(T, \ell)\) via a union bound, because there are \(\sim_{d, R} \abs{\Lambda}\) boxes \(Q_{j_n}\) with \(Q_{j_n} \cap \Lambda \neq \emptyset\) and less than \((3^d)^n\) sequences \(Q_{j_1}, \dots, Q_{j_n}\) of neighboring boxes ending at \(Q_{j_n}\), yielding
    \begin{align*}
        \mathbb{P}\big(\widetilde{E}_\Lambda(T, \ell)\big)
        &\lesssim_{d, R} \abs{\Lambda} \sum_{n \geq \lceil \ell/ R \rceil} 3^{dn} \exp\Big(-n \log\Big(\frac{n}{eT \Vert b\Vert_\infty \abs{Q_1}}\Big) \Big).
    \end{align*}
    In particular, the right-hand side decays faster than exponential in \(\ell\) for fixed \(\Lambda, T\).
\end{proof}

The preceding upper bound on the probability that the coupling of the two processes fails can be converted into an estimate of the total variation distance of the laws $(\mu T_t)_\Lambda$ and $(\mu T^{(\hat{\Lambda})}_t)_\Lambda$ as claimed in~\zcref{Lemma:finite_speed_of_propagation_version_1}. 
\begin{proof}[Proof of~\zcref{Lemma:finite_speed_of_propagation_version_1}]
    Clearly,
    \begin{align*}
        \big\vert \mu[g T_t f] - \mu[g T_t\Ssup{\widehat{\Lambda}}f] \big\vert
        \leq \int \mu(\d\eta) \, \abs{g(\eta)} \mathbb{E}\Big[\big\lvert f(X_{t, \Lambda}\Ssup{\widehat{\Lambda}, \eta}) - f(X_{t, \Lambda}\Ssup{\eta}) \big\rvert \1_{\{X_{t, \Lambda}\Ssup{\widehat{\Lambda}, \eta} \neq X_{t, \Lambda}\Ssup{\eta}\}} \Big],
    \end{align*} 
    and therefore, by \zcref{Lemma:comparison_finite_infinite},
    \begin{align*}
        \big\vert \mu[g T_t f] - \mu[g T_t\Ssup{\widehat{\Lambda}}f] \big\vert
        &\leq 2 \Vert f\Vert_\infty \Vert g\Vert_\infty \sup_{\eta \in \Omega}  \mathbb{P}(X_{t, \Lambda}\Ssup{\widehat{\Lambda}, \eta} \neq X_{t, \Lambda}\Ssup{\eta})
        \\
        &\lesssim_t \abs{\Lambda} \Vert f\Vert_\infty \Vert g\Vert_\infty  \e^{-\mathrm{dist}(\Lambda, \widehat{\Lambda}^{c}) \rho(\mathrm{dist}(\Lambda, \widehat{\Lambda}^{c}))}.
    \end{align*} 
    For the second statement, note that by the locality and boundedness assumptions on \(f_\Lambda\) as well as by Cauchy-Schwarz,
    \begin{align*}
        &\big\vert \mu T_t[f_\L] - \mu T_t\Ssup{\widehat{\Lambda}}[f_\L] \big\vert
        \lesssim \int \mu(\d\eta) \, \mathbb{E}\Big[2 \big(N_\Lambda(\eta) + \Ncal(\Lambda \times [0, t] \times \mathbb{R}_+ \times [0, \Vert b \Vert_\infty]) \big) \1_{\{X_{t, \Lambda}\Ssup{\widehat{\Lambda}, \eta} \neq X_{t, \Lambda}\Ssup{\eta}\}} \Big] \\
        &\lesssim \abs{\Lambda} \mu[N_{[0,1]^d}] \sup_{\eta \in \Omega} \mathbb{P}(X_{t, \Lambda}\Ssup{\widehat{\Lambda}, \eta} \neq X_{t, \Lambda}\Ssup{\eta}) 
        + \big( (\abs{\Lambda}\Vert b\Vert_\infty t)^2 + \abs{\Lambda}\Vert b\Vert_\infty t \big)^{1/2} \big(\sup_{\eta \in \Omega} \mathbb{P}(X_{t, \Lambda}\Ssup{\widehat{\Lambda}, \eta} \neq X_{t, \Lambda}\Ssup{\eta})\big)^{1/2}.
    \end{align*} 
    Moreover, from \zcref{Lemma:comparison_finite_infinite},
    \begin{align*}
        \sup_{\eta \in \Omega} \mathbb{P}(X_{t, \Lambda}\Ssup{\widehat{\Lambda}, \eta} \neq X_{t, \Lambda}\Ssup{\eta})
        &\lesssim \abs{\Lambda} \e^{-\mathrm{dist}(\Lambda, \widehat{\Lambda}^{c})\rho(\mathrm{dist}(\Lambda, \widehat{\Lambda}^{c}))}
        \leq \abs{\Lambda}^{- \rho(\log(\abs{\Lambda}) / C) / C + 1}
    \end{align*} 
    since \(\log(\vert\Lambda\vert) \leq C \, \mathrm{dist}(\Lambda, \widehat{\Lambda}^{c}) \), for some $C>0$, by assumption.
    We can then conclude that
    \begin{align*}
        \frac{1}{\abs{\Lambda}}\big\vert \mu T_t[f_\Lambda] - \mu T_t\Ssup{\widehat{\Lambda}}[f_\Lambda] \big\vert
        \lesssim \abs{\Lambda}^{- \rho(\log(\abs{\Lambda}) / C) / (2C) + 1/2}
        \xrightarrow[\Lambda, \widehat{\Lambda} \uparrow \R^d]{} 0.
    \end{align*}
    as desired. 
\end{proof}

\section{Characterizing reversible measures}\label{sec:reversible-Gibbs}

The main objective of this section is to prove the equivalence between reversible measures and infinite-volume Gibbs point processes described in \zcref{Prop:Reversible_Measures_equal_Gibbs_Measures}. 
The equivalences \textit{(ii)} \(\Leftrightarrow\) \textit{(iii)} \(\Leftrightarrow\) \textit{(iv)} are content of \zcref{proposition:equivalence_of_GNZ_Gibbs_infinitesimally_reversible} in the following \zcref{section:equivalence_of_infinitesimal_reversibility_and_GNZ}. The remaining equivalences are shown in \zcref{subsection:GNZ_to_reversible_to_infinitesimally_reversible}, in the form \textit{(i)} \(\Rightarrow\) \textit{(ii)} and \textit{(iv)} \(\Rightarrow\) \textit{(i)}, presented in \zcref{proposition:from_reversible_to_infinitesimally_reversible} and \zcref{proposition:from_GNZ_to_reversible}, respectively.

\subsection{Equivalence of infinitesimal reversibility and GNZ equations}\label{section:equivalence_of_infinitesimal_reversibility_and_GNZ}

We start with characterizing \(\GG_\theta\) via a symmetry of the formal generator \(\mathscr{L}\) and the identity \(b(o, \cdot) \mu = \mu_o^!\).
The proofs mainly follow the arguments in~\cite[Theorem]{Gloetzl1981}, which we report here for the convenience of the reader. 

Before we start, let us introduce the {\em reduced Campbell measure} \(C_\mu^{!}\) of $\mu\in\Pcal(\Omega)$, which is a measure on $\R^d \times \Omega$ defined by 
\begin{align*}
    C^!_\mu(B \times F) := \int_{\Omega}  \mu(d\eta) \int_{\R^d} \eta(dx)\, \1_B(x) \1_F(\eta - \delta_x), \quad B \in \mathfrak{B}(\R^d), F \in \mathscr{F}. 
\end{align*}
Its significance lies, among other things, in the following alternative representation of the GNZ equations, characterizing the Gibbs measures \(\GG_\theta\).
\begin{proposition}[{\cite[Proposition p.~219]{Gloetzl1981}}]\label{Prop:gloetzl-characterisation}
    For  $\mu \in \mathcal{P}_\theta(\Omega)$ we have \(\mu \in \GG_\theta\) if and only if
    \begin{align}\label{Equation:Reduced_Campbell_Measure_is_Birth_Rate}
        \frac{\d C_\mu^{!}}{\d x \otimes \d \mu}(x, \eta)
        = b(x, \eta),
        \quad \text{for } \d x \otimes \d \mu \text{ almost all } (x, \eta).
    \end{align} 
\end{proposition}

This in turn allows us to prove the following equivalences.
\begin{proposition}\label{proposition:equivalence_of_GNZ_Gibbs_infinitesimally_reversible}
    Let \(\mu \in \Pcal_\theta(\Omega)\). Any of the following statements is then equivalent to  \eqref{Equation:Reduced_Campbell_Measure_is_Birth_Rate}.
    \begin{enumerate}[(i)]
        \item \(\mu \in \GG_\theta\),
        \item \(b(o, \cdot) \mu = \mu_o^!\),
        \item \(\mu\) is infinitesimally reversible, i.e.\ \eqref{equation:definition_infinitesimally_reversible} holds for all \(f, g \in \Lcal\) such that \(\mathscr{L}f, \mathscr{L}g \in \Lcal\). 
        \item \eqref{equation:definition_infinitesimally_reversible} holds for all \(f, g \in \mathcal{A}\) for 
        \begin{align}\label{equation:small_test_function_class_infinitesimally_reversible}
            \mathcal{A}
            := \big\{{\1}_{F} \,\colon\,\exists B \Subset \R^d,  n \in \N_0 \text{ such that } F \in \mathcal{F}_B  \text{ and } F \subseteq \{N_B = n\} \big\}.
        \end{align}
    \end{enumerate}
\end{proposition}

Before we present the proof, let us introduce the following technical helper useful for inferring \eqref{Equation:Reduced_Campbell_Measure_is_Birth_Rate} from the infinitesimal reversibility condition \eqref{equation:definition_infinitesimally_reversible} with the test functions chosen from the more restricted test function class \(\mathcal{A}\) of \eqref{equation:small_test_function_class_infinitesimally_reversible}.
\begin{lemma}\label{le}
    Let \(\mu\in \Pcal(\Omega)\), \(B \Subset \R^d\), \(n \in \mathbb{N}_0\) and \(H, F \in \F\) local events with \(H \subseteq \{N_B = n\}\) and \(F \subseteq \{N_B = n+1\}\).
    
    Then,
    \begin{enumerate}[(i)]
        \item 
            \begin{align*}
                \mu\bigl[\1_H(\eta) (\LL \1_F)(\eta) \bigr]
                = \mu\Big[ \int_{B}\d x \,  b(x, \eta) \1_H(\eta) \1_F(\eta +\delta_x) \Big] \text{ and }
            \end{align*}

        \item 
            \begin{align*}
                \mu\bigl[\1_F(\eta) (\LL \1_H)(\eta) \bigr]
                = \int C_\mu^{!}(\d x, \d \eta)\,  \1_B(x) \1_F(\eta +\delta_x) \1_H(\eta).
            \end{align*}
    \end{enumerate}
\end{lemma}
\begin{proof}
    \textit{(i)}: By definition of the generator \(\LL\) and because of the conditions on $H$ and $F$, we have 
    \begin{align*}
        &\mu\big[\1_H(\eta) (\LL \1_F)(\eta) \big]\\
        &=
        \mu\Big[\1_H(\eta) \Big(\int_{\R^d}\d x \, b(x,\eta)\big(\1_F(\eta + \delta_x)-\1_F(\eta)\big) + \int\eta(\d x) \big(\1_F(\eta - \delta_x) - \1_F(\eta)\big) \Big)\Big]
        \\
        &=
        \mu\Big[
            \int_B \d x \, b(x,\eta) \1_{H}(\eta)\1_{F}(\eta+\delta_x) 
        \Big].
    \end{align*}
    \textit{(ii)}: Note that the roles of $F,H$ on the left-hand side are now reversed. Hence, the assumptions on the sets $F,H$ imply that 
    \begin{align*}
        \mu\big[\1_F(\eta) (\LL \1_H)(\eta) \big]
        &=
        \mu\Big[
        \1_F(\eta)\int_B \eta(\d x) \1_{H}(\eta-\delta_x)
        \Big]\\
        &=
        \int_\Omega \mu(d\eta) \int\eta(dx)\, \1_B(x) \1_F(\eta)\1_H(\eta -\delta_x).
    \end{align*}
    The claimed identity now follows from the definition of the reduced Campbell measure. 
\end{proof}

We are now in a position to give a
\begin{proof}[Proof of \zcref{proposition:equivalence_of_GNZ_Gibbs_infinitesimally_reversible}]

    \textit{(i) \(\Leftrightarrow\) \eqref{Equation:Reduced_Campbell_Measure_is_Birth_Rate}}: This equivalence is the content of \zcref{Prop:gloetzl-characterisation}.

    \textit{(ii) \(\Leftrightarrow\) \eqref{Equation:Reduced_Campbell_Measure_is_Birth_Rate}}: First assume that \(b(o, \cdot) \mu = \mu_o^!\). Then, by definition, for measurable \(f \geq 0\),
    \begin{align*}
        \int C_\mu^{!}(\d x, \d \eta)\, f(x, \eta)
        &= \int \mu(\d\eta) \, \int \eta(\d x) f(x, \eta-\delta_x)
        = \int \d x \int \mu_o^!(\d\eta) \, f(x, \theta_{-x} \eta) \\
        &= \int \d x \int \mu(\d\eta) \, b(o, \eta) f(x, \theta_{-x} \eta)
        = \int \d x \int \mu(\d\eta) \, b(o, \theta_{x}\eta) f(x, \eta) \\
        &= \int \d x \int \mu(\d\eta) \, b(x, \eta) f(x, \eta),
    \end{align*} 
    yielding \eqref{Equation:Reduced_Campbell_Measure_is_Birth_Rate}.
    For the other direction, assume that \eqref{Equation:Reduced_Campbell_Measure_is_Birth_Rate} holds.
    Let \(\varphi\) be some positive measurable function with \(\int \d x\, \varphi(x) = 1\). Then, for every measurable \(f \geq 0\),
        \begin{align*}
            \int\mu_o^!(\d\eta) \, f(\eta)
            &= \int \mu(\d\eta) \int \eta(\d x)\, \varphi(x) f(\theta_x (\eta -\delta_x))
            = \int C_\mu^{!}(\d x, \d \eta) \, \varphi(x) f(\theta_x \eta) \\
            &= \int \d x \int \mu(\d\eta) \, b(x, \eta) \varphi(x) f(\theta_x \eta)
            = \int \d x  \, \varphi(x) \int \mu(\d\eta) \, b(o, \theta_x\eta) f(\theta_x \eta) \\
            &= \int \mu(\d\eta) \, b(o,\eta) f(\eta),
        \end{align*} 
        yielding \(b(o, \cdot)\mu = \mu_o^!\).

    \textit{(iii) \(\Leftrightarrow\) \eqref{Equation:Reduced_Campbell_Measure_is_Birth_Rate}}: 
    Suppose first that \eqref{Equation:Reduced_Campbell_Measure_is_Birth_Rate} holds and let \(f, g \in \Lcal\). Applying \eqref{Equation:Reduced_Campbell_Measure_is_Birth_Rate} yields
    \begin{align*}
        &\mu[(\mathscr{L}f) g]
        = \int \mu(\d\eta) \Big[\int_{\R^d} \, b(x, \eta) [f(\eta + \delta_x) - f(\eta)] g(\eta) \, \d x  + \int \eta(\d x) [f(\eta - \delta_x) - f(\eta)] g(\eta) \Big] \\
        &= \int (b(\d x \otimes \d\mu)(\d x, \d\eta) [f(\eta + \delta_x) - f(\eta)] g(\eta) 
        + \int C_\mu^!(\d x, \d\eta) [f(\eta) - f(\eta + \delta_x)] g(\eta +\delta_x)  \\
        &= -\int (b(\d x \otimes \d\mu)(\d x, \d\eta) \, [f(\eta + \delta_x) - f(\eta)][g(\eta + \delta_x) - g(\eta) ],
    \end{align*} which is symmetric in \(f, g\) so that  \(\mu[(\mathscr{L}f) g] = \mu[f (\mathscr{L}g)]\). Integrability of all terms is justified by the boundedness of \(f, g\) and noticing that \eqref{Equation:Reduced_Campbell_Measure_is_Birth_Rate} implies \(\mu[N_B] \leq \Vert b\Vert_\infty \abs{B}\) for every \(B \Subset \R^d\).
    As \textit{(iii)} is evidently a priori stronger than \textit{(iv)}, it is enough to show in the following that \textit{(iv)} implies \eqref{Equation:Reduced_Campbell_Measure_is_Birth_Rate}.

    \textit{(iv) \(\Leftrightarrow\) \eqref{Equation:Reduced_Campbell_Measure_is_Birth_Rate}}:    
    We have already shown \eqref{Equation:Reduced_Campbell_Measure_is_Birth_Rate} \(\Rightarrow\) \textit{(iii)} \(\Rightarrow\) \textit{(iv)}.
    For the other direction, assume \textit{(iv)}, i.e.\ that \eqref{equation:definition_infinitesimally_reversible} holds for all \(f, g \in \mathcal{A}\).

    Now let \(B \Subset \R^d\), \(H \in \mathcal{F}\) be a \(D\)-local event for \(B \subseteq D \Subset \R^d\). Put
    \begin{align*}
        H_{n, k} := H \cap \{N_B = n\} \cap \{N_D = k\}, \quad
        H_{n, k}^+ = \{\eta \,\colon\, \eta = \zeta + \delta_x, \zeta \in H_{n, k}, x \in B\}.
    \end{align*}
    Indeed, we have then \({\1}_{H_{n, k}}, {\1}_{H_{n, k}^+} \in \mathcal{A}\), 
    so that using \zcref{le} and our assumption \eqref{equation:definition_infinitesimally_reversible} for all \(f, g \in \mathcal{A}\) twice yields    
    \begin{align*}
        C_\mu^{!}(B \times H_{n, k})
        &= \int C^!_\mu(\d x, \d \eta)\, \1_{B}(x) \1_{H_{n, k}}(\eta) \,  \\
        &= \int C^!_\mu(\d x, \d \eta)\, \1_{B}(x) \1_{H_{n, k}^{+}}(\eta +\delta_x) \1_{H_{n, k}}(\eta) 
        = \int \mu(\d\eta)\, \1_{H_{n, k}^{+}}(\eta) (\LL \1_{H_{n, k}})(\eta)  \\
        &= \int \mu(\d\eta)\, \1_{H_{n, k}}(\eta) (\LL \1_{H_{n, k}^{+}})(\eta) 
        = \int \mu(\d\eta) \int_{B} \d x\, b(x, \eta) \1_{H_{n, k}}(\eta) \1_{H_{n, k}^{+}}(\eta +\delta_x) \\
        &= (b(\d x \otimes \mu))(B \times H_{n, k}).
    \end{align*} 
    Summing over \(n, k \in \N_0\) proves that the two measures \(C_\mu^{!}\) and \(b(\d x \otimes \mu)\) agree on the \(\pi\)-system \(\{B \times H \,\colon\, B \Subset \R^d, H \in \mathcal{F} \text{ local}\}\) generating the whole \(\sigma\)-algebra.
    The sets \(\Lambda_n \times \Omega\) belong to this \(\pi\)-system and increase to the whole space \(\R^d \times \Omega\), and on these sets
    \begin{align*}
        C_\mu^{!}(\Lambda_n \times \Omega)
        = (b(\d x \otimes \d\mu))(\Lambda_n \times \Omega)
        \leq \Vert b \Vert_\infty \abs{\Lambda_n} 
        < \infty.
    \end{align*} Hence, \(C_\mu^{!} = b(\d x \otimes \d\mu)\)
\end{proof}

\subsection{From GNZ equations to reversibility to infinitesimal reversibility}\label{subsection:GNZ_to_reversible_to_infinitesimally_reversible}

\begin{proposition}\label{proposition:from_reversible_to_infinitesimally_reversible}
    If \(\mu \in \RR_\theta\) then it is infinitesimally reversible.
\end{proposition}
\begin{proof}
    By \zcref{proposition:equivalence_of_GNZ_Gibbs_infinitesimally_reversible}, it is enough to establish \eqref{equation:definition_infinitesimally_reversible} for all \(f, g \in \mathcal{A}\) with \(\mathcal{A}\) defined in \eqref{equation:small_test_function_class_infinitesimally_reversible}.
    Let \(f, g \in \mathcal{A}\).
    We compute, using the reversibility of \(\mu\), 
    \begin{align*}
        \mu[(\mathscr{L} f) g]
        &= \mu\Big[\Big(\lim_{t \downarrow 0} \frac{T_t f - f}{t} \Big) g \Big]
        = \lim_{t \downarrow 0}  \mu\Big[\frac{T_t f - f}{t} g \Big]
        = \lim_{t \downarrow 0}  \mu\Big[f \frac{T_t g - g}{t} \Big] \\
        &= \mu\Big[f  \Big(\lim_{t \downarrow 0} \frac{T_t g - g}{t} \Big) \Big]
        =  \mu[f (\mathscr{L} g)].
    \end{align*}
    Indeed, \(\mathscr{L} h = \lim_{t \downarrow 0} \frac{T_t h - h}{t}\) pointwise for \(h \in \mathcal{A} \subseteq \Lcal\) by \zcref{proposition:transition_semigroup}, \zcref{lemma:infinite_volume_process_forward_equation_local_bounded_measurable} and, as the following argument shows, we can apply dominated convergence, since for any \(h \in \mathcal{A}\), 
    \begin{align*}
        \sup_{t > 0} \Big\Vert \frac{T_t h - h}{t} \Big\Vert_\infty < \infty.
    \end{align*}
    Recall from \zcref{proposition:transition_semigroup}, \zcref{lemma:infinite_volume_process_forward_equation_local_bounded_measurable}, 
    \begin{align*}
        \frac{T_t h - h}{t}
        = \frac{1}{t} \int_0^t T_s \mathscr{L} h  \, \d s,
    \end{align*} so that by the triangle inequality
    \begin{align*}
        \Big\Vert \frac{T_t h - h}{t} \Big\Vert_\infty
        \leq \Vert \mathscr{L} h \Vert_\infty.
    \end{align*}
    Since \(h \in \mathcal{A}\), there are \(n \in \N_0\) and \(B \Subset \R^d\) such that \(h\) is \(B\)-local and \(\abs{h} \leq {\1}_{N_B = n}\), yielding
    \begin{align*}
         \Vert \mathscr{L} h \Vert_\infty
         \leq (\abs{B} \Vert b\Vert_\infty + n) \lor  n+1 
         < \infty.
    \end{align*}
\end{proof}

\begin{proposition}\label{proposition:from_GNZ_to_reversible}
    We have that \(\GG_\theta\subset \RR_\theta\).
\end{proposition}
\begin{proof}
    Let $\mu\in\GG_\theta$. From \zcref{proposition:finite_volume_reversibility}, we know that the GNZ equations show that, for the corresponding finite-volume dynamics \((\Tl_t)_{t \geq 0}\), the reversibility
    \begin{align*}
        \mu[(\Tl_t f) g]
        = \mu[f (\Tl_t g)]
    \end{align*} 
    holds for all bounded, measurable \(\Lambda\)-local \(f, g\).
    For fixed bounded, measurable local \(f, g\), we therefore get, thanks to the finite speed of propagation result of \zcref{Lemma:finite_speed_of_propagation_version_1}, that
    \begin{align*}
        \mu[(T_t f) g]
        = \lim_{\Lambda \uparrow \R^d} \mu[(\Tl_t f) g]
        =  \lim_{\Lambda \uparrow \R^d} \mu[f (\Tl_t g)]
        = \mu[f (T_t g)],
    \end{align*} 
    meaning that \(\mu \in \RR_\theta\).
\end{proof}

\section{Monotonicity of the specific relative entropy}\label{sec:relative-entropy-decreasing}

    In this section, we complete the proof of the monotonicity of the specific relative entropy along trajectories of the infinite volume dynamics. In \zcref{Sec:Quasi_Superadditivity_of_Entropy_wrt_Gibbs_measure}, we derive a factorization estimate for Gibbs measures, cf.~\zcref{Le:factor_property_of_the_Gibbs_measure}, which yields the quasi-superadditivity of the relative entropy of~\zcref{Le:quasi_superadditivity_entropy_wrt_gibbs_measure}, as well as the existence and lower semi-continuity of the specific relative entropy, cf.~\zcref{prop:existence-relative-entropy}. We then establish a technical convergence result, \zcref{lemma:tau_L_convergence_from_weak_convergence_and_bounded_specific_entropy}, that will be used in the proof of the attractor property. Finally, we use the finite speed of propagation of~\zcref{Lemma:comparison_finite_infinite} to compare the relative entropies under the finite- and infinite-volume dynamics, cf.~\zcref{Le:compariability_of_finite_volume_entropies_original_vs_modified}, and derive, in~\zcref{lemma:control-number-of-points}, bounds on the expected number of points in bounded volumes, needed to control the resulting boundary errors.
    \subsection{Quasi-superadditivity of the relative entropy}\label{Sec:Quasi_Superadditivity_of_Entropy_wrt_Gibbs_measure}
    
        It will be useful to see that for large boxes \(\Lambda\), any Gibbs measure \(\nu \in \Gcal_\theta\) is not too far from the measure \(\nu_{\Lambda} \otimes \nu_{\Lambda^c}\), where the particles in \(\Lambda\) and \(\Lambda^c\) are independently distributed according to the correct marginal of \(\nu\). 

       \begin{lemma}[Factor property]\label{Le:factor_property_of_the_Gibbs_measure}
            For every \(\nu \in \Gcal_\theta\) and \(\Lambda \Subset \R^d\), \(\nu\) is absolutely continuous with respect to \(\nu_{\Lambda} \otimes \nu_{\Lambda^c}\) and we have the bound
            \begin{align}\label{equation:gibbs_radon_nikodym_wrt_independent_inside_outside_bound}
                \frac{\d\nu}{\d(\nu_{\Lambda} \otimes \nu_{\Lambda^c})}
                \leq \exp\big(\log(\Vert b\Vert_\infty \Vert 1/b\Vert_\infty) N_{\partial^R \Lambda} + (\Vert b\Vert_\infty - 1/\Vert 1/b\Vert_\infty) \abs{\partial^R \Lambda}\big)
            \end{align} 
            for the Radon--Nikodym derivative of \(\nu\) w.r.t.\ \(\nu_{\Lambda} \otimes \nu_{\Lambda^c}\), with \(\partial^R \Lambda = \{x \in \Lambda \,\colon\, \mathrm{dist}(x, \Lambda^c) \leq R\}\).
        \end{lemma}
        \begin{proof}
            An easy calculation using the DLR equations gives
            \begin{align}\label{equation:gibbs_radon_nikodym_wrt_independent_inside_outside}
                \frac{\d\nu}{\d(\nu_{\Lambda} \otimes \nu_{\Lambda^c})}(\omega_\Lambda \omega_{\Lambda^c})
                = \Big(\int \nu_{\Lambda^c}(\d\zeta_{\Lambda^c}) \, \frac{Z_{\Lambda, \omega_{\Lambda^c}}}{Z_{\Lambda, \zeta_{\Lambda^c}}} \e^{H_{\Lambda, \omega_{\Lambda^c}}(\omega_{\Lambda}) - H_{\Lambda, \zeta_{\Lambda^c}}(\omega_{\Lambda})} \Big)^{-1}
            \end{align} 
            Indeed, for measurable \(f \geq 0\), by the DLR equations,
            \begin{align*}
                &\int (\nu_{\Lambda} \otimes \nu_{\Lambda^c})(\d\omega) \, f(\omega)
                = \int \nu(\d\widetilde{\omega}) \int \nu(\d\omega) \, f(\omega_\Lambda \widetilde{\omega}_{\Lambda^c})  \\
                &= \int \nu(\d\widetilde{\omega}) \int \nu_{\Lambda^c}(\d\zeta_{\Lambda^c}) \int \pi_\Lambda(\d\omega_\Lambda) \, \frac{1}{Z_{\Lambda, \zeta_{\Lambda^c}}} \e^{-H_{\Lambda, \zeta_{\Lambda^c}}(\omega_\Lambda)} \, f(\omega_\Lambda \widetilde{\omega}_{\Lambda^c})  \\
                &= \int \nu_{\Lambda^c}(\d\omega_{\Lambda^c}) \int \pi_\Lambda(\d\omega_\Lambda) \, \Big[\int \nu_{\Lambda^c}(\d\zeta_{\Lambda^c}) \frac{1}{Z_{\Lambda, \zeta_{\Lambda^c}}} \e^{-H_{\Lambda, \zeta_{\Lambda^c}}(\omega_\Lambda)}\Big] \, f(\omega_\Lambda \omega_{\Lambda^c}) 
            \end{align*} 
            and, also by the DLR equations,
            \begin{align*}
                \int \nu(\d\omega) \, f(\omega)
                = \int \nu_{\Lambda^c}(\d\omega_{\Lambda^c}) \int \pi_\Lambda(\d\omega_\Lambda) \, \Big[\frac{1}{Z_{\Lambda, \omega_{\Lambda^c}}}  \e^{-H_{\Lambda, \omega_{\Lambda^c}}(\omega_\Lambda)}\Big] \, f(\omega_\Lambda \omega_{\Lambda^c}). 
            \end{align*}
            Let us now show that \eqref{equation:gibbs_radon_nikodym_wrt_independent_inside_outside_bound} holds. For that, we bound the integrand in \eqref{equation:gibbs_radon_nikodym_wrt_independent_inside_outside} from below as follows. 
            First,
            \begin{align}\label{equation:uniform_lower_bound_energy_difference}
                H_{\Lambda, \omega_{\Lambda^c}}(\omega_{\Lambda}) - H_{\Lambda, \zeta_{\Lambda^c}}(\omega_{\Lambda})
                \geq  -\log(\Vert b\Vert_\infty \Vert 1/b\Vert_\infty) N_{\partial^R \Lambda}(\omega_\Lambda),
            \end{align} 
            by successively inserting particles of \(\omega_\Lambda\) into the configuration \(\emptyset_\Lambda \omega_{\Lambda^c}\) resp. \(\emptyset_\Lambda \zeta_{\Lambda^c}\), starting with the particles in \(\Lambda \setminus \partial^R \Lambda\), where the corresponding conditional energies of adding the points agree for the two boundary conditions \(\omega_{\Lambda^c}, \zeta_{\Lambda^c}\) by finite range of \(b\). To compare the partition functions in the integrand, the same argument of successive insertions gives
            \begin{align}\label{equation:uniform_bounds_boltzmann_factor}
                \Vert 1 / b \Vert_\infty^{-N_{\partial^R \Lambda}(\omega_{\Lambda})} \e^{-H_{\Lambda \setminus \partial^R \Lambda}(\omega_{\Lambda \setminus \partial^R \Lambda})}
                \leq \e^{-H_{\Lambda, \xi_{\Lambda^c}}(\omega_\Lambda)}
                \leq \Vert b \Vert_\infty^{N_{\partial^R \Lambda}(\omega_{\Lambda})} \e^{-H_{\Lambda \setminus \partial^R \Lambda}(\omega_{\Lambda \setminus \partial^R \Lambda})}
            \end{align} 
            for every boundary condition \(\xi_{\Lambda^c}\). Integrating \zcref{equation:uniform_bounds_boltzmann_factor} against the Poisson reference measure \(\pi_\Lambda = \pi_{\Lambda \setminus \partial^R \Lambda} \otimes \pi_{\partial^R \Lambda}\) gives
            \begin{align}\label{equation:uniform_lower_bound_partition_function_ratio}
                \frac{Z_{\Lambda, \omega_{\Lambda^c}}}{Z_{\Lambda, \zeta_{\Lambda^c}}} 
                \geq \exp\big(-(\Vert b\Vert_\infty - 1/\Vert 1/b\Vert_\infty) \abs{\partial^R \Lambda}\big).
            \end{align}
            Finally, substituting \eqref{equation:uniform_lower_bound_partition_function_ratio} and \eqref{equation:uniform_lower_bound_energy_difference} into \eqref{equation:gibbs_radon_nikodym_wrt_independent_inside_outside} yields \eqref{equation:gibbs_radon_nikodym_wrt_independent_inside_outside_bound}. 
        \end{proof}
        
        This weak dependence under \(\nu\) between regions \(\Lambda\) and \(\Lambda^c\) translates into the quasi-super\-additivity of the relative entropy claimed in~\zcref{Le:quasi_superadditivity_entropy_wrt_gibbs_measure}. 
        The proof is an adaption of~\cite[Lemma 4.1]{Sullivan1976} to our setting.
        \begin{proof}[Proof of~\zcref{Le:quasi_superadditivity_entropy_wrt_gibbs_measure}]
            Of course, if \(\RelEnt_{\Lambda \cup \Delta}(\mu |\nu) = \infty\), the statement is trivial, so we can assume it to be finite, which means there exists a density \(\d \mu_{\Lambda \cup \Delta} / \d \nu_{\Lambda \cup \Delta}\), and that \(\RelEnt_{\Lambda}(\mu |\nu)\), \(\RelEnt_{\Delta}(\mu |\nu)\) are finite too.
            The Radon--Nikodym derivative \(\d\nu_{\Lambda \cup \Delta}/\d(\nu_\L\otimes\nu_\Delta)\) is obtained from \(\d\nu / \d(\nu_{\Lambda} \otimes \nu_{\Lambda^c})\) by conditional expectation,
            \begin{align*}
                \frac{\d \nu_{\Lambda \cup \Delta}}{\d(\nu_\Lambda \otimes \nu_\Delta)}
                = (\nu_\Lambda \otimes \nu_{\Lambda^c})\Big[\frac{\d\nu}{\d(\nu_{\Lambda} \otimes \nu_{\Lambda^c})} \,\big\vert\, \mathcal{F}_{\Lambda \cup \Delta}\Big].
            \end{align*} 
            The chain rule for Radon--Nikodym derivatives gives
            \begin{align*}
                \frac{\d \mu_{\Lambda \cup \Delta}}{\d(\nu_{\Lambda} \otimes \nu_{\Delta})}
                = \frac{\d \nu_{\Lambda \cup \Delta}}{\d(\nu_\Lambda \otimes \nu_\Delta)} \frac{\d \mu_{\Lambda \cup \Delta}}{\d \nu_{\Lambda \cup \Delta}}.
            \end{align*}
            In particular, by \zcref{Le:factor_property_of_the_Gibbs_measure} and because the upper bound in it is \(\mathcal{F}_\Lambda\)-measurable, we get
            \begin{align*}
                \frac{\d \mu_{\Lambda \cup \Delta}}{\d(\nu_{\Lambda} \otimes \nu_{\Delta})}
                \leq \e^{\log(\Vert b\Vert_\infty \Vert 1/b\Vert_\infty) N_{\partial^R \Lambda} + (\Vert b\Vert_\infty - 1/\Vert 1/b\Vert_\infty) \abs{\partial^R \Lambda}} \frac{\d \mu_{\Lambda \cup \Delta}}{\d \nu_{\Lambda \cup \Delta}}.
            \end{align*}
            It follows, by taking logarithms and the \(\mu\)-expectation, and the assumption \(\RelEnt_{\Lambda \cup \Delta}(\mu |\nu) < \infty\), that \(\RelEnt(\mu_{\Lambda \cup \Delta} |\nu_\Lambda \otimes \nu_\Delta)\) is finite too and that
            \begin{align*}
                \RelEnt_{\Lambda \cup \Delta}(\mu |\nu)
                \geq \RelEnt(\mu_{\Lambda \cup \Delta} |\nu_\Lambda \otimes \nu_\Delta) - \big\{ \log(\Vert b\Vert_\infty \Vert 1/b\Vert_\infty) \mu[N_{\partial^R \Lambda}] + (\Vert b\Vert_\infty - 1/\Vert 1/b\Vert_\infty) \abs{\partial^R \Lambda} \big\}.
            \end{align*} 
            This yields the desired estimate, since
            \begin{align*}
                \RelEnt(\mu_{\Lambda \cup \Delta} |\nu_\Lambda \otimes \nu_\Delta)
                - \Big( \RelEnt_{\Lambda}(\mu |\nu) + \RelEnt_{\Delta}(\mu |\nu) \Big)
                = \RelEnt(\mu_{\Lambda \cup \Delta} |\mu_\Lambda \otimes \mu_\Delta)
                \geq 0,
            \end{align*}
            which finishes the proof.
        \end{proof}

        As a by-product, we obtain the existence (and lower semi-continuity) of the specific relative entropy with respect to $\nu$. Let us note that in the continuum this existence usually requires quite some work, see e.g.~\cite[Section 4.4]{JKSZ24}, whereas the proof given here is very short and elementary. We also note that, in our setting, existence can also be deduced from~\cite{Dereudre_2016_variationalprinciple}.        

        \begin{proof}[Proof of~\zcref{prop:existence-relative-entropy}]        
            \textit{(i) Existence of the limit: }
            Let us first record that translation invariance, together with the entropy inequality \zcref{equation:entropic_inequality} and \eqref{equation:gibbs_marginal_exponential_moment}, yields
            \begin{equation}\label{equation:first_moment_bound_by_finite_entropy}
            \begin{split}
                \mu[N_{[0, 1)^d}]
                &= \liminf_{n \to \infty }\frac{\mu[N_{\Lambda_n}]}{\abs{\Lambda_n}} 
                \leq \liminf_{n \to \infty} \Big( \frac{\RelEnt_{\Lambda_n}(\mu \vert \nu)}{\abs{\Lambda_n}} + \frac{\log \nu_{\Lambda_n}[\e^{N_{\Lambda_n}}]}{\abs{\Lambda_n}} \Big)\\
                &\leq \SpecEnt(\mu \vert \nu) + (\e - 1)\Vert b\Vert_\infty.
            \end{split}
            \end{equation}
             If \(\SpecEnt(\mu\vert\nu) = \infty\), the claimed convergence \eqref{equation:existence-relative-entropy-limit} is a direct consequence of the definition as limes inferior. Hence, we can assume \(\SpecEnt(\mu\vert\nu) < \infty\) now,  in particular, \(\RelEnt_{\Lambda_n}(\mu \vert \nu) < \infty\) and \(\mu[N_{[0, 1)^d}] < \infty\) by \eqref{equation:first_moment_bound_by_finite_entropy}.
            Now fix \(n \in \N\) and let us show
            \begin{align}\label{equation:lower_bound_spec_ent_fixed_volume_error_term}
                \SpecEnt(\mu|\nu)
                \geq 
                \frac{\RelEnt_{\Lambda_n}(\mu \vert \nu)}{\abs{\Lambda_n}}
                - \frac{\abs{\partial^R \Lambda_n}}{\abs{\Lambda_n}} \big\{ \log(\Vert b\Vert_\infty \Vert 1/b\Vert_\infty) \mu[N_{[0,1)^d}] + (\Vert b\Vert_\infty - 1/\Vert 1/b\Vert_\infty) \big\}.
            \end{align} Taking the \(\limsup_{n \to \infty}\) then yields the limit \eqref{equation:existence-relative-entropy-limit} by the definition of \(\SpecEnt(\mu\vert\nu)\) as corresponding \(\liminf_{n \to \infty}\).
            Decompose the box $\Lambda_m$, \(m \geq n\), into 
            \begin{align*}
                \Lambda_m = \Lambda^1 \cup \Lambda^2 \cup \dots \cup \Lambda^{\ell_m} \cup \Delta_m, 
            \end{align*}
            where $\Lambda^i$, $i=1,\dots,\ell_m$ are disjoint translates of $\Lambda_n$ and $\Delta_m \subset \Lambda_m$ is the possibly empty remainder. These decompositions can be arranged such that
            \begin{equation}\label{partition-property} 
            \begin{split}
                \frac{\abs{\Delta_m}}{\abs{\Lambda_m}} \xrightarrow[m \to \infty]{} 0,
                \quad \frac{\ell_m \abs{\Lambda_n}}{\abs{\Lambda_m}} \xrightarrow[m \to \infty]{} 1.
            \end{split}
            \end{equation}
            We can now repeatedly apply the estimate from \zcref{Le:quasi_superadditivity_entropy_wrt_gibbs_measure} to such a partition, and use the translation invariance of $\mu$ and $\nu$ to obtain
            \begin{align*}
                &I_{\Lambda_m}(\mu | \nu) \\
                &\geq \ell_m I_{\Lambda_n}(\mu | \nu) 
                - \ell_m \abs{\partial^R \Lambda_n} \big\{ \log(\Vert b\Vert_\infty \Vert 1/b\Vert_\infty) \mu[N_{[0,1)^d}] + (\Vert b\Vert_\infty - 1/\Vert 1/b\Vert_\infty) \big\}
            \end{align*}
            Dividing by \(\abs{\Lambda_m}\) and taking the \(\liminf_{m \to \infty}\) yields then, recalling \eqref{partition-property}, the sufficient \eqref{equation:lower_bound_spec_ent_fixed_volume_error_term}.

            \textit{(ii) Lower semicontinuity: } 
            First note that for every finite box $\Lambda_n \Subset \R^d$ and every \(\mu \in \Pcal_\theta(\Omega)\), by translation invariance we have that
            \begin{align*}
                \mu(N_{\partial \Lambda_n})
                = 0.
            \end{align*} 
            Consequently, if \(\mu_k \Rightarrow \mu\) for a sequence \((\mu_k)_{k \in \N}\) in \(\Pcal_\theta(\Omega)\), the restrictions to \(\Lambda_n\) converge weakly too, \((\mu_k)_{\Lambda_n} \Rightarrow \mu_{\Lambda_n}\).
            The Donsker--Varadhan variational formula \eqref{equation:donsker_varadhan_variational_formula_probability_measures}, replacing the class of test functions with the bounded continuous functions, implies that the functional $\alpha \mapsto I(\alpha | \nu_{\Lambda_n})$ is lower semicontinuous with respect to the weak topology on \(\Pcal(\Omega_{\Lambda_ n})\). Combining these two facts,
            \begin{align}\label{equation:local_entropy_lower_semicontinuity}
                \RelEnt_{\Lambda_n}(\mu \vert \nu)
                \leq \liminf_{k \to \infty} \RelEnt_{\Lambda_n}(\mu_k \vert \nu).
            \end{align}

            If \(\liminf_{k \to \infty}\mathscr{I}(\mu_k | \nu)  = \infty\), then there is nothing to show. Hence, we can assume that there is a subsequence \((\mu_{k_j})_{j \in \N}\) with bounded specific relative entropies such that 
            \begin{align*}
                \lim_{j \to \infty}\mathscr{I}(\mu_{k_j} | \nu)
                =\liminf_{k \to \infty}\mathscr{I}(\mu_k | \nu)
                < \infty.
            \end{align*}
            In particular, \(C := \sup_{j \in \N} \mu_{k_j}[N_{[0, 1)^d}] < \infty\) by \zcref{equation:first_moment_bound_by_finite_entropy}.
            We can therefore again take \eqref{equation:lower_bound_spec_ent_fixed_volume_error_term} as a starting point and see that for every \(j \in \N\) and $n \in \N$ we have 
            \begin{align*}
                &\mathscr{I}(\mu_{k_j} | \nu)
                \\
                &\geq \frac{I_{\Lambda_n}(\mu_{k_j} | \nu)}{\abs{\Lambda_n}} - \frac{\abs{\partial^R \Lambda_n}}{\abs{\Lambda_n}} \big\{ \log(\Vert b\Vert_\infty \Vert 1/b\Vert_\infty) C + (\Vert b\Vert_\infty - 1/\Vert 1/b\Vert_\infty) \big\}. 
            \end{align*}
            The claimed lower semicontinuity now follows from first taking the limit \(j \to \infty\), applying \eqref{equation:local_entropy_lower_semicontinuity}, and then taking \(n \to \infty\).
        \end{proof}

        \subsection{Upgrading the mode of convergence}
        We now provide the proof of an important technical helper for proving the attractor property. 
        
        \begin{proof}[Proof of~\zcref{lemma:tau_L_convergence_from_weak_convergence_and_bounded_specific_entropy}]
            The lemma is inspired by \cite[Proposition 2.6]{GeorgiiZessin1993}, and one could also adapt its proof to show that the entropy level sets \(\{\SpecEnt(\cdot \vert \nu) \leq c\}\), \(c \geq 0\), are sequentially compact w.r.t.\ convergence of expectations of the function class of all local, measurable \(f\) with \(\abs{f} \lesssim (1+N_\Lambda)\) for some \(\Lambda \Subset \R^d\), now called \emph{tame} functions.
            Then, weak convergence of a sequence \((\rho_k)_{k \in \N}\) with uniformly bounded specific relative entropy as above would also force this convergence by uniqueness of limits.
            We give a slightly more direct argument using the entropy inequality \zcref{equation:entropic_inequality} to first uniformly approximate tame functions \(f\) with bounded local measurable functions, and then further with bounded continuous functions.
            
            First fix \(f \in \Lcal\) and let \(\Lambda \Subset \R^d\) be a cube such that \(f\) is \(\Lambda\)-measurable, and also let \((f_n)_{n \in \N}\) be a sequence of bounded continuous functions on \(\Omega_\Lambda\) with
            \begin{align}\label{equation:approximation_by_continuous_functions_with_same_bounds}
                \Vert f_n\Vert_\infty \leq \Vert f\Vert_\infty,
                \quad f_n \xrightarrow[n \to \infty]{} f \text{ in } \nu_\Lambda\text{-probability}.
            \end{align} 
            Notice that this is possible via Lusin's theorem; \(\nu_\Lambda\) is a Borel probability measure on the Borel subset \(\Omega_\Lambda\) of the Polish space \(\Omega\) and therefore Radon.
            Now, by the entropy inequality \zcref{equation:entropic_inequality} we get
            \begin{align}\label{equation:uniform_approximation_by_continuous_functions_from_entropy_bounds}
                \sup_{k \in \N} \rho_k[\abs{f - f_n}]
                \leq \delta \sup_{k \in \N} I_{\Lambda}(\rho_k \,\vert\, \nu) + \delta \log \nu_\Lambda[\e^{\abs{f-f_n} / \delta}],
            \end{align} 
            for every \(\delta > 0\) and \(n \in \N\).
            The same upper bound holds for \(\rho[\abs{f - f_n}]\) by the lower semicontinuity of the local relative entropy \eqref{equation:local_entropy_lower_semicontinuity}. Observe that, by \eqref{equation:lower_bound_spec_ent_fixed_volume_error_term} combined with \eqref{equation:first_moment_bound_by_finite_entropy} and the assumption \eqref{equation:assumption_bounded_specific_entropy}, 
            \begin{align*}
                &\sup_{k \in \N} I_{\Lambda}(\rho_k \,\vert\, \nu) \\
                &\leq \abs{\Lambda} \sup_{k \in \N} \SpecEnt(\rho_k \vert \nu)  \\
                &\qquad + \abs{\partial^R \Lambda} \big\{ \log(\Vert b\Vert_\infty \Vert 1/b\Vert_\infty) [\sup_{k \in \N} \SpecEnt(\rho_k \vert \nu) + (\e -1)\Vert b\Vert_\infty] + (\Vert b\Vert_\infty - 1/\Vert 1/b\Vert_\infty) \big\} \\
                &< \infty.
            \end{align*}
            Hence, \eqref{equation:convergence_of_expectations_bounded_local_measurable} follows for \(f \in \Lcal\) by an \(\epsilon /3\)-argument, \eqref{equation:uniform_approximation_by_continuous_functions_from_entropy_bounds} and \zcref{equation:approximation_by_continuous_functions_with_same_bounds}, first letting \(k \to \infty\), then \(n \to \infty\) and finally \(\delta \downarrow 0\).

            Similarly, we upgrade the convergence to tame \(f\). Let again \(\Lambda \Subset \R^d\) be a cube such that \(f\) is \(\Lambda\)-measurable and such that \(\abs{f} \leq c(1+ N_\Lambda)\), \(c > 0\). Then, 
            \begin{align}\label{equation:uniform_integrability_particle_numbers_entropy_bound}
                \rho[N_\Lambda {\1}_{N_\Lambda > M}] \lor \big(\sup_{k \in \N} \rho_k[N_\Lambda {\1}_{N_\Lambda > M}]\big)
                \leq \delta \sup_{k \in \N} I_{\Lambda}(\rho_k \,\vert\, \nu) + \delta \log \nu_\Lambda[\e^{N_\Lambda {\1}_{N_\Lambda > M} / \delta}]
            \end{align} as above. 
            By the corresponding exponential moment bound with rate \(1/\delta\), completely analogous to \eqref{equation:gibbs_marginal_exponential_moment}, 
            \begin{align*}
                \nu_\Lambda[\e^{N_\Lambda {\1}_{N_\Lambda > M} / \delta}]
                \xrightarrow[M \to \infty]{} 1
            \end{align*} 
            for fixed \(\delta > 0\). 
            Again, the convergence \eqref{equation:convergence_of_expectations_bounded_local_measurable} follows from an \(\epsilon/3\) argument, applying the bounded case above to \(f {\1}_{N_\Lambda \leq M}\) and handling the error of this approximation by \eqref{equation:uniform_integrability_particle_numbers_entropy_bound}, letting first \(k \to \infty\), then \(M \to \infty\) and finally \(\delta \downarrow 0\).
        \end{proof}

    \subsection{Comparison of relative entropies under global and local dynamics}\label{Sec:compariability_of_finite_volume_entropies_original_vs_modified}
    
        \begin{proof}[Proof of~\zcref{Le:compariability_of_finite_volume_entropies_original_vs_modified}]
            By the assumption \(\RelEnt_{\Lambda}(\mu T_t | \nu) > r\) and the Donsker--Varadhan formula \eqref{equation:donsker_varadhan_variational_formula_probability_measures}, there is a \(\delta > 0\) and a bounded measurable \(\Lambda\)-local \(F\) such that
            \begin{align*}
                (\mu T_t)[F] - \log(\nu[\e^F])
                > r + \delta.
            \end{align*}
            Now let \(\Lambda^* \Subset \R^d\) with \(\Lambda \subseteq \Lambda^*\) and \(\widetilde{\Lambda} \subseteq \Lambda^*\) be a translate of \(\Lambda\).
            Set \(G = F \circ \theta_x\), where \(x \in \mathbb{R}^d\) is such that \(\Lambda = \theta_x\widetilde{\Lambda}\).
            We then have, again by the Donsker--Varadhan formula \eqref{equation:donsker_varadhan_variational_formula_probability_measures} and by translation invariance,
            \begin{align*}
                \RelEnt_{\widetilde{\Lambda}}\bigl(\mu T_t\Ssup{\Lambda^\ast} \big| \nu \bigr) 
                &\geq (\mu T_t\Ssup{\Lambda^\ast})[G] - \log(\nu[\exp(G)]) \\
                &\geq (\mu T_t)[G] - \log(\nu[\exp(G)]) - \big\vert \mu T_t\Ssup{\Lambda^\ast}[G] - \mu T_t[G] \big\vert \nonumber \\
                &= (\mu T_t)[F] - \log(\nu[\exp(F)]) - \big\vert \mu T_t\Ssup{\Lambda^\ast}[G] - \mu T_t[G] \big\vert \nonumber \\
                &> r + \delta - \big\vert \mu T_t[G] - \mu T_t\Ssup{\Lambda^\ast}[G] \big\vert. \nonumber
            \end{align*} 
            It remains to show that there is some \(\ell > 0\) which is independent of \(\Lambda^*\) and the concrete translate \(\widetilde{\Lambda} \subseteq \Lambda^*\) of \(\Lambda\), such that
            \begin{align}\label{equation:finite_volume_discrepancy_applied_to_test_function_entropies}
                \big\vert \mu T_t[G] - \mu T_t\Ssup{\Lambda^\ast}[G] \big\vert
                \leq \delta
            \end{align} if \(\mathrm{dist}(\widetilde{\Lambda}, (\Lambda^*)^c) \geq \ell\).
            Indeed, we can apply \zcref{Lemma:finite_speed_of_propagation_version_1} for that.
            For general \(\Lambda \Subset \R^d\), we first fix a box of side lengths at least \(1\) with \(\Lambda \subseteq Q\). Then, for every translate \(\widetilde{\Lambda}\) of \(\Lambda\), the accordingly translated box \(\widetilde{Q}\) has \(\mathrm{dist}(\widetilde{Q}, (\Lambda^*)^c) \geq \mathrm{dist}(\widetilde{\Lambda}, (\Lambda^*)^c) - \mathrm{diam}(Q)\). In particular, by \zcref{Lemma:finite_speed_of_propagation_version_1},
            \begin{align}\label{equation:finite_volume_discrepancy_applied_to_test_function_entropies_bound}
                \big\vert  \mu T_t[G] - \mu T_t\Ssup{\Lambda^\ast}[G]\big\vert
                \leq c_t \abs{Q} \Vert F\Vert_\infty \e^{- \mathrm{dist}(\widetilde{Q}, (\Lambda^*)^c) \rho(\mathrm{dist}(\widetilde{Q}, (\Lambda^*)^c))}.
            \end{align} 
            If we now take \(L\) such that, for \(l \geq L\), \(\rho(l) \geq 1\) and
            \begin{align}
                \ell
                = (\mathrm{diam}(Q) + L) \lor \log \frac{c_t \abs{Q} \e^{\mathrm{diam}(Q)} \Vert F\Vert_\infty}{\delta},
            \end{align}
            then the translates \(\widetilde{Q}\)
            lie inside \(\Lambda^*\) with positive distance to the boundary and the right-hand side of \eqref{equation:finite_volume_discrepancy_applied_to_test_function_entropies_bound} is smaller than \(\delta\), finally meaning that~\eqref{equation:finite_volume_discrepancy_applied_to_test_function_entropies} holds for \(\mathrm{dist}(\widetilde{\Lambda}, (\Lambda^*)^c) \geq \ell\), independent of \(\Lambda^*\) and the concrete translate \(\widetilde{\Lambda} \subseteq \Lambda^*\) of \(\Lambda\).
        \end{proof}

\subsection{Controlling the number of points}

\begin{proof}[Proof of~\zcref{lemma:control-number-of-points}]
        The graphical construction gives for every measurable \(B \subseteq \Lambda\) the bounds
        \begin{align}\label{equation:proof-bounds_particle_numbers_local_dynamics} 
            \e^{-t} \mu[N_B]
            \leq (\mu \Tl_t)[N_B]
            \leq \e^{-t}\mu[N_B] + \Vert b\Vert_\infty(1-\e^{-t}) \abs{B}
        \end{align} 
        via stochastic majorization by the surviving initial particles plus the contribution from
        \begin{align*}
            \int_{B \times (0, t] \times [0, \infty) \times [0, \Vert b\Vert_\infty]} \Ncal(\d x, \d s, \d r, \d u) {\1}_{s+r > t} \delta_{x},
        \end{align*} 
        i.e., the birth proposals surviving until time \(t\), and stochastic minorization by the surviving initial particles.
        The lower bound in \eqref{equation:proof-bounds_particle_numbers_local_dynamics} together with the entropy inequality \eqref{equation:entropic_inequality} and the exponential moment bound \eqref{equation:gibbs_marginal_exponential_moment} yield that \(\mu[N_{[0,1)^d}] = \infty\) already implies \(\RelEnt_{\Lambda_k}(\mu T_t\Ssup{\Lambda_k} |\nu) = \infty\) for all \(k \in \N\). 
\end{proof}

\section{Entropy dissipation and Fisher information}\label{sec:attractor-property}
\subsection{Finite-volume entropy-dissipation identity}\label{Sec:Calculations_in_Finite_Volume_with_Stochastic_Boundary_Conditions}

    We first prove the finite-volume entropy-dissipation identity.
    \begin{proof}[Proof of~\zcref{Prop:representation_of_finite_volume_derivative_of_entropy_for_modified_semigroup}]
        By reversibility,
        \begin{align}
            f_t
            :=\frac{\d (\mu T_t\Ssup{\L})}{\d \nu_\Lambda}
            = T_t\Ssup{\L} \frac{\d \mu_\Lambda}{\d \nu_\Lambda},
        \end{align} so that
        \begin{align}
            \RelEnt_{\Lambda}(\mu \Tl_t |\nu)
            = \nu_\Lambda[\phi(f_t)]
        \end{align} with \(\phi(x) := x \log(x) - x + 1 \geq 0\).

        From \zcref{lemma:L1_differentiability_from_finite_local_entropy} and the assumption \(\nu_\Lambda[f_0 \log f_0] = \RelEnt_{\Lambda}(\mu  |\nu) < \infty\), we get that \(s \mapsto f_s\) is continuously differentiable in \(L^1(\nu_\Lambda)\) with derivative \(s \mapsto \mathscr{L}_\Lambda f_s\). Formally, the chain rule gives
        \begin{align*}
            \frac{\d}{\d s} \phi(f_s)
            = \phi'(f_s) \mathscr{L}_\Lambda f_s
            = \log(f_s) \mathscr{L}_\Lambda f_s,
        \end{align*} 
        and integration together with GNZ suggests
        \begin{align*}
            \RelEnt_{\Lambda}(\mu |\nu) - \RelEnt_{\Lambda}(\mu \Tl_t |\nu)
            &= \nu_\Lambda[\phi(f_0)] - \nu_\Lambda[\phi(f_t)]
            = \int_0^{t} \d s \, \nu_\Lambda[\log(f_s) (-\mathscr{L}_\Lambda f_s)] \\
            &= \int_0^{t} \d s \, \mathcal{J}_\L(\mu \Tl_s\vert\nu). \nonumber 
        \end{align*}

        To justify this calculation, we use the fact from \zcref{lemma:L1_differentiability_from_finite_local_entropy} that for any continuously differentiable function \(\psi \colon [0, \infty) \to \R\) with \emph{bounded} derivative \(\psi'\), we have indeed
        \begin{align}\label{equation:integral_identity_approximation_finite_volume_relative_entropy_differences}
            &\nu_\Lambda[\psi(f_0)] - \nu_\Lambda[\psi(f_t)]
            = \int_0^{t} \d s \, \nu_\Lambda[\psi'(f_s) (-\mathscr{L}_\Lambda f_s)] \\
            &= \int_0^{t} \d s  \int_\Lambda \d x \int \nu_\Lambda(\d\eta_\Lambda) \, b_\Lambda^\nu(x,\eta_\Lambda) [f_s(\eta_\Lambda + \delta_x) - f_s(\eta_\Lambda)][\psi'(f_s(\eta_\Lambda + \delta_x)) - \psi'(f_s(\eta_\Lambda))]. \nonumber
        \end{align}
        It now just remains to approximate \(\phi\) by \(\phi_n\) such that
        \begin{enumerate}
            \item \(\phi_n \geq 0\) and \(\phi_n \uparrow \phi\) as \(n \to \infty\),
            \item \(\phi_n\) is differentiable with bounded derivative and satisfies
            \begin{align*}
                (b-a)(\phi_n'(b) - \phi_n'(a)) \uparrow  (b-a)(\log b - \log a)
            \end{align*} as \(n \to \infty\) for all \(b \geq a\),
        \end{enumerate} because then, by monotone convergence, the left-hand side of \eqref{equation:integral_identity_approximation_finite_volume_relative_entropy_differences} with \(\psi = \phi_n\) converges to \(\RelEnt_{\Lambda}(\mu |\nu) - \RelEnt_{\Lambda}(\mu \Tl_t |\nu) = \nu_\Lambda[\phi(f_0)] - \nu_\Lambda[\phi(f_t)]\) and the right-hand side of \eqref{equation:integral_identity_approximation_finite_volume_relative_entropy_differences} with \(\psi = \phi_n\) converges to \(\int_0^{t} \d s \, \mathcal{J}_\L(\mu \Tl_s\vert\nu)\).
        Such an approximation \((\phi_n)_{n \geq 2}\) was already used in proving a similar integral identity for the time evolution of relative entropy under the Ornstein--Uhlenbeck semigroup (i.e.\ birth-and-death dynamics with birth rate and death rate \(1\)) in \cite[Proof of Theorem 3.5]{Dello2024Wasserstein} and is e.g.\ given by defining
        \begin{align}
            \phi_n(x)
            := \int_1^x \d y \, \log [n \land (y \lor n^{-1})],
        \end{align}  
        cf.~\cite[(3.13)]{Dello2024Wasserstein}.
    \end{proof}

\subsection{The double-layer representation}\label{Sec:Fisher}
    We now identify the modified Fisher information with the double-layer entropy, starting with the required facts about Palm measures.
    \begin{lemma}
        Fix \(\mu \in \Pcal_{1,\theta}(\Omega)\) and let \(\Lambda \Subset \R^d\). Then, there exists a (jointly measurable) kernel
        \begin{align}\label{equation:jointly_measurable_Palm_kernel_space_time}
            (s, x) \mapsto (\mu \Tl_s )^x, 
            \quad s \geq 0, \, x \in \Lambda,
        \end{align} such that each \((\mu \Tl_s )^x\) is a finite measure on \(\Omega_\Lambda\) satisfying \((\mu \Tl_s )^x(\{\eta_\Lambda \in \Omega_\Lambda \,\colon\, x \not\in \eta_\Lambda\}) = 0\), and such that, for all \(s \geq 0\),
        \begin{align}\label{equation:defining_property_local_Palm_distribution}
            \int (\mu \Tl_s )(\d\eta_\Lambda)\int \eta_\Lambda(\d x)\, f(x, \eta_\Lambda)
            = \int_{\Lambda}\d x \int (\mu \Tl_s)^x(\d\eta_\Lambda) \, f(x, \eta_\Lambda)
        \end{align}
        for every measurable \(f \colon \Lambda \times \Omega_\Lambda \to [0, \infty]\).
    \end{lemma}
    \begin{remark}
        In particular, we get (jointly measurable) kernels
        \begin{align*}
            (s, x) \mapsto (\mu \Tl_s ) \circ G_x^{(\Lambda)},
            \quad  (s, x) \mapsto (\mu \Tl_s ) \ast G_x^{(\Lambda)}.
        \end{align*}
        The joint measurability of the corresponding kernels in infinite volume
        \begin{align*}
            (s, x) \mapsto (\mu T_s ) \circ G_x
            \quad  (s, x) \mapsto (\mu T_s ) \ast G_x
        \end{align*} is clear by their definition, i.e.\ only shifting \(x\) to the origin \(o\), and by the measurability in the initial conditions and c\`adl\`ag paths of the Markov process underlying \((T_s)_{s \geq 0}\).
    \end{remark}
    \begin{proof}
        First notice that by \eqref{equation:bounds_particle_numbers_local_dynamics_2}, every \(\mu \Tl_s\) has finite intensity and intensity measure absolutely continuous to Lebesgue measure on \(\Lambda\). Therefore, for each fixed \(s \geq 0\), corresponding local Palm distributions \(((\mu \Tl_s )^x)_{x \in \Lambda}\) of \(\mu \Tl_s \) exist such that \eqref{equation:defining_property_local_Palm_distribution} holds, cf.~\cite[Equations 13.1.4 and 13.1.5]{DaleyVereJones2008}. To construct versions that are jointly measurable in \((s, x)\), we only have to use universal disintegration instead of disintegration for each fixed \(s \geq 0\).
        Indeed, consider the family of Campbell measures
        \begin{align}
            C_{\mu \Tl_s }(\d x, \d \eta_\Lambda)
            := (\mu \Tl_s )(\d\eta_\Lambda) \, \eta_\Lambda(\d x),
        \end{align} which are finite measures on \(\Lambda \times \Omega_\Lambda\) with their first marginals absolutely continuous w.r.t.\ Lebesgue measure \(\d x\) on \(\Lambda\) by the discussion above. Moreover, the dependence on \(s\) is measurable, in the sense that \(s \mapsto C_{\mu \Tl_s}(A)\) is measurable for every measurable \(A \subseteq \Lambda \times \Omega_\Lambda\), by the measurability in the initial conditions of the underlying Markov process and its c\`adl\`ag paths.
        The universal disintegration theorem \cite[Corollary 1.26]{Kal17} therefore gives existence of a (jointly measurable) kernel as in \eqref{equation:jointly_measurable_Palm_kernel_space_time} such that for every \(s \geq 0\) we have
        \begin{align*}
            C_{\mu \Tl_s }(\d x, \d \eta_\Lambda)
            = \d x \, (\mu \Tl_s )^x(\d\eta_\Lambda),
        \end{align*} that is, \eqref{equation:defining_property_local_Palm_distribution} holds. Since \(C_{\mu \Tl_s }\) is finite and \(C_{\mu \Tl_s}(\{(x, \eta_\Lambda) \in \Lambda \times \Omega_\Lambda \,\colon\, x \notin \eta_\Lambda\}) = 0\), for each \(s \geq 0\), we can modify the kernel such that it stays measurable and with \eqref{equation:defining_property_local_Palm_distribution}, but we have for all \(s \geq 0\) and \(x \in \Lambda\) that \((\mu \Tl_s )^x\) is finite and \((\mu \Tl_s )^x(\{\eta_\Lambda \in \Omega_\Lambda \,\colon\, x \not\in \eta_\Lambda\}) = 0\).
    \end{proof}

    We make the following observation.
    \begin{lemma}\label{le_reversible_measure_cancellation_in_density}
        Let \(\mu_\Lambda \in \Pcal(\Omega_\Lambda)\) with \(\mu_\Lambda \ll \nu_\Lambda\) and \(\mu_\Lambda[N_\Lambda] < \infty\).
        The measures $\mu_\Lambda \ast G\Ssup{\Lambda}_x$ and $\mu_\Lambda \circ G\Ssup{\Lambda}_x$ are supported on 
        \begin{align*}
            A_x:= \{(\eta_\Lambda, \zeta_\Lambda)\colon \eta_\Lambda = \zeta_\Lambda \pm \delta_x\} \subset \Omega_\Lambda \times \Omega_\Lambda.
        \end{align*}

        Moreover, for such configurations, we have
        \begin{align}\label{equation:radon_nikodym_birth_death_flux_per_site}
            \frac{\d(\mu_\Lambda \ast G\Ssup{\Lambda}_x)}{\d(\nu_\Lambda \ast G\Ssup{\Lambda}_x)}(\eta_\Lambda, \zeta_\Lambda)
            = \frac{\d \mu_\Lambda}{\d \nu_\Lambda}(\eta_\Lambda),
            \quad \frac{\d(\mu_\Lambda \circ G\Ssup{\Lambda}_x)}{\d(\nu_\Lambda \circ G\Ssup{\Lambda}_x)}(\eta_\Lambda, \zeta_\Lambda)
            = \frac{\d \mu_\Lambda}{\d \nu_\Lambda}(\zeta_\Lambda)
        \end{align} 
        for \(\d x\)-a.e.\ \(x \in \Lambda\).
        Finally, for \(\d x\)-a.e.\ \(x \in \Lambda\), we have \(\nu_\Lambda \ast G\Ssup{\Lambda}_x = \nu_\Lambda \circ G\Ssup{\Lambda}_x\).
    \end{lemma}
    \begin{proof}
        The absolute continuity \(\mu_\Lambda \ll \nu_\Lambda\) implies absolute continuity of the intensity measure of \(\mu_\Lambda\) w.r.t.\ Lebesgue measure on \(\Lambda\), so that \((\mu_\Lambda^x)_{x \in \Lambda}\) are well-defined.
        The support assertion is immediate from the definition of $\mu_\Lambda \ast G\Ssup{\Lambda}_x$ and $\mu_\Lambda \circ G\Ssup{\Lambda}_x$.
        
        Now let \(f:= \d\mu_\L/\d\nu_\L\).
        For \zcref{equation:radon_nikodym_birth_death_flux_per_site}, we only have to show that for \(\d x\)-a.e.\ \(x \in \Lambda\) the following holds:
        \begin{align}\label{equation:radon_nikodym_local_palm}
            \mu_\Lambda^x
            = f \nu_\Lambda^x.
        \end{align}
        Indeed, let \(g\) be any bounded measurable function on \(\Omega_\Lambda\), and compute for arbitrary \(x \in \Lambda\) and \(\epsilon > 0\) that
        \begin{align*}
            &\int_{B_{\epsilon}(x) \cap \Lambda} \d y \int \mu_\Lambda^y(\d\eta_\Lambda) g(\eta_\Lambda)
            = \int \mu_\Lambda(\d\eta_\Lambda) \int \eta_\Lambda(\d y)\, g(\eta_\Lambda) {\1}_{B_\epsilon(x)}(y) \\
            &= \int \nu_\Lambda(\d\eta_\Lambda) \int \eta_\Lambda(\d y)\, g(\eta_\Lambda) f(\eta_\Lambda) {\1}_{B_\epsilon(x)}(y) 
            =\int_{B_{\epsilon}(x) \cap \Lambda} \d y \int \nu_\Lambda^y(\d\eta_\Lambda) f(\eta_\Lambda) g(\eta_\Lambda).
        \end{align*}
        Hence, 
        \begin{align*}
            \int \mu_\Lambda^x(\d\eta_\Lambda) g(\eta_\Lambda)
            = \int \nu_\Lambda^x(\d\eta_\Lambda) f(\eta_\Lambda) g(\eta_\Lambda)
        \end{align*}
        holds for \(\d x\)-a.e. \(x\).
        As we only have to check countably many \(g\)'s to decide whether two measures on \(\Omega_\Lambda\) (with the canonical $\sigma$-algebra) coincide (cf.~\cite[Proposition 2.8]{JansenGibbsPP}), the equality \eqref{equation:radon_nikodym_local_palm} is shown for \(\d x\)-a.e. \(x\).

        The symmetry  \(\nu_\Lambda \ast G\Ssup{\Lambda}_x = \nu_\Lambda \circ G\Ssup{\Lambda}_x\) follows from the GNZ equations for \(\nu\), compare with the proof of \zcref{lemma:finite_volume_infinitesimal_reversibility}.
    \end{proof}

    We can now put these ingredients together to prove the double-layer representation of the modified Fisher information.

\begin{proof}[Proof of~\zcref{proposition:double-layer-to-modified-fisher}]
    Assume first that \(\mu_\Lambda \ll \nu_\Lambda\). Put \(f = \d \mu_\Lambda / \ d\nu_\Lambda\). By \zcref{le_reversible_measure_cancellation_in_density}, we have
    \begin{align*}
            \frac{\d(\mu_\Lambda \ast G\Ssup{\Lambda}_x)}{\d(\nu_\Lambda \ast G\Ssup{\Lambda}_x)}(\eta_\Lambda, \zeta_\Lambda)
            = f(\eta_\Lambda),
            \quad \frac{\d(\mu_\Lambda \circ G\Ssup{\Lambda}_x)}{\d(\nu_\Lambda \circ G\Ssup{\Lambda}_x)}(\eta_\Lambda, \zeta_\Lambda)
            = \frac{\d(\mu_\Lambda \circ G\Ssup{\Lambda}_x)}{\d(\nu_\Lambda \ast G\Ssup{\Lambda}_x)}(\eta_\Lambda, \zeta_\Lambda)
            = f(\zeta_\Lambda)
    \end{align*} 
    for \(\d x\)-a.e.\ \(x\).
    Therefore, by symmetry of \(\nu_\Lambda \ast G\Ssup{\Lambda}_x = \nu_\Lambda \circ G\Ssup{\Lambda}_x \), and GNZ it follows
    \begin{align*}
         &\RelEnt\bigl( \mu_\Lambda \ast G\Ssup{\Lambda}_x \,\big\vert\, \mu_\Lambda \circ G\Ssup{\Lambda}_x \bigr)
         = \int (\nu_\Lambda \ast G\Ssup{\Lambda}_x)(\d\eta_\Lambda, \d\zeta_\Lambda) f(\eta_\Lambda) \log \frac{f(\eta_\Lambda)}{f(\zeta_\Lambda)}  \\
         &= \frac{1}{2} \int (\nu_\Lambda \ast G\Ssup{\Lambda}_x)(\d\eta_\Lambda, \d\zeta_\Lambda) \Big[f(\eta_\Lambda) \log \frac{f(\eta_\Lambda)}{f(\zeta_\Lambda)} + f(\zeta_\Lambda) \log \frac{f(\zeta_\Lambda)}{f(\eta_\Lambda)} \Big] \\
         &= \frac{1}{2}  \int (\nu_\Lambda \ast G\Ssup{\Lambda}_x)(\d\eta_\Lambda, \d\zeta_\Lambda) [f(\eta_\Lambda) -f(\zeta_\Lambda)][\log f(\eta_\Lambda) - \log f(\zeta_\Lambda)] \\
         &= \int \nu_\Lambda(\d\eta_\Lambda) \, b_\Lambda^\nu(x, \eta_\Lambda) [f(\eta_\Lambda + \delta_x) -f(\eta_\Lambda)][\log f(\eta_\Lambda + \delta_x) - \log f(\eta_\Lambda)] .
    \end{align*}
    Integrating over \(x \in \Lambda\) gives \eqref{equation:finite_volume_fisher_information_birth_death_flux_per_site_representation}.

    Let us finally discuss the case  \(\mu_\Lambda \not\ll \nu_\Lambda\), where the left-hand side of \eqref{equation:finite_volume_fisher_information_birth_death_flux_per_site_representation} is \(\infty\) by definition. We show that 
    \begin{align*}
        \mu_\Lambda \ast G\Ssup{\Lambda}_x 
        \not\ll \mu_\Lambda  \circ G\Ssup{\Lambda}_x 
    \end{align*} on a set of \(x\) of positive Lebesgue measure, so that the integral on the right-hand side takes the value \(\infty\) too.
    It suffices to exhibit an event \(H \in \mathcal{F}_\Lambda\) such that
    \begin{align}\label{equation:event_H_exhibiting_non_abs_cont}
        \int_\Lambda \d x \, \int \mu_\Lambda(\d\eta_\Lambda) b_\Lambda^\nu(x, \eta_\Lambda) {\1}_{H}(\eta_\Lambda + \delta_x) = 0,
        \quad \int_\Lambda \d x \, \int \mu_\Lambda^x(\d\eta_\Lambda) {\1}_{H}(\eta_\Lambda)
        > 0.
    \end{align}
    Indeed, then w.l.o.g.\ \(H \subseteq \{N_\Lambda = n\}\) for some \(n \in \N\) and 
    \begin{align*}
        \int_{\Lambda} \d x \, \mu_\Lambda \ast G\Ssup{\Lambda}_x(H \times \{N_\Lambda = n-1\})
        = 0 + \int \mu_\Lambda(\d \eta_\Lambda) \int \eta_\Lambda(\d x)\, {\1}_{H}(\eta_\Lambda)
        > 0,
    \end{align*} 
    as well as
    \begin{align}
        \int_{\Lambda} \d x \, \mu_\Lambda \circ G\Ssup{\Lambda}_x(H \times \{N_\Lambda = n-1\})
        = \int_{\Lambda} \d x \int \mu_\Lambda(\d\eta_\Lambda)\, b_\Lambda^\nu(x, \eta_\Lambda) {\1}_{H}(\eta_\Lambda + \delta_x)  + 0
        = 0,
    \end{align} 
    so that \(\mu_\Lambda \ast G\Ssup{\Lambda}_x \not\ll \mu_\Lambda  \circ G\Ssup{\Lambda}_x \) on a set of \(x\) of positive Lebesgue measure.

    To find \(H\) as in \zcref{equation:event_H_exhibiting_non_abs_cont}, let \(n \in \N\) be the smallest \(n\) such that \(\mu_\Lambda(\cdot \cap \{N_\Lambda = n\}) \not\ll \nu_\Lambda(\cdot \cap \{N_\Lambda = n\})\). A minimal \(n\) exists as \(\mu_\Lambda \not\ll  \nu_\Lambda\) by assumption and it is \(n \geq 1\) since \(\nu_\Lambda(\emptyset) > 0\). Now let \(H \subseteq \{N_\Lambda = n\}\) be any event with \(\mu_\Lambda(H) > 0\) but  \(\nu_\Lambda(H) = 0\).
    Then \zcref{equation:event_H_exhibiting_non_abs_cont} holds, as
    \begin{align*}
        \int_\Lambda \d x \, \int \mu_\Lambda^x(\d\eta_\Lambda) {\1}_{H}(\eta_\Lambda)
        = \int \mu_\Lambda(\d\eta_\Lambda) \int \eta_\Lambda(\d x)\, {\1}_{H}(\eta_\Lambda)
        = n \mu_\Lambda(H) > 0
    \end{align*} 
    but, because \(g := \frac{\d \mu_\Lambda(\cdot \cap \{N_\Lambda = n-1\})}{\d\nu_\Lambda(\cdot \cap \{N_\Lambda = n-1\})}\) exists by minimality of \(n\), and by the GNZ equations,
    \begin{align*}
        &\int_\Lambda \d x \, \int \mu_\Lambda(\d\eta_\Lambda) b_\Lambda^\nu(x, \eta_\Lambda) {\1}_{H}(\eta_\Lambda + \delta_x)
        = \int_\Lambda \d x \, \int \mu_\Lambda(\d\eta_\Lambda) b_\Lambda^\nu(x, \eta_\Lambda) {\1}_{H}(\eta_\Lambda + \delta_x) {\1}_{\{N_\Lambda = n-1\}}(\eta_\Lambda) \\
        &= \int_\Lambda \d x \, \int \nu_\Lambda(\d\eta_\Lambda) b_\Lambda^\nu(x, \eta_\Lambda) g(\eta_\Lambda) {\1}_{H}(\eta_\Lambda + \delta_x) 
        = \int \nu_\Lambda(\d\eta_\Lambda) \Big(\int \eta_\Lambda(\d x)\, g(\eta_\Lambda - \delta_x) \Big)  {\1}_{H}(\eta_\Lambda) 
        = 0.
    \end{align*}
    
\end{proof}

\subsection{Bounding the Fisher information from below}\label{Sec:Proof_of_Comparison_of_xi_Lambda_and_xi}
     
    \begin{proof}[Proof of~\zcref{proposition:lower_bound_xi_mu_fisher_information_finite_volume_dynamics}]
        First note that, by definition, we have
        \begin{align*}
            \sigma(\mu Ts)
            =  \RelEnt(\mu T_s \ast G_o | \mu T_s \circ G_o) 
            = \RelEnt(\mu T_s \ast G_x |\mu T_s \circ G_x) 
        \end{align*} 
        for all \(s \geq 0\) and \(x \in \mathbb{R}^d\).
        Now fix \(\epsilon > 0\). Suppose that \(\sigma(\mu T_s) < \infty\). Then, by the Donsker--Varadhan variational formula \eqref{equation:donsker_varadhan_variational_formula_finite_measures} and approximation, there exists a bounded and \(\Delta\)-local measurable function \(F\), for a box \(\Delta \supseteq B_{R}(o)\), such that
        \begin{align}\label{equation:upper_bound_test_function_relative_entropy_global_dynamics}
            &(\mu T_s \ast G_x)[1] + (\mu T_s \ast G_x)[F_x] 
            - (\mu T_s \circ G_x)[\exp(F_x)] \\
            &= (\mu T_s \ast G_o)[1] + (\mu T_s \ast G_o)[F] 
            - (\mu T_s \circ G_o)[\exp(F)] \nonumber
            \geq \RelEnt(\mu T_s \ast G_o | \mu T_s \circ G_o) - \epsilon \\
            &= \RelEnt(\mu T_s \ast G_x | \mu T_s \circ G_x) - \epsilon, \nonumber
        \end{align} 
        with \(F_x = F\circ \theta_x\) for all \(x \in \mathbb{R}^d\). We will go on to show that
        \begin{align}\label{equation:lower_bound_fisher_information_finite_volume_dynamics_goal}
            \liminf_{n \uparrow \infty}  n^{-d} \mathcal{J}_{\Lambda_n}(\mu T\Ssup{\L_n}_s\vert\nu)
            \geq (\mu T_s \ast G_o)[1] + (\mu T_s \ast G_o)[F] 
            - (\mu T_s \circ G_o)[\exp(F)].
        \end{align} If we had \(\sigma(\mu T_s) = \infty\), we could use the same argument below to find for every \(N \in \mathbb{N}\) a bounded and \(\Delta\)-local measurable function \(F\) such that
        \begin{align*}
            (\mu T_s \ast G_o)[1] + (\mu T_s \ast G_o)[F] 
            - (\mu T_s \circ G_o)[\exp(F)]  
            &\geq N
        \end{align*} and hence
        \begin{align*}
            \liminf_{n \uparrow \infty}  n^{-d} \mathcal{J}_{\Lambda_n}(\mu T\Ssup{\L_n}_s\vert\nu)\geq N,
        \end{align*} or ultimately
        \begin{align*}
            \liminf_{n \uparrow \infty}  n^{-d} \mathcal{J}_{\Lambda_n}(\mu T\Ssup{\L_n}_s\vert\nu) \geq \infty 
            = \sigma(\mu T_s).
        \end{align*}
        Now, let \(\Lambda_n^{\ominus}\) denote the set
        \begin{align*}
            \Lambda_n^{\ominus}
            := \Big\{ x \in \Lambda_n \colon \mathrm{dist}(\Delta+x, \Lambda_n^c) \geq \vert \Lambda_n \vert^{1/(2d)} \Big\},
        \end{align*} 
        which is such that 
        \(\vert \Lambda_n^{\ominus}\vert/\vert \Lambda_n\vert \xrightarrow[n \to \infty]{} 1\) but
        \begin{align}\label{equation:puffer_zone_is_valid}
            \log(\abs{\Lambda_n^{\ominus} \oplus \Delta})
            \leq \log(\abs{\Lambda_n})
            \leq \abs{\Lambda_n}^{1/(2d)}
            \leq \mathrm{dist}(\Lambda_n^{\ominus} \oplus \Delta, \Lambda_n^c).
        \end{align}

        Then, again by the Donsker--Varadhan variational formula \eqref{equation:donsker_varadhan_variational_formula_finite_measures},
        \begin{equation}\label{equation:lower_bound_fisher_information_finite_volume_dynamics}
        \begin{split}
            &\frac{1}{\vert\Lambda_n\vert} \mathcal{J}_{\Lambda_n}(\mu T\Ssup{\L_n}_s\vert\nu)
            = \frac{1}{\vert\Lambda_n\vert} \int_{\Lambda_n} \d x\,   \RelEnt\bigl( \mu T\Ssup{\L_n}_s \ast G\Ssup{\Lambda_n}_x \,\big\vert\, \mu T\Ssup{\L_n}_s  \circ G\Ssup{\Lambda_n}_x \bigr)  \\
            &\geq \Big(\frac{\vert\Lambda_n^{\ominus}\vert}{\vert\Lambda_n\vert}\Big) \cdot \frac{1}{\vert\Lambda_n^{\ominus}\vert} \int_{\Lambda_n^{\ominus}} \d x\,   \big\{ (\mu T\Ssup{\L_n}_s \ast G\Ssup{\Lambda_n}_x)[1] \\
            &\qquad\qquad\qquad\qquad + (\mu T\Ssup{\L_n}_s \ast G\Ssup{\Lambda_n}_x)[F_x] - (\mu T\Ssup{\L_n}_s \circ G\Ssup{\Lambda_n}_x)[\exp(F_x)]\big\}.
        \end{split}
        \end{equation}
        By definition, 
        \begin{align*}
            &(\mu T\Ssup{\L_n}_s \ast G\Ssup{\Lambda_n}_x)[F_x] \\
            &= \int \mu T\Ssup{\L_n}_s(\d \eta_{\Lambda_n}) \, b^\nu_{\Lambda_n}(x, \eta_{\Lambda_n}) \, F(\theta_x\eta_{\Lambda_n}, \theta_x\eta_{\Lambda_n} +\delta_o)
                + \int (\mu T\Ssup{\L_n}_s)^{x}(\d \eta_{\Lambda_n}) \, F( \theta_x\eta_{\Lambda_n},  \theta_x\eta_{\Lambda_n} -\delta_o),
        \end{align*} 
        and therefore
        \begin{equation}\label{equation:expansion_entropy_test_function_average_local_evolution}
        \begin{split}
            &\frac{1}{\vert\Lambda_n^{\ominus}\vert} \int_{\Lambda_n^{\ominus}} \d x\, (\mu T\Ssup{\L_n}_s \ast G\Ssup{\Lambda_n}_x)[F_x]  \\
            &= \frac{1}{\vert\Lambda_n^{\ominus}\vert} \int_{\Lambda_n^{\ominus}} \d x \int \mu T\Ssup{\L_n}_s(\d \eta_{\Lambda_n}) \, b^\nu_{\Lambda_n}(x, \eta_{\Lambda_n}) \, F(\theta_x\eta_{\Lambda_n}, \theta_x\eta_{\Lambda_n} +\delta_o) \\
            &\qquad + \frac{1}{\vert\Lambda_n^{\ominus}\vert} \int \mu T\Ssup{\L_n}_s(\d \eta_{\Lambda_n})\int \eta_{\Lambda_n^\ominus}(\d x) \, F( \theta_x\eta_{\Lambda_n},  \theta_x\eta_{\Lambda_n} -\delta_o).
        \end{split}
        \end{equation}
        We have to compare \eqref{equation:expansion_entropy_test_function_average_local_evolution} with
        \begin{align}\label{equation:expansion_entropy_test_function_average_global_evolution}
            &(\mu T_s \ast G_o)[F] 
            = \frac{1}{\vert\Lambda_n^{\ominus}\vert} \int_{\Lambda_n^{\ominus}} \d x\, (\mu T_s \ast G_x)[F_x] \\
            &= \frac{1}{\vert\Lambda_n^{\ominus}\vert} \int_{\Lambda_n^{\ominus}} \d x  \int \mu T_s(\d \eta) \, b(x, \eta) \, F(\theta_x\eta, \theta_x\eta +\delta_o)
            + \frac{1}{\vert\Lambda_n^{\ominus}\vert} \int \mu T_s(\d \eta) \int \eta_{\Lambda_n^\ominus}(\d x) \, F( \theta_x\eta,  \theta_x\eta -\delta_o). \nonumber
        \end{align}
        For the first summands on the right-hand side of \eqref{equation:expansion_entropy_test_function_average_local_evolution}, \eqref{equation:expansion_entropy_test_function_average_global_evolution},  we see by the finite-range property of \(b\), that
        \begin{equation}\label{equation:expansion_entropy_test_function_average_local_global_comparison}
        \begin{split}
            &\Big\vert \frac{1}{\vert\Lambda_n^{\ominus}\vert} \int_{\Lambda_n^{\ominus}} \d x \int \mu T\Ssup{\L_n}_s(\d \eta_{\Lambda_n}) \, b^\nu_{\Lambda_n}(x, \eta_{\Lambda_n}) \, F(\theta_x\eta_{\Lambda_n}, \theta_x\eta_{\Lambda_n} +\delta_o) \\
            &\qquad\qquad - \frac{1}{\vert\Lambda_n^{\ominus}\vert}  \int_{\Lambda_n^{\ominus}} \d x  \int \mu T_s(\d \eta) \, b(x, \eta) \, F(\theta_x\eta, \theta_x\eta +\delta_o) \Big\vert \\
            &= \frac{1}{\vert\Lambda_n^{\ominus}\vert}  \big\vert \mu T\Ssup{\L_n}_s[f_{\Lambda_n^{\ominus} \oplus \Delta}] - \mu T_s[f_{\Lambda_n^{\ominus} \oplus \Delta}] \big\vert 
        \end{split}
        \end{equation} 
        with the \((\Lambda_n^{\ominus} \oplus \Delta)\)-local function
        \begin{align*}
            f_{\Lambda_n^{\ominus} \oplus \Delta}(\eta)
            = \int_{\Lambda_n^{\ominus}} \d x \, b(x, \eta) \, F(\theta_x\eta, \theta_x\eta +\delta_o)
        \end{align*} 
        obeying \(\vert f_{\Lambda_n^{\ominus} \oplus \Delta} \vert \lesssim_{b, \Vert F\Vert_\infty} \abs{\Lambda_n^{\ominus}} \leq \abs{\Lambda_n^{\ominus} \oplus \Delta}\).
        It follows from \zcref{Lemma:finite_speed_of_propagation_version_1} and \eqref{equation:puffer_zone_is_valid} that the right-hand side of \eqref{equation:expansion_entropy_test_function_average_local_global_comparison} vanishes as \(n \to \infty\),
        \begin{align*}
            \frac{1}{\vert\Lambda_n^{\ominus}\vert}  \big\vert \mu T\Ssup{\L_n}_s[f_{\Lambda_n^{\ominus} \oplus \Delta}] - \mu T_s[f_{\Lambda_n^{\ominus} \oplus \Delta}] \big\vert  
            \xrightarrow[n \to \infty]{} 0.
        \end{align*}
        For the second summands on the right-hand side of \eqref{equation:expansion_entropy_test_function_average_local_evolution}, \eqref{equation:expansion_entropy_test_function_average_global_evolution}, we apply similarly \zcref{Lemma:finite_speed_of_propagation_version_1} with the \((\Lambda_n^{\ominus} \oplus \Delta)\)-local function
        \begin{align*}
            \widetilde{f}_{\Lambda_n^{\ominus} \oplus \Delta}(\eta)
            = \int \eta_{\Lambda_n^\ominus}(\d x) \, \, F(\theta_x\eta, \theta_x\eta - \delta_o)
        \end{align*} 
        obeying \(\vert \widetilde{f}_{\Lambda_n^{\ominus} \oplus \Delta}\vert \lesssim_{\Vert F\Vert_\infty} N_{\Lambda_n^{\ominus} \oplus \Delta}\), to also get
        \begin{align*}
            &\frac{1}{\vert\Lambda_n^{\ominus}\vert} \Big\vert \int \mu T\Ssup{\L_n}_s(\d \eta_{\Lambda_n}) \, \sum_{x \in \eta_{{\Lambda_n}^{\ominus}}} F( \theta_x\eta_{\Lambda_n},  \theta_x\eta_{\Lambda_n} -\delta_o) - \int \mu T_s(\d \eta) \int \eta_{\Lambda_n^\ominus}(\d x) \, F( \theta_x\eta,  \theta_x\eta -\delta_o) \Big\vert \\
            &\leq \frac{1}{\vert\Lambda_n^{\ominus}\vert}  \big\vert \mu T\Ssup{\L_n}_s[\widetilde{f}_{\Lambda_n^{\ominus} \oplus \Delta}] - \mu T_s[\widetilde{f}_{\Lambda_n^{\ominus} \oplus \Delta}] \big\vert  
            \xrightarrow[n \to \infty]{} 0. \nonumber
        \end{align*}
        This concludes the proof, as the other terms in \eqref{equation:lower_bound_fisher_information_finite_volume_dynamics} can be dealt with completely analogously, and we can send \(\epsilon \downarrow 0\) after taking the \(\liminf_{n \to \infty}\) to get first \eqref{equation:lower_bound_fisher_information_finite_volume_dynamics_goal} and then, with \eqref{equation:upper_bound_test_function_relative_entropy_global_dynamics}, finally \eqref{equation:lower_bound_xi_mu_fisher_information_finite_volume_dynamics}.
    \end{proof}

    \section{Time-stationary measures are Gibbs}\label{sec:time-stationary}

    \begin{proof}[Proof of~\zcref{corollary:time-stationary-iff-gibbs}]
    Since time-stationarity is weaker than reversibility, the inclusion \(\GG_\theta = \RR_\theta \subseteq \mathcal{S}_\theta\) follows already from \zcref{Prop:Reversible_Measures_equal_Gibbs_Measures}. 
    For the inclusion \(\mathcal{S}_\theta \subseteq \GG_\theta\) we notice that \zcref{Theorem:Main_Theorem_Limit_Points} directly yields that every time-stationary \(\mu \in \Pcal_\theta(\Omega)\) with \(\SpecEnt(\mu \,\vert\, \nu) < \infty\) has to be Gibbsian too, i.e., \(\mu \in \GG_\theta\). Hence, it suffices to show that a time-stationary \(\mu \in \Pcal_\theta(\Omega)\) automatically has finite specific relative entropy. The idea is that initial particles die with unit exponential rate and the newly born and still alive particles are stochastically dominated by a Poisson point process with intensity \((1-\e^{-t})\Vert b \Vert_\infty\), as can easily be seen via the graphical construction. Time-stationarity then forces \(\mu\) to be stochastically dominated by a Poisson point process of intensity \(\Vert b\Vert_\infty\).
    Indeed, in the context of the graphical construction, we decompose \(X_t^\eta\) into the surviving initial particles \(\eta_t\) and the newly born and still alive particles \(Y^\eta_t\). Then, we have the stochastic domination
    \begin{align*}
        Y_t^\eta
        &\leq \overline{Y}_{t} := \int_{\R^d \times (0, t] \times [0, \infty) \times [0, \Vert b\Vert_\infty]} \Ncal(\d x, \d s, \d r, \d u) \, \mathbbm{1}_{r+s > t} \, \delta_x 
    \end{align*} and we get for all \(\Lambda \Subset \R^d\) and all bounded increasing \(\mathcal{F}_\Lambda\)-measurable \(f \geq 0\) , that
    \begin{align*}
        \mu[f]
        = (\mu T_t)[f]
        = \mu[ \mathbb{E}[f(Y_t^\eta + \eta_t) ] ]
        \leq \mu[\mathbb{E}[f(\overline{Y}_t)] + \Vert f \Vert_\infty \mathbb{P}(N_\Lambda(\eta_t) > 0)], 
        \quad t \geq 0,
    \end{align*} by time-stationarity. Notice that each \(\overline{Y}_t\), \(t \geq 0\), is a Poisson point process of intensity \((1-\e^{-t})\Vert b\Vert_\infty\).
    In particular, we see that \(\mu\) is stochastically dominated by a Poisson process with intensity \(\Vert b\Vert_\infty\) by letting \(t \to \infty\) and using bounded convergence. Consequently, the \(n\)-th factorial moment measure of \(\mu\) is bounded by \(\Vert b\Vert_\infty^n \d x_1 \dots \d x_n\), and therefore \(j_{n, \Lambda}^\mu \leq \Vert b\Vert_\infty^n\) for the \(n\)-th Janossy densities \(j_{n, \Lambda}^\mu\), cf.~\cite[Definition 5.4.IV]{DaleyVereJones2003}, of \(\mu\), defined through
    \begin{align*}
        \mu_\Lambda[f]
        = \sum_{n = 0}^\infty \frac{1}{n!} \int_{\Lambda^n} \d x_1 \dots \d x_n \, j_{n, \Lambda}^\mu(x_1, \dots, x_n) \, f\big(\sum_{j = 1}^{n} \delta_{x_j}\big)
    \end{align*} for measurable \(f \geq 0\). For the Janossy densities \(j_{n, \Lambda}^\nu\) of \(\nu\) we get
    \begin{align*}
        j_{n, \Lambda}^\nu
        \geq \Vert 1/b\Vert_\infty^{-n} \e^{-\abs{\Lambda} \Vert b\Vert_\infty},
    \end{align*} using the DLR equations and bounding as in the proof of \zcref{Le:factor_property_of_the_Gibbs_measure}.
    Together, we get
    \begin{align*}
        \frac{\d \mu_\Lambda}{\d \nu_\Lambda} \mathbbm{1}_{N_\Lambda = n}
        \leq (\Vert b\Vert_\infty \Vert 1/b\Vert_\infty)^{n} \e^{\abs{\Lambda} \Vert b\Vert_\infty},
        \quad n \geq 0,
    \end{align*} and therefore
    \begin{align*}
        \RelEnt_\Lambda(\mu \,\vert\, \nu)
        &\leq \mu[\log (\Vert b\Vert_\infty \Vert 1/b\Vert_\infty)^{N_\Lambda} \e^{\abs{\Lambda} \Vert b\Vert_\infty}]
        = \abs{\Lambda}  \Vert b\Vert_\infty + \log(\Vert b\Vert_\infty \Vert 1/b\Vert_\infty) \mu[N_\Lambda] \\
        &\leq \abs{\Lambda} \Vert b\Vert_\infty \Big( 1+ \log(\Vert b\Vert_\infty \Vert 1/b\Vert_\infty) \Big),
    \end{align*} yielding
    \begin{align}
        \SpecEnt(\mu \,\vert\, \nu)
        \leq \Vert b\Vert_\infty \Big( 1+ \log (\Vert b\Vert_\infty \Vert 1/b\Vert_\infty) \Big) < \infty.
    \end{align}

    \end{proof}


\subsection*{Acknowledgments}
AZ thanks Lorenzo Dello Schiavo for the useful discussion.
BJ and JK gratefully received support by the Leibniz Association within the Leibniz Junior Research Group on \textit{Probabilistic Methods for Dynamic Communication Networks} as part of the Leibniz Competition (grant no.\ J105/2020).
BJ gratefully received support from Deutsche Forschungsgemeinschaft through DFG Project no.\ P27 within the SPP 2265. 
AZ is also affiliated with the University of Potsdam. The authors wish to thank the two anonymous referees for their careful reading and constructive comments, which greatly improved the presentation of this work.

\appendix

\section{Donsker--Varadhan variational formulas}

\begin{lemma}[Donsker-Varadhan variational formulas]\label{lemma:donsker_varadhan_variational_formulas}
    Suppose \(\alpha, \beta\) are finite measures defined on a measurable space \((\mathcal{X}, \mathcal{F})\) and
    \begin{align}
        I(\alpha \vert \beta)
        = \begin{cases}
            \int \d\alpha \log \frac{\d \alpha}{\d\beta}, & \alpha \ll \beta, \\
            +\infty, &\text{ otherwise}.
        \end{cases}
    \end{align}
    Let \(\mathfrak{B}_b(\mathcal{X}, \mathcal{F})\) denote the bounded measurable real-valued functions on \((\mathcal{X}, \mathcal{F})\).
    Then, the following variational formula holds
    \begin{align}\label{equation:donsker_varadhan_variational_formula_finite_measures}
        I(\alpha \vert \beta)
        = \alpha[1] + \sup_{f \in \mathfrak{B}_b(\mathcal{X}, \mathcal{F})} \big\{\alpha[f] - \beta[\e^f] \big\}
        = \sup_{f \in \mathfrak{B}_b(\mathcal{X}, \mathcal{F})} \big\{\alpha[f] - \beta[\e^{f-1}] \big\}
    \end{align} where the supremum ranges over bounded measurable functions.
    In the case that \(\alpha, \beta\) are probability measures, we additionally have
    \begin{align}\label{equation:donsker_varadhan_variational_formula_probability_measures}
        I(\alpha \vert \beta)
        = \sup_{f \in \mathfrak{B}_b(\mathcal{X}, \mathcal{F})} \big\{\alpha[f] - \log \beta[\e^{f}] \big\}.
    \end{align}
    In particular, if  \(\alpha, \beta\) are probability measures, we have the entropic inequality
    \begin{align}\label{equation:entropic_inequality}
        \alpha[f]
        \leq I(\alpha \vert \beta) + \log \beta[\e^{f}]
    \end{align} for all \(f \in \mathfrak{B}_b(\mathcal{X}, \mathcal{F})\), and then also all \(f\) such that \(\alpha[f]\) is well-defined (which is e.g.\ automatic when \(I(\alpha \vert \beta) < \infty\) and \(\e^f\) is \(\beta\)-integrable).
\end{lemma}
\begin{proof}
    Cf.~\cite[Theorem 2.1]{Donsker1976Asymptotic} or \cite[Corollary 4.15]{boucheron_concentration_2013} for the case of probability measures, \eqref{equation:donsker_varadhan_variational_formula_probability_measures}. The extension to finite measures follows by normalization.
\end{proof}

\section{Infinite-volume dynamics}\label{sec:appendix:infiniteconstruction}

In this section we prove the existence result of \cite{garcia2006spatial} as presented in \zcref{sec:construction_dynamics} above.

\subsection{Solution maps for the SDE}

We reformulate the existence result of \cite{garcia2006spatial} and its proof via \emph{solution maps}.

\begin{proposition}[Realization of the birth-and-death dynamics via solution maps]\label{proposition:birth_death_dynamics_via_solution_map}
    There exists a measurable map \(\mathcal{X} \colon \widehat{\Omega} \times \mathtt{N} \to E\) and a measurable map \(\widehat{h} \colon \Omega \times (0, \infty)^\mathbb{N} \to \widehat{\Omega}^{\ast}\) such that the following holds.
    \begin{enumerate}[(i)]
        \item On every filtered probability space as in \zcref{definition:SDE}, the stochastic process
        \begin{align*}
            X := \mathcal{X(\widehat{\eta}, \mathcal{N})}
        \end{align*} is a strong solution to \eqref{equation:SDE} with initial condition \(\mu \in \Pcal(\Omega)\). Moreover, any other solution of  \eqref{equation:SDE} on this filtered probability space, w.r.t.\ the same driving noise and with the same initial condition, is indistinguishable from \(X\). If it is the case that \(\mu(\Omega_{\mathrm{s}}) = 1\), then the strong solution \(X\) admits an indistinguishable \(E_{\mathrm{s}}\)-valued version.
        
        \item Any solution \(X = (X_t)_{t \geq 0}\) of \zcref{definition:SDE} is Markov w.r.t.\ its natural filtration.

        \item Fixing any complete probability space \((\Xi, \mathcal{G}, \mathbb{P})\) supporting independent unit-exponentials \((R_i)_{i \in \N}\) and an independent Poisson random measure \(\Ncal\) with the intensity from \zcref{definition:SDE}, we equip it with the filtration given by \(\mathcal{G}_t := \sigma((R_i)_{i \in \N}, \Ncal_{\R^d \times [0, t] \times [0, \infty)^2}) \lor \mathfrak{N}\), where \(\mathfrak{N}\) denotes all \(\mathbb{P}\)-null sets.
        Put 
        \begin{align}
            X^\eta := \mathcal{X}(\widehat{h}(\eta, (R_i)_{i \in \N}), \Ncal), \quad \eta \in \Omega.
        \end{align}
        Then, \((\eta, \xi) \mapsto X^\eta(\xi)\) is (jointly) measurable as map \(\Omega \times \Xi \to E\). 
        The random variable \(\widehat{\eta} = \widehat{h}(\eta, (R_i)_{i \in \N})\) is \(\eta\) independently marked with unit-exponentials and there exists an event of \(\mathbb{P}\)-probability one, on which, simultaneously for all \(\eta \in \Omega\), \(X^\eta\) is a strong solution to \eqref{equation:SDE} with the initial condition \(\widehat{\eta}\), on the here constructed filtered probability space.
    \end{enumerate}
\end{proposition}

\subsubsection{Deterministic solution maps}

To construct such solution maps, we need to introduce a few objects.
Let \(\mathtt{N}\) be the space of all simple locally finite counting measures on \(\R^d \times [0, \infty)^3\), equipped with the Borel $\sigma$-algebra of the vague topology, equivalently the sigma-algebra generated by all local counting variables, analogously to \(\Omega\). 
We define \(\widehat{\Omega}^\ast := \{\widehat{\eta} \in \widehat{\Omega} \,\colon\, \widehat{\eta}(\R^d \times \{0\}) = 0\}\) and a measurable subset \(\mathtt{N}^*\) of \(\mathtt{N}\) such that for given \(\mathcal{M} \in \mathtt{N}_\ast\) and \(\widehat{\eta} \in \widehat{\Omega}^\ast\) one can deterministically construct a solution to the equation
\begin{align}\label{equation:deterministic_version_of_SDE}
    Y_t(B)
    &= \int_{B \times (0, t] \times [0, \infty)^2} \Mcal(\d x, \d s, \d r, \d u) {\1}_{(t-s, \infty)}(r) {\1}_{[0, q(x, Y_{s-})]}(u)
    + \int_{B \times [0, \infty)} \widehat{\eta}(\d x, \d r) {\1}_{(t, \infty)}(r),
\end{align}
Indeed, let \(\mathtt{N}^*\) consist of all \(\mathcal{M} \in \mathtt{N}\) such that 
\begin{enumerate}[(i)]
    \item For every \(\Lambda \Subset \R^d\) and \(T > 0\),
    \begin{align*}
        \mathcal{M}(\Lambda \times [0, T] \times [0, \infty) \times [0, \Vert q\Vert_\infty])
        < \infty.
    \end{align*}

    \item For every \(\Lambda \Subset \R^d\) and every \(T > 0\) there exists an \(M \in \N\) such that every chain of successive birth proposals in \(\mathcal{M}\), with acceptance marks \(u \in [0, \Vert q\Vert_\infty]\), ending in \(\Lambda\) before or at time \(T\), has length at most \(M\).
    More precisely, for every \(m > M\) and \(0 < t_1 < \dots < t_m \leq T\), \(x_1, \dots, x_{m-1} \in \R^d\), \(x_m \in \Lambda\) and \(\vert x_{i+1} - x_i\vert \leq R\) for \(1 \leq i < m\), we have 
    \begin{align*}
        \mathcal{M}(\{x_i\} \times \{t_i\} \times [0, \infty) \times [0, \Vert q\Vert_\infty])
        = 0
    \end{align*} 
    for at least one \(1 \leq i \leq m\).
\end{enumerate}

We show in \zcref{lemma:almost_surely_good_birth_proposals} that \(\mathtt{N}^\ast\) is measurable and that the Poisson random measures considered for the SDE \eqref{equation:SDE} give probability one to this event.
It is also clearly a translation-invariant event.

\begin{proposition}[Deterministic solution maps]\label{proposition:deterministic_solution_maps}
    There exists a measurable map \(\mathcal{X} \colon \widehat{\Omega} \times \mathtt{N} \to E\) with the following properties.
    \begin{enumerate}[(i)]
        \item \emph{Deterministic existence and uniqueness.} For every \(\widehat{\eta} \in \widehat{\Omega}^\ast\) and \(\mathcal{M} \in \mathtt{N}^\ast\), \(Y = \mathcal{X}(\widehat{\eta}, \mathcal{M})\) is the unique \(Y \in E\) solving the equation \eqref{equation:deterministic_version_of_SDE} simultaneously for all \(B \Subset \R^d\) and \(t \geq 0\).
        For \((\widehat{\eta}, \mathcal{M}) \in (\widehat{\Omega}^\ast \times \mathtt{N}^\ast)^c\) we choose \(\mathcal{X}(\widehat{\eta}, \mathcal{M}) = (t \mapsto \widehat{\eta}(\cdot \times [0, \infty)))\). 
        Then, for every \((\widehat{\eta}, \mathcal{M}) \in \widehat{\Omega} \times \mathtt{N}\), we have that \((\mathcal{X}(\widehat{\eta}, \mathcal{M}))_0 = \widehat{\eta}(\cdot \times [0, \infty))\), and that \((\mathcal{X}(\widehat{\eta}, \mathcal{M}))_{\Lambda}\) has only finitely many jumps in \([0, T]\) for every \(\Lambda \Subset \R^d\) and \(T > 0\).
        
        \item \emph{Deterministic non-anticipativity}. 
        Write
        \begin{align*}
            \mathcal{M}_{[0, T]} = \mathcal{M}\vert_{\R^d \times [0, T] \times [0, \infty)^2} = \mathcal{M}(\cdot \cap (\R^d \times [0, T] \times [0, \infty)^2))
        \end{align*} for the version of \(\mathcal{M} \in \mathtt{N}\) that keeps only the birth proposals from the time interval \([0, T]\), \(T > 0\).
        We have, for \(\widehat{\eta} \in \widehat{\Omega}\), \(\mathcal{M} \in \mathtt{N}^\ast\), and for every \(T > 0\),
        \begin{align*}
            \mathcal{X}(\widehat{\eta}, \mathcal{M})\vert_{[0, T]}
            = \mathcal{X}(\widehat{\eta}, \mathcal{M}_{[0, T]})\vert_{[0, T]}
        \end{align*}

        \item \emph{Deterministic translation covariance}. If the birth rate \(q\) is translation invariant as in \eqref{equation:translation_covariance_birth_rate_q}, then the solution map is translation covariant:
        \begin{align}
            \mathcal{X}(\theta_z \widehat{\eta}, \theta_z \mathcal{M})
            = \theta_z \mathcal{X}(\widehat{\eta}, \mathcal{M})
        \end{align} for all \(z \in \R^d\), \(\widehat{\eta} \in \widehat{\Omega}\) and \(\mathcal{M} \in \mathtt{N}\).
    \end{enumerate}
\end{proposition}

\begin{proof}[Proof of \zcref{proposition:deterministic_solution_maps}]
    We follow the idea of Picard-iteration from the proof of \cite[Theorem 2.13]{garcia2006spatial}, and make its properties explicit for the (measurability) statements of our proposition. The finite-range property of \(q\) and the deterministic definition of \(\mathtt{N}^\ast\) replace their use of measuring influence via \(a(x, y) := \sup_{\eta \in \Omega} \abs{q(x, \eta + \delta_y) - q(x, \eta)}\) and probabilistic estimates. Similarly, a Gronwall estimate is not needed as in our setting the Picard-iterates will not only converge but locally stabilize eventually.
    
    \textit{(i)}: 
    Fix \(\widehat{\eta} \in \widehat{\Omega}^\ast\) and \(\mathcal{M} \in \mathtt{N}^\ast\).
    Let
    \begin{align}
        Y^{(0)}_t
        := \int_{\R^d \times [0, \infty)} \widehat{\eta}(\d x, \d r) {\1}_{r > t} \delta_{x}, \quad t \geq 0,
    \end{align} and recursively
    \begin{align}
        Y^{(n+1)}_t
        := \int_{\R^d \times (0, t] \times [0, \infty)^2} \mathcal{M}(\d x, \d s, \d r, \d u) {\1}_{r > t-s} {\1}_{[0, q(x, Y^{(n)}_{s-})]}(u) \delta_{x}
        + \int_{\R^d \times [0, \infty)} \widehat{\eta}(\d x, \d r) {\1}_{r > t} \delta_{x}
    \end{align} for all \(t \geq 0\).
    It is then clear by the local finiteness of the initial condition \(\widehat{\eta}\) and \textit{(i)} that each \(Y^{(n)}\) is locally finite and has  c\`adl\`ag paths with locally finitely many jumps. For any \(\Lambda \Subset \R^d\) and \(T > 0\), the number of possible jumps in \(\Lambda\) during \([0, T]\) does not depend on \(n\), but only on \(\mathcal{M}(\Lambda \times (0, T] \times [0, \infty) \times [0, \Vert q\Vert_\infty])\) and \(\widehat{\eta}\) on \(\Lambda \times [0, T]\).

    We also have local stabilization in the following sense. For \(n \in \N\), a difference between successive Picard-iterates \(Y^{(n+1)}\), \(Y^{(n)}\) at time \(t\) can only be caused by a birth proposal coming from \(\mathcal{M}\) at time \(s \leq t\) with different acceptance decisions, since the death jumps from the initial particles agree. By the finite-range property of \(q\) these different decisions require a discrepancy between \(Y^{(n)}\) and \(Y^{(n-1)}\) within distance \(R\) of that birth proposal immediately before \(s\), and hence at a strictly earlier time. Recursion of this argument yields the following fact.
    For any \(n \in \N\) and \(\Lambda \Subset \R^d\), \(T > 0\), if
    \begin{align*}
        (Y^{(n)}_t)_\Lambda \neq (Y^{(n-1)}_t)_\Lambda, \quad \text{for some } t \in [0, T],
    \end{align*} then \(\mathcal{M}\) has to contain a chain of at least \(n\) successive birth proposals with acceptance marks \(u \in [0, \Vert q\Vert_\infty]\), ending in \(\Lambda\) before or at time \(T\). 
    Hence, by \textit{(ii)} in the definition of \(\mathtt{N}^\ast\), for every \(\Lambda \Subset \R^d\), \(T > 0\), there is an \(M \in \N\) such that
    \begin{align*}
        (Y^{(n)}_t)_\Lambda = (Y^{(M)}_t)_\Lambda, \quad \text{for all } t \in [0, T],
    \end{align*} holds for all \(n \geq M\).
    These stabilized restrictions are consistent too as \(\Lambda\) increases to \(\R^d\) and \(T \to \infty\).
    
    Call the limiting path \(Y\), which is then by stabilization an element of \(E\) with locally finitely many jumps, and with \(Y_0 = \widehat{\eta}(\cdot \times [0, \infty))\).
    It also solves the deterministic equation \eqref{equation:deterministic_version_of_SDE} by stabilization and the finite-range property of \(q\).

    Regarding uniqueness, if \(Z \in E\) is any other solution of \eqref{equation:deterministic_version_of_SDE}, it automatically has locally finitely many jumps by being a solution and property \textit{(i)} in the definition of \(\mathtt{N}^\ast\). One can inductively see as above, that if there is some discrepancy between \(Z_t\) and \(Y_t^{(n)}\) in \(\Lambda \Subset \R^d\) for \(t \in [0, T]\), then \(\mathcal{M}\) has to contain a chain of at least \(n\) successive birth proposals with acceptance marks \(u \in [0, \Vert q\Vert_\infty]\), ending in \(\Lambda\) before or at time \(T\). 
    Hence, by \textit{(ii)} in the definition of \(\mathtt{N}^\ast\), \(Z\) and \(Y_t^{(n)}\) agree on \(\Lambda \times [0, T]\) for sufficiently large \(n\), so that \(Z\) and \(Y\) agree there as well.

    Now define \(\mathcal{X}(\widehat{\eta}, \mathcal{M}) := Y = Y(\widehat{\eta}, \mathcal{M})\) as constructed above for \(\widehat{\eta} \in \widehat{\Omega}^\ast\) and \(\mathcal{M} \in \mathtt{N}^\ast\). On the complement of \(\widehat{\Omega}^\ast \times \mathtt{N}^\ast\) we just take the path which is constantly equal to \( \widehat{\eta}(\cdot \times [0, \infty))\), as already stated in the proposition.
    Measurability follows because \(\widehat{\Omega}^\ast\), \(\mathtt{N}^\ast\) are measurable events and every Picard-iterate \(Y^{(n)}\) is a measurable map of \((\widehat{\eta}, \mathcal{M})\), so that their limit is measurable too.

    \textit{(ii)}: If \(\mathcal{M}\) belongs to \(\mathtt{N}^\ast\), then \(\mathcal{M}_{[0, T]}\) belongs to \(\mathtt{N}^\ast\) too, as we are at most dropping birth proposals. The statement is then clear for \(\widehat{\eta} \in \widehat{\Omega}^\ast\) from the construction of  \(\mathcal{X}\) in \textit{(i)}, as the corresponding Picard-iterates have to agree on the time interval \([0, T]\). For \(\widehat{\eta} \not\in \widehat{\Omega}^\ast\), both \(\mathcal{X}(\widehat{\eta}, \mathcal{M})\)  and \(\mathcal{X}(\widehat{\eta}, \mathcal{M}_{[0, T]})\) are the same constant path by definition.

    \textit{(iii)}: Suppose \(q\) is translation invariant. The sets \(\widehat{\Omega}^\ast\) and \(\mathtt{N}^\ast\) are translation invariant. For \(\widehat{\eta} \in \widehat{\Omega}^\ast\) and \(\mathcal{M} \in \mathtt{N}^\ast\), and for all of the Picard-iterates \(Y^{(n)}(\widehat{\eta}, \mathcal{M})\) in the construction of \(\mathcal{X}\) in \textit{(i)}, we have that
    \begin{align*}
        Y^{(n)}(\theta_z \widehat{\eta}, \theta_z \mathcal{M})
        = \theta_z  Y^{(n)}(\widehat{\eta}, \mathcal{M}), \quad z \in \R^d,
    \end{align*} 
    by the translation invariance of \(q\), so that this property also passes to the limit.
    On the complement of \(\widehat{\Omega}^\ast \times \mathtt{N}^\ast\), translation covariance follows from the constant path definition.
\end{proof}

As in \cite{garcia2006spatial}, we will also have to consider versions \(\widehat{\mathcal{X}}\) of the solution map \(\mathcal{X}\) which retain as marks the residual lifespans of particles.
We define \(\widehat{\mathcal{X}} \colon \widehat{\Omega} \times \mathtt{N} \to D([0, \infty); \widehat{\Omega})\) by
\begin{align}\label{equation:definition_marked_solution}
    \widehat{\mathcal{X}}_t(\widehat{\eta}, \mathcal{M})
    &:= \int_{\R^d \times (0, t] \times [0, \infty)^2} \mathcal{M}(\d x, \d s, \d r, \d u) {\1}_{(t-s, \infty)}(r) \, {\1}_{[0, q(x, \mathcal{X}(\widehat{\eta}, \mathcal{M})_{s-})]}(u) \, \delta_{(x, (s+r) - t)}   \\
    &\qquad + \int_{\R^d \times (t, \infty)} \widehat{\eta}(\d x, \d r) \delta_{(x, r-t)} \nonumber
\end{align} for \(\widehat{\eta} \in \widehat{\Omega}^\ast\), \(\mathcal{M} \in \mathtt{N}^\ast\).
For \((\widehat{\eta}, \mathcal{M}) \in (\widehat{\Omega}^\ast \times \mathtt{N}^\ast)^c\), we set \(\widehat{\mathcal{X}}(\widehat{\eta}, \mathcal{M}) := (t \mapsto \widehat{\eta})\).

The following lemma will be important to see the Markov property of the solution map with random inputs. It follows directly from uniqueness of the solution to the deterministic equation.
\begin{lemma}[Deterministic restart]\label{lemma:deterministic_restart}
    Let \(\widehat{\eta} \in \widehat{\Omega}^\ast\), \(\mathcal{M} \in \mathtt{N}^\ast\) and \(s, t \geq 0\).
    Write \(\mathcal{M}^{(s)}\) for \(\mathcal{M}\) with birth proposal times shifted by \(s\),
    \begin{align*}
        \mathcal{M}^{(s)}(A)
        := \int \mathcal{M}(\d x, \d v, \d r, \d u) \, {\1}_{v > s} {\1}_A(x, v - s, r, u)
    \end{align*}
    Then,
    \begin{align}\label{equation:deterministic_restart_marked}
        \widehat{\mathcal{X}}(\widehat{\eta}, \mathcal{M})_{s+t}
        = \Big(\widehat{\mathcal{X}}\big((\widehat{\mathcal{X}}(\widehat{\eta}, \mathcal{M}))_s, \mathcal{M}^{(s)} \big)\big)_{t}
    \end{align} and
    \begin{align}\label{equation:deterministic_restart_unmarked}
        \mathcal{X}(\widehat{\eta}, \mathcal{M})_{s+t}
        = \Big(\mathcal{X}\big((\widehat{\mathcal{X}}(\widehat{\eta}, \mathcal{M}))_s, \mathcal{M}^{(s)} \big)\big)_{t}
    \end{align}
\end{lemma}
\begin{proof}
    Write \(X := \mathcal{X}(\widehat{\eta}, \Mcal)\) and \(\widehat{X} := \widehat{\mathcal{X}}(\widehat{\eta}, \Mcal)\). Notice that \(\Mcal^{(s)} \in \mathtt{N}^\ast\) for every \(s \geq 0\), since the conditions \textit{(i)}, \textit{(ii)} in the definition of \(\mathtt{N}^\ast\) on a time interval \([0, T]\) are satisfied for \(\Mcal^{(s)}\) since \(\Mcal\) satisfies the corresponding conditions on \([0, s + T]\). Moreover, \(\widehat{X}_s \in \widehat{\Omega}^\ast\).
    It follows by construction of \(\mathcal{X}\) in \zcref{proposition:deterministic_solution_maps}, that simultaneously for all \(B \Subset \R^d\) and \(t \geq 0\), 
    \begin{align}\label{equation:proof_unmarked_restart_identity}
        &\mathcal{X}(\widehat{\eta}, \Mcal)_{s+t} (B) 
        = X_{s+t} (B) \\
        &= \int_{B \times (0, s+t] \times [0, \infty)^2} \Mcal(\d x, \d v, \d r, \d u) {\1}_{v+r > s+t} {\1}_{u \leq q(x, X_{v-})}
        + \int_{B \times [0, \infty)} \widehat{\eta}(\d x, \d r) {\1}_{r > s+t}\nonumber \\
        &= \int_{B \times (s, s+t] \times [0, \infty)^2} \Mcal(\d x, \d v, \d r, \d u) {\1}_{v+r > s+t} {\1}_{u \leq q(x, X_{v-})} \\
        &\qquad + \Big(\int_{B \times (0, s] \times [0, \infty)^2} \Mcal(\d x, \d v, \d r, \d u) {\1}_{v+r > s+t} {\1}_{u \leq q(x, X_{v-})} \nonumber \\
        &\qquad\qquad + \int_{B \times [0, \infty)} \widehat{\eta}(\d x, \d r) {\1}_{r > s+t} \Big) \nonumber \\
        &= \int_{B \times (0, t] \times [0, \infty)^2} \Mcal^{(s)}(\d x, \d v, \d r, \d u) {\1}_{v+r > t} {\1}_{u \leq q(x, X_{(s+v)-})} 
        + \int_{B \times [0, \infty)} \widehat{X}_s(\d x, \d r) {\1}_{r > t}. \nonumber
    \end{align}
    Hence, \(t \mapsto X_{s+t}\) solves the deterministic equation \eqref{equation:deterministic_version_of_SDE} with the input \((\widehat{X}_s, \Mcal^{(s)})\) and therefore, by uniqueness in \textit{(i)} of \zcref{proposition:deterministic_solution_maps}, has to equal
    \begin{align*}
        t \mapsto \mathcal{X}(\widehat{X}_s, \Mcal^{(s)})_t
        = \Big(\mathcal{X}\big((\widehat{\mathcal{X}}(\widehat{\eta}, \mathcal{M}))_s, \mathcal{M}^{(s)} \big)\big)_{t}.
    \end{align*}
    This proves \eqref{equation:deterministic_restart_unmarked}. Using in the defining equation \eqref{equation:definition_marked_solution} of \(\widehat{X}_{s+t}\) the same decomposition into birth proposals at times before and including \(s\) and into birth proposals strictly after time \(s\) as in \eqref{equation:proof_unmarked_restart_identity}, and then the already proved identity \eqref{equation:deterministic_restart_unmarked} to identify the accepted birth proposals yields \eqref{equation:deterministic_restart_marked}.
\end{proof}

\subsubsection{Solution maps with random input}

\begin{lemma}\label{lemma:almost_surely_good_birth_proposals}
    The set \(\mathtt{N}^\ast\) is measurable and, for every Poisson random measure \(\Ncal\) on \(\R^d \times [0, \infty)^3\) with intensity measure \(\d x \otimes \d s \otimes \e^{-r} \d r \otimes \d u\), we have
    \begin{align}
        \mathbb{P}(\mathcal{N} \in \mathtt{N}^\ast)
        = 1.
    \end{align}
\end{lemma}
\begin{proof}
    Let, for \(\Lambda \Subset \R^d\), \(T > 0\) and \(M \in \N\),
    \begin{align*}
        \mathtt{C}^\ast_{\Lambda, T, M}
        := \{\Mcal \in \mathtt{N} \,\colon\, &\exists 0 < t_1 < \cdots < t_M \leq T, x_1, \dots, x_M \in \R^d \text{ such that} \\
            &\quad \Mcal(\{x_i\} \times \{t_i\} \times \mathbb{R}_+ \times [0, \Vert q\Vert_\infty]) \geq 1 \text{ for } 1\leq i \leq M, \\
            &\quad x_M \in \Lambda, \abs{x_i - x_{i+1}} \leq R \text{ for } 1 \leq i < M
        \}
    \end{align*}
    and
    \begin{align*}
        \mathtt{N}^\ast_{\Lambda, T, M}
        := \{\Mcal \in \mathtt{N} \,\colon\, \mathcal{M}(\Lambda \times [0, T] \times [0, \infty) \times [0, \Vert q\Vert_\infty]) < \infty \} \cap (\mathtt{C}^\ast_{\Lambda, T, M})^c
    \end{align*}
    Then, by definition of \(\mathtt{N}^\ast\) and monotonicity in \(\Lambda, T\),
    \begin{align*}
        \mathtt{N}^\ast
        = \bigcap_{n \in \N} \bigcup_{M \in N} \mathtt{N}^\ast_{\Lambda_n, n, M}.
    \end{align*}
    The set \(\{\Mcal \in \mathtt{N} \,\colon\, \mathcal{M}(\Lambda \times [0, T] \times [0, \infty) \times [0, \Vert q\Vert_\infty]) < \infty \}\) is clearly measurable and letting
    \begin{align*}
        &\mathfrak{c}_{\Lambda, T, M}((x_1, s_1, r_1, u_1), \dots, (x_M, s_M, r_M, u_M)) \\
        &:= {\1}_{0 < s_1 < s_2 < \dots < s_M \leq T} {\1}_{\abs{x_{i+1} - x_{i}} \leq R \,\, \forall 1 \leq i < M} {\1}_{x_M \in \Lambda} {\1}_{u_i \in [0, \Vert q \Vert_\infty] \,\, \forall 1 \leq i \leq M},
    \end{align*} we see that
    \begin{align*}
        \mathtt{C}^\ast_{\Lambda, T, M}
        = \Big \{\Mcal \in \mathtt{N} \,\colon\, \int \d\mathcal{M}^{\otimes M} \, \mathfrak{c}_{\Lambda, T, M}  > 0 \Big\}
    \end{align*} is measurable too.

    Now, we have already
    \begin{align*}
        \mathbb{E}[\Ncal(\Lambda \times [0, T] \times [0, \infty) \times [0, \Vert q\Vert_\infty])] 
        = \abs{\Lambda} T \Vert q\Vert_\infty
        < \infty,
    \end{align*} so that \(\mathbb{P}(\Ncal(\Lambda \times [0, T] \times [0, \infty) \times [0, \Vert q\Vert_\infty]) < \infty) = 1\).
    It remains to show \(\mathbb{P}(\mathcal{N} \in  \mathtt{C}^\ast_{\Lambda, T, M}) \xrightarrow[M \to \infty]{} 0\) for fixed \(\Lambda \Subset \R^d\), \(T > 0\).
    Let
    \begin{align*}
        &\mathfrak{C}_{\Lambda, T, M}  \\
        &:= \int \Ncal^{\otimes M}(\d(x_1, s_1, r_1, u_1), \dots, \d(x_M, s_M, r_M, u_M)) \, \mathfrak{c}((x_1, s_1, r_1, u_1), \dots, (x_M, s_M, r_M, u_M)) 
    \end{align*} with \(\mathfrak{c} = \mathfrak{c}_{\Lambda, T, M}\).
     We first bound
    \begin{align*}
        {\1}_{\mathtt{C}^\ast_{\Lambda, T, M}}
        \leq \mathfrak{C}_{\Lambda, T, M} 
    \end{align*} and now the expectation of the random variable on the right-hand side.
    Indeed, by the multivariate Mecke equation (which is applicable because \(\mathfrak{c}\) is zero on tuples with non-pairwise-different entries) and iterative integration in every variable,
    \begin{align*}
        &\mathbb{P}(\mathcal{N} \in  \mathtt{C}^\ast_{\Lambda, T, M}) 
        \leq \mathbb{E}[\mathfrak{C}_{\Lambda, T, M} ] 
        = \int (\d x \otimes \d s \otimes \e^{-r} {\1}_{[0, \infty)}(r) \d r \otimes \d u)^{\otimes M} \, \mathfrak{c}_{\Lambda, T, M} \\
        &= \Vert q \Vert_\infty^M \Big(\int_{\Lambda \times (\R^d)^{M-1}} \d x_M \dots \d x_1 {\1}_{\abs{x_{M} - x_{M-1}} \leq R}  \dots {\1}_{\abs{x_{2} - x_{1}} \leq R} \Big) \Big(\int_{0}^{T} \d s_1 \int_{s_1}^{T} \d s_2 \dots \int_{s_{M-1}}^T \d s_M \Big)  \\
        &=  \Vert q \Vert_\infty^M \abs{\Lambda} \abs{B_R(0)}^{M-1} \frac{T^M}{M!}
        \xrightarrow[M \to \infty]{} 0.
    \end{align*}

\end{proof}

Using measurable enumerations of the (occurrences of) points in the elements of \(\Omega\), cf.\cite[Proposition 6.3]{Last_Penrose_2017}, one can construct a measurable marking with i.i.d.\ unit-exponentials as follows, also compare with the definition and corresponding Laplace formula in \cite[Section 5.2, Proposition 5.4]{Last_Penrose_2017}.
\begin{lemma}\label{lemma:exponential_lifespan_marking}
    Let \(\mathfrak{r} := \mathrm{Exp}(1)^{\otimes \N}\).
    There exists a measurable map \(\widehat{h} \colon \Omega \times (0, \infty)^\mathbb{N} \to \widehat{\Omega}^{\ast}\) and a probability kernel \(\widehat{\alpha}\) from \(\Omega \to \widehat{\Omega}\) such that the following holds.
    \begin{enumerate}[(i)]
        \item For every \(\eta \in \Omega\) and \(r \in (0,\infty)^\N\),
        \begin{align*}
            \widehat{h}(\eta, r)(\cdot \times [0, \infty))
            = \eta.
        \end{align*}

        \item For every \(\eta \in \Omega\), 
        \begin{align*}
            \widehat{\alpha}(\eta, \cdot)
            := \mathfrak{r} \circ \widehat{h}(\eta, \cdot)^{-1}
        \end{align*} 
        is the law arising from assigning independent unit-exponential marks to each (occurrence of a) point in \(\eta\). 
        More precisely, for every \(g \in C_c(\R^d \times [0, \infty))\), \(g \geq 0\), we have
        \begin{align*}
            \int_{\widehat{\Omega}} \widehat{\alpha}(\eta, \d\widehat{\eta}) \e^{-\int \widehat{\eta}(\d x, \d r) g(x, r)}
            = \prod_{x \in \eta} \int_0^\infty \d r \, \e^{-g(x, r)} \e^{-r},
        \end{align*} 
        with the product counting multiplicities.
    \end{enumerate}
\end{lemma}

Finally, we needed the following lemma as input in proving the Markov property of \((X_t)_{t \geq 0}\), which is really just the memorylessness of the exponential distribution.
\begin{lemma}\label{lemma:memorylessness}
    Consider the setting of \zcref{definition:SDE}.
    Let \(X := \mathcal{X}(\widehat{\eta}, \Ncal)\) and \(\widehat{X} := \widehat{\mathcal{X}}(\widehat{\eta}, \Ncal)\).
    Then,
    \begin{align}\label{equation:conditional_law}
        \mathrm{Law}(\widehat{X}_t \,\vert\, \sigma(X_s \colon 0 \leq s \leq t))
        = \widehat{\alpha}(X_t, \cdot).
    \end{align}
\end{lemma}
\begin{proof}
    As \(X\), \(\widehat{X}\) are measurable functions of \((\widehat{\eta}, \Ncal)\), the validity of \eqref{equation:conditional_law} only depends on the joint distribution of \((\widehat{\eta}, \Ncal)\). Consequently, we can work with the following representation using measurable enumerations and independent markings.
    Fix \(t \geq 0\).
    Let \(\Pi_t\) be a Poisson random measure on \(\R^d \times [0, t] \times [0, \Vert q \Vert_\infty]\) with intensity \(\d x \otimes \d s \otimes \d u\), independent of \(\eta\).
    We measurably enumerate 
    \begin{align*}
        \eta = \sum_i \delta_{x_i^{0}},
        \quad \Pi_t = \sum_{i} \delta_{(x_i^1, s_i^1, u_i)}.
    \end{align*}
    Now independently mark with independent unit-exponential lifespans
    \begin{align*}
        \widehat{\eta} := \sum_i \delta_{(x_i^0, R_i^0)},
        \quad \mathcal{N}_{\vert \R^d \times [0,t] \times [0, \infty) \times [0, \Vert q\Vert_\infty]} := \sum_{i} \delta_{(x_i^1, s_i^1, R_i^1, u_i)}.
    \end{align*}
    and combine the marks into one countable family \((x_i, s_i, R_i)_i\), with \(s_i = 0\) for the initial particles.

    Consider 
    \begin{align*}
        \widetilde{R}_i := R_i {\1}_{s_i + R_i \leq t}  + \infty {\1}_{s_i + R_i > t},
    \end{align*} which forgets the exact lifespan of particles which are still alive after time \(t\),
    and define the sigma-algebra
    \begin{align*}
        \widetilde{G}_t
        := \sigma(\eta, \Pi_t) \lor \sigma(\widetilde{R}_i \,\colon\, i) \lor \mathfrak{N}.
    \end{align*} Notice that in particular \(\sigma(X_s \colon 0 \leq s \leq t) \subseteq \widetilde{\mathcal{G}}_t\) by the construction of \(\mathcal{X}\) using Picard-iterates from \zcref{proposition:deterministic_solution_maps}.
    Conditional on \(\widetilde{G}_t\), the lifespans \((R_i)_{i}\) remain independent, with conditional laws
    \begin{align}\label{equation:conditional_laws}
        \mathrm{Law}(R_i \,\vert\, \widetilde{G}_t)(\d r)
        = \begin{cases}
            \delta_{\widetilde{R}_i}(\d r), & \widetilde{R}_i < \infty, \\
            {\1}_{r > t - s_i} \e^{-(r-(t-s_i))} \d r,  &\widetilde{R}_i = \infty,
        \end{cases}
    \end{align} by the memorylessness of the exponential distribution.
    Let
    \begin{align*}
        A_i 
        := \begin{cases}
            1, &i \text{ corresponds to an initial particle } x_i \text{ of } \eta, \\ 
            1 &i \text{ corresponds to a birth proposal} (x_i, s_i, R_i, u_i) \in \Ncal \text{ with } u_i \in [0, q(x_i, X_{s_i -})], \\
            0, &\text{otherwise}
        \end{cases}
    \end{align*}
    and let
    \(I_t := \{i \,\colon\, A_i = 1, \, s_i + R_i > t\}\) be the set of indices corresponding to particles which are still alive at time \(t\). Then,
    \begin{align*}
        X_t = \sum_{i \in I_t} \delta_{x_i},
        \quad \widehat{X}_t = \sum_{i \in I_t} \delta_{(x_i, s_i + R_i - t)}.
    \end{align*}
    Let \(g \in C_c(\R^d \times [0, \infty))\), \(g \geq 0\) and \(\mathcal{K} = \{x \in \R^d \,\colon\, \exists r \geq 0 \colon  g(x, r) \neq 0\} \Subset \R^d\).
    For every \(i\), the indicator \({\1}_{i \in I_t, x_i \in \mathcal{K}}\) is \(\widetilde{G}_t\)-measurable, and almost surely only finitely many of these indicators are non-zero, so that we get by conditional independence and \eqref{equation:conditional_laws} that
    \begin{align*}
        \mathbb{E}[\e^{-\int \widehat{X}_t(\d x, \d r) g(x, r)} \,\vert\, \widetilde{G}_t]
        &= \mathbb{E}\Big[\prod_{\substack{i \in I_t, \\ x_i \in \mathcal{K}}} \e^{-g(x_i, s_i + R_i - t)} \,\Big\vert\, \widetilde{G}_t \Big] 
        = \prod_{\substack{i \in I_t, \\ x_i \in \mathcal{K}}} \int_0^\infty\d r \, \e^{-g(x_i, r)} \e^{-r}  \\
        &= \int \widehat{\alpha}(X_t, \d\widehat{\zeta}) \, \e^{-\int \widehat{\zeta}(\d x, \d r) g(x, r)}.
    \end{align*}
    As countably many functions of this form determine the conditional law, by further conditioning on \(\sigma(X_s \colon 0 \leq s \leq t) \subseteq \widetilde{\mathcal{G}}_t\), we get \(\mathrm{Law}(\widehat{X}_t \,\vert\, \sigma(X_s \colon 0 \leq s \leq t)) = \widehat{\alpha}(X_t, \cdot)\).
\end{proof}

\begin{proof}[Proof of \zcref{proposition:birth_death_dynamics_via_solution_map}]
    \emph{Strong solutions, pathwise uniqueness:} Consider an arbitrary filtered probability space as in \zcref{definition:SDE}. Almost surely, by \zcref{lemma:almost_surely_good_birth_proposals}, \(\Ncal\) takes values in \(\mathtt{N}^\ast\), and \(\widehat{\eta}\) takes values in \(\widehat{\Omega}^\ast\) by the assumed exponential distribution. Consequently, \zcref{proposition:deterministic_solution_maps} shows then that 
    \begin{align*}
        X := \mathcal{X}(\widehat{\eta}, \Ncal)
    \end{align*} satisfies \eqref{equation:SDE} almost surely simultaneously for every \(B \Subset \R^d\) and \(t \geq 0\). Furthermore, \(X\) is a strong solution. Indeed, \zcref{proposition:deterministic_solution_maps}\textit{(ii)} gives
    \begin{align}\label{equation:strong_solution_by_solution_map}
        X_t 
        = \mathcal{X}(\widehat{\eta}, \Ncal)_t
        = \mathcal{X}(\widehat{\eta}, \Ncal_{[0, t]})_t, \quad \text{ almost surely}, t \geq 0,
    \end{align} and
    the right-hand side of \eqref{equation:strong_solution_by_solution_map} is measurable w.r.t.\ \(\sigma(\widehat{\eta}, \Ncal_{[0, t]})\), so that \(X_t\) is measurable w.r.t\ the completion \(\sigma(\widehat{\eta}, \Ncal_{[0, t]}) \lor \mathfrak{N}\).

    If \(Y\) is any other solution to \eqref{equation:SDE} on this filtered probability space (with the same initial condition), deterministic uniqueness in \zcref{proposition:deterministic_solution_maps} gives \(Y = X\) almost surely as \(E\)-valued random variables, so that pathwise uniqueness holds.

    \emph{Preservation of simplicity:} Suppose \(\mu(\Omega_{\mathrm{s}}) = 1\) for the distribution \(\mu\) of \(\eta\) now. Since almost surely, simultaneously for all \(x \in \R^d\), \(\Ncal(\{x\} \times [0, \infty)^3) \leq 1\) and also, by independence, \(\Ncal(\{x\} \times [0, \infty)^3) = 0\) for all \(x \in \eta\), we have the domination by an almost-surely simple-valued random variable
    \begin{align}\label{equation:stochastic_domination_simple_random_variable}
        X_t = \mathcal{X}(\widehat{\eta}, \Ncal)_t
        \leq \eta + \int_{\R^d \times (0, t] \times [0, \infty) \times [0, \Vert q\Vert_\infty]} \Ncal(\d x, \d s, \d r, \d u) \, \delta_x,
    \end{align} simultaneously for all \(t \geq 0\).
    Hence, almost surely, \(t \mapsto X_t\) is a path in \(\Omega_{\mathrm{s}}\). It is c\`adl\`ag in \(\Omega_{\mathrm{s}}\) too, since \(t \mapsto X_t\) is c\`adl\`ag in \(\Omega\) and the limits in the vague topology \(X_{t-}\) (and \(X_{t+} = X_t\)) are elements of \(\Omega_{\mathrm{s}}\) (which is not closed in \(\Omega\)) too, e.g.\ by the domination \eqref{equation:stochastic_domination_simple_random_variable} or already by the property that \(t \mapsto X_t\) has locally only finitely many jumps.
    We obtain an indistinguishable version of \(X\) as element in \(E_{\mathrm{s}}\) by modifying \(X\) on a null set to be e.g.\ the constant path \(t \mapsto \emptyset \in \Omega_{\mathrm{s}}\). The resulting process is still a strong solution to \zcref{equation:SDE} since, by assumption, the filtration contains all null sets.

    \emph{Markov property}:
    As any other solution to \eqref{equation:SDE} is indistinguishable by the step above, we may work with
    \begin{align*}
        X := \mathcal{X}(\widehat{\eta}, \Ncal)
    \end{align*} as above.
    For \(s \geq 0\) let \(\Ncal^{(s)}\) denote the birth proposals after time \(s\), with their times shifted by \(s\), as in \zcref{lemma:deterministic_restart}. The deterministic restart identity from  \zcref{lemma:deterministic_restart} then gives
    \begin{align}
        X_{s+t}
        = \mathcal{X}(\widehat{X}_s,  \mathcal{N}^{(s)})_t.
    \end{align} 
    It follows for bounded, measurable \(f \colon \Omega \to \R\), because \(\Ncal^{(s)}\) is independent from \(\mathcal{G}_s\) and has the same law as \(\mathcal{N}\) and because \(\widehat{X}_s\) is measurable w.r.t.\ \(\mathcal{G}_s\), that
    \begin{align*}
        \mathbb{E}[f(X_{s+t}) \,\vert\, \mathcal{G}_s]
        = \mathbb{E}[f(\mathcal{X}(\widehat{X}_s,  \mathcal{N}^{(s)})_t)  \,\vert\, \mathcal{G}_s]
        = \mathbb{E}[f(\mathcal{X}(\widehat{\zeta},  \mathcal{N})_t)]_{\vert \widehat{\zeta} = \widehat{X}_s}.
    \end{align*}
    Conditioning further on \(\sigma(X_r \,\colon 0 \leq r \leq s)\) and \zcref{lemma:memorylessness} give
    \begin{align}
        &\mathbb{E}[f(X_{s+t}) \,\vert\, \sigma(X_r \,\colon\, 0 \leq r \leq s)]
        = \mathbb{E}[\mathbb{E}[f(\mathcal{X}(\widehat{\zeta},  \mathcal{N})_t)]_{\vert \widehat{\zeta} = \widehat{X}_s}\,\vert\, \sigma(X_r \,\colon\, 0 \leq r \leq s)] \\
        &= \int \widehat{\alpha}(X_s, \d\widehat{\zeta}) \mathbb{E}[f(\mathcal{X}(\widehat{\zeta},  \mathcal{N})_t)]. \nonumber
    \end{align}
    Hence, the Markov property of \((X_t)_{t \geq 0}\) w.r.t.\ its natural filtration is established, and also the stronger statement of \((X_t^\eta)_{t \geq 0}\), \(\eta \in \Omega\), constituting a time-homogeneous Markov family, since then
    \begin{align*}
        \mathbb{E}[f(X_{s+t}^\eta) \,\vert\, \sigma(X_r^\eta \,\colon\, 0 \leq r \leq s)]
        = \int\widehat{\alpha}(X_s^\eta, \d\widehat{\zeta}) \, \mathbb{E}[f(\mathcal{X}(\widehat{\zeta},  \mathcal{N})_t) \nonumber
        = \mathbb{E}[f(X^\zeta_t)]_{\vert \zeta = X_s^\eta}.
    \end{align*}

    \emph{Jointly measurable realization:} We can take \((R_i)_{i \in \N}\) to be \((0, \infty)^\N\)-valued, so that then the map
    \begin{align*}
        (\eta, \xi) \mapsto \mathcal{X}(\widehat{h}(\eta, (R_i(\xi))_{i \in \N}), \Ncal(\xi)), 
        \quad \eta \in \Omega, \xi \in \Xi,
    \end{align*} is jointly measurable by \zcref{lemma:exponential_lifespan_marking} and \zcref{proposition:deterministic_solution_maps}.
    For every \(\eta \in \Omega\), \(\widehat{\eta} = \widehat{h}(\eta, (R_i)_{i \in \N})\) is a marking of \(\eta\) with independent unit-exponentials by \zcref{lemma:exponential_lifespan_marking}. It is \(\Gcal_0\)-measurable by construction and \(\Ncal\) is compatible with \((\Gcal_t)_{t \geq 0}\) in the sense of \zcref{definition:SDE} too, due to its independent increments and independence of \((R_i)_{i \in \N}\).
    Therefore, we can apply \textit{(i)} of this proposition to see that \(X^\eta = \mathcal{X}(\widehat{\eta}, \Ncal)\) is a strong solution of \eqref{equation:SDE} with initial condition \(\eta\).
    Finally, we have \(\mathbb{P}(\Ncal \in \mathtt{N}^\ast) = 1\) by \zcref{lemma:almost_surely_good_birth_proposals}, so that \zcref{proposition:deterministic_solution_maps}\textit{(i)} shows that on this event, \(X^\eta\) satisfies \eqref{equation:SDE} simultaneously for all \(\eta \in \Omega\) and all \(B \Subset \R^d\), \(t \geq 0\).
\end{proof}

We are now in a position to give a
\begin{proof}[Proof of \zcref{proposition:SDE}\textit{(i),(ii),(iv),(v),(vi)}]

    \textit{(i)}: Existence of a strong solution and pathwise uniqueness are both proved in \zcref{proposition:birth_death_dynamics_via_solution_map}\textit{(i)}.

    \textit{(ii)}: Uniqueness in law is most easily seen again from \zcref{proposition:birth_death_dynamics_via_solution_map}\textit{(i)}. The law of any solution is that of \(\mathcal{X}(\widehat{\eta}, \mathcal{N})\) and this is evidently independent of the underlying filtered probability space.

    \textit{(iv)}: The Markov property is proved in \zcref{proposition:birth_death_dynamics_via_solution_map}\textit{(ii)}. 

    \textit{(v)}: This is also proved in \zcref{proposition:birth_death_dynamics_via_solution_map}\textit{(i)}.

    \textit{(vi)}: This follows from uniqueness in law of the solution established in \textit{(ii)}. Indeed, fix \(z \in \R^d\) and let \(X\) be any solution of \zcref{equation:SDE} with initial distribution \(\mu\), driven by \(\widehat{\eta}\) and \(\mathcal{N}\). Set \(\widetilde{X} := \theta_z X\), \(\widetilde{\widehat{\eta}} = \theta_z \widehat{\eta}\) and \(\widetilde{\mathcal{N}} := \theta_z \mathcal{N}\). Then, by translation invariance of \(q\), \(\widetilde{X}\) satisfies \zcref{equation:SDE} but driven by \(\widetilde{\widehat{\eta}}\) and \(\widetilde{\mathcal{N}}\). Since also \(\widetilde{X}_0 \sim \mu \circ \theta_z^{-1}\), uniqueness in law gives
    \begin{align*}
        P^{\mu} \circ \theta_z^{-1}
        =\mathrm{Law}(X) \circ \theta_z^{-1}
        = \mathrm{Law}(\widetilde{X})
        = P^{\mu \circ \theta_z^{-1}}. 
    \end{align*}
\end{proof}

\subsection{The semigroup and forward equation}

\begin{lemma}\label{lemma:infinite_volume_process_forward_equation_local_bounded_measurable}
    Let \(f \in \Lcal\).
    Then, for every \(\eta \in \Omega\), \(t \mapsto T_t f(\eta)\) is continuously differentiable on \([0, \infty)\) with derivative
    \begin{align*}
        t \mapsto T_t (\mathscr{L} f)(\eta).
    \end{align*} The derivative in \(t = 0\) is interpreted as right-derivative.
    In particular, the forward equation
    \begin{align}
        T_t f(\eta)
        = f(\eta) + \int_0^t \, T_s(\mathscr{L} f)(\eta) \, \d s, \quad t \geq 0,
    \end{align}
    holds.
\end{lemma}
\begin{proof}
     \newcommand{\deathsinLfrometa}[1]{D^{\eta}_\L(#1)}
    \newcommand{\deathsinLfrombirths}[1]{D^{2}_\L(#1)}

    Let \(f\) be bounded and \(\Lambda\)-measurable, \(\Lambda \Subset \R^d\). Fix \(\eta \in \Omega\) and consider the solution-map representation \(X = X^\eta = \mathcal{X}(\widehat{\eta}, \Ncal)\) from \zcref{proposition:birth_death_dynamics_via_solution_map}.
    
    We denote by \(\birthsinL{t}\), resp.\ \(\deathsinL{t}\), the set of times \(s \leq t\) in which a birth, resp.\ death, occurs in \(\L\). Furthermore, we split up \(\deathsinL{t}\) in the deaths \(\deathsinLfrometa{t}\) coming from the initial condition \(\widehat{\eta}\) and the deaths \(\deathsinLfrombirths{t}\) of particles which are not coming from the initial condition but are born in \([0, t)\). Also write
    \begin{align*}
        \eta_t := \int_{\R^d \times [0, \infty)} \widehat{\eta}(\d x, \d r) {\1}_{r > t} \delta_x,
        \quad Y_t := X_t - \eta_t.
    \end{align*}
    We will sometimes use the short-hand
    \begin{align*}
        D_x^{-}f(\zeta)
        = \begin{cases}
            f(\zeta - \delta_x) - f(\zeta), &\zeta(\{x\}) \geq 1, \\
            0, &\text{otherwise}.
        \end{cases}
    \end{align*}

    We decompose into the individual jump contributions
    \begin{align*}
        &f(X_t) - f(\eta) \\
        &= \sum_{s \in \birthsinL{t}}(f(X_s) - f(X_{s-})) + \sum_{s \in \deathsinLfrombirths{t}}(f(X_s) - f(X_{s-}))
        + \sum_{s \in \deathsinLfrometa{t}}(f(X_s) - f(X_{s-})).
    \end{align*}
    Taking expectations gives
    \begin{align}
        (T_t f)(\eta) - f(\eta)
        = C_B(t) + C_{D^2}(t) + C_{D\eta}(t),
    \end{align} where
    \begin{align*}
        &C_B(t) 
        := \mathbb{E}\Big[\sum_{s \in \birthsinL{t}} (f(X_s) - f(X_{s-}))\Big]
    \end{align*} and  \(C_{D^2}(t)\), \(C_{D\eta}(t)\) are defined analogously.
    Notice that each term is integrable, as the number of jumps in \(\Lambda\) during \([0, T]\) is dominated by \(2 \Ncal(\Lambda \times [0, T] \times [0, \infty) \times [0, \Vert q\Vert_\infty]) + N_{\Lambda}(\eta)\).

    \textit{1. Birth term \(C_B(t)\):} 
    Since \(X\) solves the graphical SDE \eqref{equation:SDE}, applying the Mecke formula yields
    \begin{align}\label{equation:formula_C_B}
        C_B(t) &= \mathbb{E}\Big[\int_{\Lambda \times (0, t] \times [0, \infty) \times [0, \Vert q\Vert_\infty]} \Ncal(\d x, \d s, \d r, \d u) [f(X_{s-} + \delta_x) - f(X_{s-})] {\1}_{u \leq q(x, X_{s-})} \Big] \\
        &= \int_0^t \d s \int_\Lambda \d x \int_0^\infty \d r \, \e^{-r} \int_0^{\Vert q\Vert_\infty} \d u \mathbb{E}\big[ [f(X_{s-} + \delta_x) - f(X_{s-})] {\1}_{u \leq q(x, X_{s-})} \big] \nonumber \\
        &= \int_0^t \d s \, \mathbb{E}\Big[\int_\Lambda \d x \, q(x, X_s) [f(X_s + \delta_x) - f(X_s)]  \Big], \nonumber
    \end{align} 
    since
    \begin{align*}
        \mathcal{X}(\widehat{\eta}, \mathcal{N} + \delta_{(x, s, r, u)})_{s-} = \mathcal{X}(\widehat{\eta}, \mathcal{N})_{s-} = X_{s-}
    \end{align*} 
    by non-anticipativity of the solution map, and replacing \(X_{s-}\) by \(X_{s}\) doesn't change the Lebesgue integral over \(s\).

    \textit{2. Second-order death term \(C_{D^2}(t)\):} Again by the graphical construction and the Mecke equation,
    \begin{align}\label{equation:C_D_square}
        C_{D^2}&(t)
        = \mathbb{E}\Big[\int_{\Lambda \times (0, t] \times [0, \infty) \times [0, \Vert q\Vert_\infty]} \Ncal(\d x, \d s, \d r, \d u) \, D_x^{-}f(X_{(s+r)-}) {\1}_{u \leq q(x, X_{s-})} {\1}_{s + r \leq t} \Big] \\
        &= \int_0^t \d s \int_\Lambda \d x \int_0^\infty \d r \, \e^{-r} \int_0^{\Vert q\Vert_\infty} \d u \, \mathbb{E}\big[ D_x^{-}f(\mathcal{X}(\widehat{\eta}, \Ncal + \delta_{(x,s,r,u)})_{(s+r)-}) {\1}_{u \leq q(x, X_{s-})} {\1}_{s + r \leq t}\big]. \nonumber
    \end{align}
    Taking right difference quotients yields for \(t > 0\),
    \begin{align*}
        \frac{\d}{\d t^+}C_{D^2}(t)
        = \int_0^t \d s \, \e^{-(t-s)} \int_\Lambda \d x \int_0^{\Vert q\Vert_\infty} \d u \, \mathbb{E}\big[ D_x^{-}f(\mathcal{X}(\widehat{\eta}, \Ncal + \delta_{(x,s,t-s,u)})_{t-}) {\1}_{u \leq q(x, X_{s-})} \big].
    \end{align*} which can then be rewritten using the Mecke equation, and \(\e^{-(t-s)} = \int_{0}^\infty \d r \e^{-r} {\1}_{s+r > t}\), as
    \begin{align*}
        &\frac{\d}{\d t^+}C_{D^2}(t) \\
        &= \int_0^t \d s \, \int_\Lambda \d x \int_0^\infty \d r \, \e^{-r} \int_0^{\Vert q\Vert_\infty} \d u \, \mathbb{E}\big[ D_x^{-}f(\mathcal{X}(\widehat{\eta}, \Ncal + \delta_{(x,s,t-s,u)})_{t-}) {\1}_{u \leq q(x, X_{s-})} {\1}_{s+r > t} \big] \\
        &=  \int_0^t \d s \, \int_\Lambda \d x \int_0^\infty \d r \, \e^{-r} \int_0^{\Vert q\Vert_\infty} \d u \,  \mathbb{E}\big[ D_x^{-}f(\mathcal{X}(\widehat{\eta}, \Ncal + \delta_{(x,s,r,u)})_{t}) {\1}_{u \leq q(x, X_{s-})} {\1}_{s+r > t} \big]  \\
        &= \mathbb{E}\Big[\int_{\Lambda \times (0, t] \times [0, \infty) \times [0, \Vert q\Vert_\infty]} \Ncal(\d x, \d s, \d r, \d u) \, D_x^{-}f(X_{t}) {\1}_{u \leq q(x, X_{s-})} {\1}_{s + r > t} \Big] \\
        &= \mathbb{E}\Big[\int_{\Lambda} Y_t(\d x) \, D_x^{-}f(X_t) \Big].
    \end{align*}
    It follows, using \(C_{D^2}(0) = 0\) and the absolute continuity already provided by the integral representation \eqref{equation:C_D_square}, that
    \begin{align}\label{equation_formula_C_D}
        C_{D^2}(t) = \int_0^t \d s \, \mathbb{E}\Big[\int_{\Lambda} Y_s(\d x) \, D_x^{-}f(X_s) \Big].
    \end{align}

    \textit{3. First-order death term \(C_{D^\eta}(t)\):} 
    Write
    \begin{align*}
        \eta = \sum_i \delta_{x_i},
        \quad \widehat{\eta} = \sum_i \delta_{(x_i, R_i)}
    \end{align*} where \(R_i\) are independent unit-exponential random variables, independent of \(\Ncal\).
    Then,
    \begin{align*}
        C_{D^\eta}(t)
        = \sum_{i \colon x_i \in \Lambda} \mathbb{E}[D_{x_i}^-f(X_{R_i -}) {\1}_{R_i \leq t}].
    \end{align*}
    For each such \(i\), writing
    \begin{align*}
        \widehat{\eta}^{(i, r)}
        := \sum_{j \neq i} \delta_{(x_j, R_j)} + \delta_{(x_i, r)},
    \end{align*} we get by construction of the solution map with random inputs (or in other words, conditioning and using that the lifespans of the initial condition are i.i.d.\ unit exponentials) that
    \begin{align*}
        \mathbb{E}[D_{x_i}^-f(X_{R_i -}) {\1}_{R_i \leq t}]
        = \int_0^t \d r \, \e^{-r} \mathbb{E}[D_{x_i}^{-} f(\mathcal{X}(\widehat{\eta}^{(i, r)}, \Ncal)_{r-})].
    \end{align*}
    This already shows that \(C_{D^\eta}(t)\) is locally absolutely continuous and we have
    \begin{align*}
        \frac{\d}{\d t^+} \mathbb{E}[D_{x_i}^-f(X_{R_i -}) {\1}_{R_i \leq t}]
        = \mathbb{E}[D_{x_i}^-f(X_{t}) {\1}_{R_i > t}],
    \end{align*} 
    yielding
    \begin{align}\label{equation:formula_C_D_eta}
        C_{D^\eta}(t)
        = \int_0^t \d s  \, \mathbb{E}\Big[\sum_{i \colon x_i \in \Lambda} D_{x_i}^{-}f(X_s) {\1}_{R_i > s}\Big]
        = \int_0^t  \d s  \, \mathbb{E}\Big[\int_\Lambda \eta_s(\d x) D_x^{-}f(X_s) \big]
    \end{align}

    \textit{4. Summing up the individual contributions:}  Adding everything up, \(X_s = \eta_s + Y_s\) and \eqref{equation:formula_C_B}, \eqref{equation_formula_C_D},\eqref{equation:formula_C_D_eta} imply
    \begin{align}\label{equation:integral_formula}
        &T_t f(\eta) - f(\eta) 
        = C_B(t) + C_{D^2}(t) + C_{D^\eta}(t) \\
        &=  \int_0^t \d s \, \mathbb{E}\Big[\int_\Lambda \d x \, q(x, X_s) [f(X_s + \delta_x) - f(X_s)]  \Big]
        + \int_0^t \d s \, \mathbb{E}\Big[\int_\Lambda X_s(\d x) \, [f(X_s - \delta_x) - f(X_s)] \Big]  \nonumber\\
        &= \int_0^t \d s \, \mathbb{E}[(\mathscr{L}f)(X_s)] 
        = \int_0^t \d s \, T_s(\mathscr{L}f)(\eta).  \nonumber
    \end{align}

    \textit{5. Pointwise continuous differentiability:}   
    For every fixed \(t \geq 0\), almost surely no birth or death occurs in \(\Lambda \oplus B_R(o)\) exactly at time \(t\), so that the locally finitely many birth/death jumps of \(X\) imply
    \begin{align*}
        (\mathscr{L}f)(X_s)
        \xrightarrow[s \to t]{}  (\mathscr{L}f)(X_t)
    \end{align*} almost surely.
    Dominated convergence (recall the domination of jumps from the beginning of this proof and that \(\Vert f\Vert_\infty < \infty\) by assumption) shows that \(t \mapsto T_t(\mathscr{L} f)(\eta) = \mathbb{E}[(\mathscr{L}f)(X_t)]\) is continuous, so that the asserted continuous differentiability follows from \eqref{equation:integral_formula}.
\end{proof}

We are now in a position to give a 
\begin{proof}[Proof of \zcref{proposition:transition_semigroup}]
    
    \textit{(i)}: These processes \(X^\eta\) were constructed in \zcref{proposition:birth_death_dynamics_via_solution_map}\textit{(iii)}, and the solution map construction ensures that
    \begin{align*}
        \eta \mapsto (T_t f)(\eta) 
        = \mathbb{E}[f(X_t^\eta)]
        = \mathbb{E}[f(\mathcal{X}(\widehat{h}(\eta, (R_i)_{i \in \N}), \mathcal{N})_t)]
    \end{align*} is measurable as map \(\Omega \to \R\). The time-homogeneous Markov property of \(X\) shown in  \zcref{proposition:birth_death_dynamics_via_solution_map} ensures that \((T_t)_{t \geq 0}\) is a Markov semigroup. 
    
    \textit{(ii)}: Translation covariance was already shown in \zcref{proposition:SDE}\textit{(vi)}.

    \textit{(iii)}: This is the content of \zcref{lemma:infinite_volume_process_forward_equation_local_bounded_measurable}.
\end{proof}

\subsection{The martingale problem}

Let us finally turn to the martingale problem from \zcref{definition:martingale_problem}, and hence the proofs of \zcref{proposition:martingale_problem} and \zcref{proposition:SDE}\textit{(iii)}.
We need the following two lemmas.

\begin{lemma}\label{lemma:SDE_solution_law_solves_martingale_problem}
    For every \(\mu \in \Pcal(\Omega)\), the martingale problem for \((\mathscr{L}, \mathcal{D}, \mu)\) introduced in \zcref{definition:martingale_problem} is solved by the law on \(E\) of any solution \(X\) of the SDE \eqref{equation:SDE} with initial condition \(\mu\), whose existence was established in \zcref{proposition:SDE}.
\end{lemma}
\begin{proof}
    By \zcref{remark:martingale_problem_markov_process}, this follows from the Markov property of \(X\) established in \zcref{proposition:SDE} and the forward equation in \zcref{lemma:infinite_volume_process_forward_equation_local_bounded_measurable} for \(f \in \mathcal{D} \subseteq \Lcal\).
\end{proof}

\begin{lemma}\label{lemma:uniqueness_in_martingale_problem}
    For every \(\mu \in \Pcal(\Omega)\), there is at most one solution to the martingale problem for \((\mathscr{L}, \mathcal{D}, \mu)\) introduced in \zcref{definition:martingale_problem}.
\end{lemma}
\begin{proof}
    We apply \cite[Corollary 2.11]{garcia2006spatial}, which infers uniqueness of the martingale problem for \((\mathscr{L}, \mathcal{D}, \mu)\) from weak uniqueness to the lifespan-marked SDE described in the following.

    Fix \(\mu \in \Pcal(\Omega)\).
    Let \(\widehat{\alpha}(\eta, \d\widehat{\eta})\) denote the kernel assigning independent unit-exponential lifespan marks to the particles of \(\eta\) from \zcref{lemma:exponential_lifespan_marking} and define \(\widehat{\mu} \in \Pcal(\widehat{\Omega})\) by
    \begin{align*}
        \widehat{\mu}(\d\widehat{\eta}) = \int_\Omega \mu(\d\eta) \widehat{\alpha}(\eta, \d\widehat{\eta}).
    \end{align*}
    A process \(\widetilde{Y}\) is a \emph{weak solution} to the SDE
    \begin{equation}\label{equation:lifespan_marked_SDE}
    \begin{split}
       \widehat{Y}_t
        &= \int_{\R^d \times [0, t] \times [0, \infty)^2} \mathcal{N}(\d x, \d s, \d r, \d u) {\1}_{(s+r) > t} {\1}_{[0, q(x, \widehat{Y}_{s-}(\cdot \times [0, \infty)))]}(u) \delta_{(x, (r+s) - t)}  \\
        &\qquad + \int \widehat{Y}_0(\d x , \d r) {\1}_{r > t} \delta_{(x, r-t)}
    \end{split}
    \end{equation} 
    if there exists a filtered probability space \((\Xi, \mathcal{G}, (\mathcal{G}_t)_{t \geq 0}, \mathbb{P})\) supporting the following random variables:
    \begin{enumerate}[1.]
        \item A Poisson random measure \(\Ncal\) on \(\R^d \times [0, \infty)^3\) with intensity \(\d x \otimes \d s \otimes\e^{-r}\d r \otimes \d u\) compatible with the filtration,
        \begin{align*}
            \sigma(\Ncal_{\R^d \times [0, t] \times [0, \infty)^2}) \subseteq \mathcal{G}_t, \quad
            \sigma(\Ncal_{\R^d \times (t, \infty) \times [0, \infty)^2}) \perp\!\!\!\perp \mathcal{G}_t, \quad t \geq 0.
        \end{align*}
        \item An \(\widehat{\Omega}\)-valued, \((\mathcal{G}_t)_{t \geq 0}\)-adapted process \(\widehat{Y} = (\widehat{Y}_t)_{t \geq 0}\) with  c\`adl\`ag paths such that \(\widehat{Y}_0 \sim \widehat{\mu}\) and such that \zcref{equation:lifespan_marked_SDE} holds \(\mathbb{P}\)-a.s.\ simultaneously for every \(t \geq 0\), and such that \(\widetilde{Y}\) and \(\widehat{Y}\) have the same law as random variables on \(\widehat{E}\).
    \end{enumerate}    
    We say that \emph{weak uniqueness} holds for \zcref{equation:lifespan_marked_SDE} if the laws of any two weak solutions to \eqref{equation:lifespan_marked_SDE} as random variables on \(\widehat{E}\) coincide.

    In our language, \cite[Corollary 2.11]{garcia2006spatial} then implies that there is at most one solution in the martingale problem for \((\mathscr{L}, \mathcal{D}, \mu)\) from \zcref{definition:martingale_problem} in case weak uniqueness holds for the SDE \eqref{equation:lifespan_marked_SDE} with initial condition \(\mu\).
    Let us now proceed by first checking in their notation that our setting is covered by their results, and then go on and demonstrate that weak uniqueness.

    \textit{Our setting is covered.}
    In the notation of \cite[Section 2]{garcia2006spatial}, we have
    \begin{align}
        S := \R^d,
        \quad \beta(\d x) := \d x,
        \quad \lambda(x, \eta) := q(x, \eta),
        \quad \delta(x, \eta) := 1.
    \end{align}
    Furthermore, we choose, also in their notation, the closed balls \(K_k := \{x \in \R^d \,\colon\, \abs{x} \leq k\}\), \(k \in \N\), and non-negative \(c_k \in C_c(\R^d)\) with
    \begin{align*}
        c_k = 1 \text{ on } K_k, \quad c_1 \leq c_2 \leq \dots.
    \end{align*}
    They then choose \(\mathcal{C}\) to be the set of all bounded continuous functions \(f\) on \(S\) with \(\abs{f} \leq a c_k\) for some \(a > 0\) and \(k \in \N\).
    For us this reads \(\mathcal{C} = C_c(\R^d)\), and their \(\mathcal{S}\) is, in our context,
    \begin{align}
        \mathcal{S} = \Omega, \text{ equipped with the vague topology}.
    \end{align}
    We now show that their standing assumptions hold in our setting.
    Their \cite[Condition 2.1]{garcia2006spatial} holds, as even
    \begin{align}
        \sup_{\eta \in \Omega} \int c_k(x) q(x, \eta) \d x
        \leq \Vert q\Vert_\infty \int c_k(x) \d x 
        < \infty, \quad k \in \N,
    \end{align} and \(\delta \equiv 1\).
    For their \cite[Condition 2.2]{garcia2006spatial}, suppose \(\eta \in \Omega\) and let \((\eta_n)_{n \in \N}\) be a sequence in \(\Omega\) such that \(\lim_{n \to \infty} \int_{\R^d} c_k(x) \abs{\eta_n - \eta}(\d x) = 0\) for all \(k \in \N\).
    Then,
    \begin{align}
        \abs{\eta_n - \eta}(K_k)
        \leq \int_{\R^d} c_k(x) \abs{\eta_n - \eta}(\d x)
        \xrightarrow[n \to \infty]{} 0,
    \end{align} 
    so that the \(\eta\) and \(\eta_n\) eventually agree on \(K_k\). In particular, since for every \(x \in \R^d\) there is a \(k \in \N\) such that \(K_k \supseteq B_R(x)\), we have that
    \begin{align}
        q(x, \eta)
        = q(x, \eta_{B_R(x)})
        = \lim_{n \to \infty} q(x, (\eta_n)_{B_R(x)})
        = \lim_{n \to \infty} q(x, \eta_n).
    \end{align} 
    Since \(\delta \equiv 1\), the equality \(\delta(x, \eta) = \lim_{n \to \infty} \delta(x, \eta_n)\) is automatic.

    Finally, we note the following important observation about the applicability of \cite[Corollary 2.11]{garcia2006spatial}.
    The authors never detail the specific martingale problem they consider for the generator \(\mathscr{L}\), but their proofs of \cite[Proposition 2.9, Theorem 2.10]{garcia2006spatial} apply to our \zcref{definition:martingale_problem} of the martingale problem for \((\mathscr{L}, \mathcal{D}, \mu)\).
    Indeed, let, as right before \cite[Section 2.1]{garcia2006spatial}, \(\widehat{C}\) be the set of bounded and continuous functions \(f \colon \R^d \times [0, \infty) \to \mathbb{R}\) such that:
    \begin{enumerate}[1.]
        \item The partial derivative \(\partial_r f(x, r)\) in the second variable exists and \(\partial_r f\) is a bounded and continuous function \(\R^d \times [0, \infty) \to \mathbb{R}\);
        \item \(f(\cdot, 0) = 0\);
        \item The two functions \(\sup_r \abs{f(\cdot, r)}\), \(\sup_r \abs{\partial_r f(\cdot, r)}\) are elements of \(\mathcal{C} = C_c(\R^d)\);
        \item There exists \(r_f > 0\) such that \(\partial_r f(\cdot, r) = 0\) for \(r > r_f\).
    \end{enumerate}
    Then, the domain directly suggested by the proof of \cite[Proposition 2.9]{garcia2006spatial} is
    \begin{align*}
        \widetilde{\mathcal{D}}
        = \{\eta \mapsto \e^{-\int \eta(\d x) g(x)} \,\colon\, \exists f \in \widehat{C}, \, f \geq 0 \text{ such that } \e^{-g(x)} = \int_0^\infty \, \e^{-f(x, r)} \e^{-r} \, \d r \text{ for all } x \in \R^d \}.
    \end{align*}
    Let us show that
    \begin{align*}
        \widetilde{\mathcal{D}}
        \subseteq \{\eta \mapsto \e^{-\int \eta(\d x) g(x)} \,\colon\, g \in C_c(\R^d),\, g \geq 0\}
        = \mathcal{D},
    \end{align*}
    with \(\mathcal{D}\) from \zcref{definition:martingale_problem} of the martingale problem for \((\mathscr{L}, \mathcal{D}, \mu)\). In fact, \(\widetilde{\mathcal{D}} = \mathcal{D}\), but of course uniqueness in the martingale problem with a smaller domain already implies uniqueness in the martingale problem with a bigger domain.
    
    Let \(f \in \widehat{C}\), \(f \geq 0\) and \(g \colon \R^d \to \R\) such that
    \begin{align*}
        \e^{-g(x)} = \int_0^\infty \, \e^{-f(x, r)} \e^{-r} \, \d r, \quad x \in \R^d.
    \end{align*}
    Then \(g\) satisfies
    \begin{align*}
        g(x)
        = -\log \int_0^\infty \, \e^{-f(x, r)} \e^{-r} \, \d r \,\,\geq\,\, 0, \quad x \in \R^d,
    \end{align*} and is evidently continuous by dominated convergence. Furthermore, it is compactly supported by the assumption that \(\sup_r \abs{f(\cdot, r)}\) is compactly supported. Hence \(g \in C_c(\R^d)\), \(g \geq 0\).

    \textit{Weak uniqueness for the lifespan-marked SDE.}
    Consider any weak solution \(\widetilde{Y}\) of \eqref{equation:lifespan_marked_SDE} and a corresponding \(\widehat{Y}\) from the definition of weak solution.
    Complete the considered probability space and augment the filtration with null-sets to match our \zcref{definition:SDE}. This does not change the law of any of the involved variables, nor the validity of the SDE, and adaptation to the filtration as well as compatibility of \(\Ncal\) with it is preserved.
    
    Write \(Y = (Y_t)_{t \geq 0}\) for an indistinguishable, adapted version of  the unmarked projection \(\widehat{Y}(\cdot \times [0, \infty))\) which has c\`adl\`ag paths everywhere, obtained by modifying \(\widehat{Y}(\cdot \times [0, \infty))\) on a \(\mathcal{G}_0\)-null set. This is possible as on every spatial region \(\Lambda \Subset \R^d\) and on every finite time interval \([0, T]\), there are only finitely many initial particles and \(\mathcal{N}\) has there, on the probability one event \(\mathtt{N}^\ast\), only finitely many birth proposals \((x, s, r, u)\) with \(u \in [0, \Vert q\Vert_\infty]\).
    
    Then, \(Y\) solves the SDE \eqref{equation:SDE} and
    \begin{align}\label{equation:solution_of_SDE_is_almost_surely_solution_map_with_random_input}
        Y = \mathcal{X}(\widehat{Y}_0, \mathcal{N})
    \end{align} 
    holds almost surely by \zcref{proposition:birth_death_dynamics_via_solution_map}\textit{(i)}.
    Since \(\widehat{Y}\) solves the SDE \eqref{equation:lifespan_marked_SDE}, substituting \eqref{equation:solution_of_SDE_is_almost_surely_solution_map_with_random_input} into the definition of \(\widehat{\mathcal{X}}\), we have
    \begin{align}
        \widehat{Y}
        = \widehat{\mathcal{X}}(\widehat{Y}_0, \mathcal{N})
    \end{align} 
    almost surely.
    Hence, the law of \(\widehat{Y}\) is determined by the law of \((\widehat{Y}_0, \mathcal{N})\) alone, and that is \(\widehat{\mu} \otimes \mathrm{Law}(\mathcal{N})\), independently of the concrete choice of filtered probability space.
\end{proof}
\begin{remark}
    \begin{enumerate}
        \item Notice that we do not need to distinguish between local martingale problem and martingale problem here, by the chosen test function class \(\mathcal{D}\) and by the bounded birth rate \(q\). For \(f \in \mathcal{D}\) it holds \(f(\eta) = \e^{-\int \eta(\d x) g(x)}\) for some \(g \in C_c(\R^d)\) with \(g \geq 0\), so that we have
        \begin{align*}
            (\mathscr{L}f)(\eta)
            =  \int \, q(x, \eta) \, f(\eta)[\e^{-g(x)} - 1] \, \d x  + \int \, f(\eta)[\e^{g(x)} - 1] \, \eta(\d x)
        \end{align*} and therefore
        \begin{align*}
            \vert (\mathscr{L}f)(\eta) \vert
            &\leq 
            \Vert q\Vert_\infty \int [1-\e^{-g(x)}] \d x + \int \, \e^{-\int \eta(\d y) g(y)} \abs{\e^{g(x)} - 1} \, \eta(\d x) \\
        \end{align*} with
        \begin{align*}
            &\int \, \e^{-\int \eta(\d y) g(y)} \abs{\e^{g(x)} - 1} \, \eta(\d x)
            \leq \int \, \e^{-\int \eta(\d y) g(y)} g(x) \e^{g(x)} \, \eta(\d x)  \\
            &\leq \e^{\Vert g\Vert_\infty} (\textstyle\int \eta(\d x) g(x))\e^{-\int \eta(\d x) g(x)} 
            \leq \e^{\Vert g\Vert_\infty -1 }. 
        \end{align*}
        Hence,
        \begin{align*}
            \sup_{0 \leq t \leq T} \abs{M_t^f}
            \leq 1 + T \Vert \mathscr{L} f\Vert_\infty
            < \infty,
        \end{align*} so that \((M_t^f)_{t \geq 0}\) is automatically a martingale if it is a local martingale.

        \item Alternatively, in the proof above, it would be possible to derive weak uniqueness of the SDE \eqref{equation:lifespan_marked_SDE} through (the proof) of \cite[Theorem 2.13]{garcia2006spatial}. Its additional assumptions would be fulfilled in our setting if we let \(a(x, y) := \Vert q \Vert_\infty {\1}_{\abs{x-y} \leq R}\), \(c \equiv 1\).

        \item Regarding \(\mathcal{D} \subseteq \widetilde{\mathcal{D}}\) in the proof above, make the product ansatz \(f(x, r) = h(x) \varphi(r)\).
        Let \(\varphi \in C^\infty([0, \infty))\) increase strictly from \(\varphi(0) = 0\) to \(\varphi(1) = 1\) and with \(\varphi(r) = 1\) for \(r > 1\). 
        Consider the smooth and strictly increasing function \(\Phi \colon [0, \infty) \to [0, \infty)\) given by
        \begin{align*}
            \Phi(y)
            = - \log \int_0^\infty \e^{-y \varphi(r)} \e^{-r} \d r, \quad y \geq 0.
        \end{align*} For given \(g \in C_c(\R^d)\), \(g \geq 0\), set
        \begin{align*}
            h(x) := \Phi^{-1}(g(x)),
            \quad f(x, r) := h(x) \varphi(r).
        \end{align*}
        Then, by construction, \(f \in \widehat{C}\) and \(\e^{-g(x)} = \int_0^\infty \, \e^{-f(x, r)} \e^{-r} \, \d r\) for all \(x \in \R^d\).
    \end{enumerate}
\end{remark}

We can now give a 
\begin{proof}[Proof of \zcref{proposition:martingale_problem}]
    The combination of \zcref{lemma:SDE_solution_law_solves_martingale_problem} and \zcref{lemma:uniqueness_in_martingale_problem} gives existence and uniqueness in the martingale problem for \((\mathscr{L}, \mathcal{D}, \mu)\), for any \(\mu \in \Pcal(\Omega)\). That any process with law \(P^\mu\) is a time-homogeneous Markov process follows from the fact that \zcref{lemma:SDE_solution_law_solves_martingale_problem} identifies \(P^\mu\) with the law of a solution \(X\) of the SDE \eqref{equation:SDE}, and \(X\) has this Markov property by \zcref{proposition:SDE}\textit{(iv)}.
\end{proof}

Finally, we also present a
\begin{proof}[Proof of \zcref{proposition:SDE}\textit{(iii)}]
    This is just the combination of \zcref{lemma:SDE_solution_law_solves_martingale_problem} and \zcref{lemma:uniqueness_in_martingale_problem}.
\end{proof}

\section{Finite-volume dynamics}\label{sec:appendix:finite_volume_dynamics}

In this section we construct the finite-volume dynamics and then prove the statements of \zcref{section:finite_volume_dynamics} relating to its semigroup properties and the comparison with the original dynamics.

\subsection{Construction and coupling with global dynamics}\label{section:finite_volume_dynamics_construction}

We will now provide a proof of \zcref{proposition:SDE_local} using the already established infinite-volume results from \zcref{proposition:SDE}. To do so, we embed the finite-volume dynamics into an infinite-volume dynamics as follows.

Consider the auxiliary infinite-volume birth rate \(\widetilde{b}_\Lambda\) given by
\begin{align}
    \widetilde{b}_\Lambda(x, \eta)
    = b_\Lambda^\nu(x, \eta_\Lambda) {\1}_{\Lambda}(x), \quad x \in \R^d, \eta \in \Omega.
\end{align}
Then, with birth rate \(q = \widetilde{b}_\Lambda\), the infinite-volume SDE \eqref{equation:SDE}, reads
\begin{equation}\label{equation:finite_volume_SDE_on_Omega_frozen_exterior}
\begin{split}
    X^{\Lambda, \eta}_t(B)
    = \int_{B \times (0, t] \times [0, \infty)^2} \Ncal(\d x, \d s, \d r, \d u) &{\1}_{(t-s, \infty)}(r) {\1}_{[0, \widetilde{b}_\Lambda(x, (X^{\Lambda, \eta}_{s-})_\Lambda)]}(u) \\
    &+ \int_{B \times [0, \infty)} \widehat{\eta}(\d x, \d r) {\1}_{(t, \infty)}(r),
\end{split}
\end{equation} for \(B \Subset \R^d\) and \(t \geq 0\), with \(\widehat{\eta}\) a marking of \(\eta\) with i.i.d.\ unit exponentials and an independent Poisson random measure \(\Ncal\) on \(\R^d \times [0, \infty)^3\) with intensity \(\d x \otimes \d s \otimes \e^{-r} \d r \otimes \d u\).
Outside \(\Lambda\), a solution \(X^{\Lambda, \eta}\) consists only of the surviving initial particles of \(\widehat{\eta}\). Its restriction to \(\Lambda\),
\begin{align}
    X_t^{(\Lambda, \eta)}
    := (X_t^{\Lambda, \eta})_\Lambda,
\end{align} 
satisfies the SDE \eqref{equation:SDE_local} with initial configuration \(\eta_\Lambda \in \Omega_\Lambda\) and initial lifespans given by the restriction of \(\widehat{\eta}\) to \(\Lambda \times [0, \infty)\).
Conversely, every solution to the SDE \eqref{equation:SDE_local} lifts to a solution of \eqref{equation:finite_volume_SDE_on_Omega_frozen_exterior} by augmenting with the surviving particles of \(\widehat{\eta}\) in \(\Lambda^c\).

We will later need, for any initial condition \(\eta \in \Omega\), not only existence of solutions to the SDE \eqref{equation:SDE_local}, but a coupling between the infinite-volume Markov process \((X_t\Ssup{\eta})_{t \geq 0}\) and the Markov process \((X_t\Ssup{\Lambda, \eta})_{t \geq 0}\) in the finite volume \(\Lambda\) with stochastic boundary conditions, using a common Poisson random measure \(\mathcal{N}\) for the birth proposals and common i.i.d.\ exponentially distributed lifespans for shared particles of the initial configuration.
\begin{proposition}[Coupling of local and global dynamics via solution maps]
    There exists a measurable map \(\mathcal{X}^{\Lambda} \colon \widehat{\Omega} \times \mathtt{N} \to E\) with the following properties.
    On the probability space of \zcref{proposition:birth_death_dynamics_via_solution_map}\textit{(iii)}, set
    \begin{align*}
        X^{\Lambda, \eta}
        := \mathcal{X}^\Lambda(\widehat{h}(\eta, (R_i)_{i \in \N}), \Ncal) \quad \eta \in \Omega.
    \end{align*}
    Then, the map \((\eta, \xi) \mapsto X^{\Lambda, \eta}(\xi)\) is jointly measurable \(\Omega \times \Xi \to E\) and there is an event of probability one on which, simultaneously for every \(\eta \in \Omega\), \(X^{\Lambda, \eta}\) satisfies \eqref{equation:finite_volume_SDE_on_Omega_frozen_exterior} with \(\widehat{\eta} = \widehat{h}(\eta, (R_i)_{i \in \N})\) simultaneously for all \(B \Subset \R^d\) and \(t \geq 0\).
    Together with 
    \begin{align*}
        X^{\eta}
        := \mathcal{X}(\widehat{h}(\eta, (R_i)_{i \in \N}), \Ncal) \quad \eta \in \Omega.
    \end{align*} this realizes the infinite-volume and finite-volume dynamics using a common Poisson random measure \(\mathcal{N}\) for the birth proposals and common i.i.d.\ exponentially distributed lifespans for shared particles of the initial configuration.
\end{proposition}
\begin{proof}
    The rate \(\widetilde{b}_\Lambda\) inherits measurability and boundedness from \(b\), and it is of finite range by the boundedness of \(\Lambda\). The solution map \(\mathcal{X}^\Lambda\) is therefore provided by \zcref{proposition:birth_death_dynamics_via_solution_map}. Using the same marking map \(\widehat{h}\) and the same driving noise as for \(b\), and intersecting the probability-one events for the two birth rates, yields the assertion.
\end{proof}

We can now give a 
\begin{proof}[Proof of \zcref{proposition:SDE_local}]
    \textit{(i) and (ii)}: Strong existence, pathwise uniqueness, uniqueness in law and the Markov properties follow from \zcref{proposition:SDE} via the correspondence between solutions of \eqref{equation:finite_volume_SDE_on_Omega_frozen_exterior} and \eqref{equation:SDE_local}. There are only finitely many initial particles and only finitely many relevant (i.e.\ with \(u \leq \Vert b\Vert_\infty\) for \((x, s, r, u) \in \Ncal\)) birth proposals in \(\Lambda\) during a bounded time interval \([0, T]\), hence only finitely many jumps.

    \textit{(iii)}: We use the same correspondence. Let \(\widetilde{\mathscr{L}}_\Lambda\) denote the infinite-volume generator associated with the birth rate \(\widetilde{b}_\Lambda\). If we view a solution of the finite-volume martingale problem as an \(\Omega\)-valued process (not adding any particles to \(\Lambda^c\)) and similarly an initial distribution \(\mu_\Lambda\) on \(\Omega_\Lambda\) as a probability measure \(\mu\) on \(\Omega\), then it solves the martingale problem for \((\widetilde{\mathscr{L}}_\Lambda, \mathcal{D}, \mu)\) as in \eqref{definition:martingale_problem}. This already yields uniqueness in the finite-volume martingale problem via \zcref{proposition:martingale_problem} then. The identification of the solution of the infinite-volume martingale problem with the law of the solution to the corresponding SDE in \zcref{proposition:SDE}\textit{(iii)}, and using the correspondence in the other direction (i.e.\ \(X_t^{(\Lambda, \eta)}
    = (X_t^{\Lambda, \eta})_\Lambda\)) gives that
    \begin{align*}
        P^{\Lambda, \mu}
        = \int \mu_\Lambda(\d\eta) \, \mathrm{Law}(X^{(\Lambda, \eta)})
    \end{align*} solves the martingale problem for \((\mathscr{L}_\Lambda, \mathcal{D}_\Lambda, \mu_\Lambda)\).
\end{proof}

\subsection{Semigroup properties of the finite-volume dynamics}\label{section:finite_volume_dynamics_semigroup}

In this section we will investigate some of the properties of the transition semigroup \((\Tl_t)_{t \geq 0}\) and reversibility w.r.t.\ \(\nu_\Lambda\).
We can already give a simple
\begin{proof}[Proof of \zcref{proposition:finite_volume_dynamics_semigroup}\textit{(i)}]
    This follows directly from \zcref{proposition:SDE_local}\textit{(ii)}.
\end{proof}

Let us now turn to the more elaborate statement of \zcref{proposition:finite_volume_dynamics_semigroup}\textit{(ii)}. For that, we need a couple of intermediary observations.

First, we record the pointwise forward and backward equations for bounded, measurable functions.
\begin{lemma}\label{lemma:finite_volume_process_forward_backward_equation_local_bounded_measurable}
    Let \(f \colon \Omega_\Lambda \to \R\) be bounded and measurable.
    Then, for every \(\eta \in \Omega_\Lambda\), \(t \mapsto \Tl_t f(\eta)\) is continuously differentiable on \([0, \infty)\) with derivative
    \begin{align*}
        t \mapsto \Tl_t (\mathscr{L}_\Lambda f)(\eta) = \mathscr{L}_\Lambda(\Tl_t f)(\eta).
    \end{align*} 
    The derivative in \(t = 0\) is interpreted as right-derivative.
    In particular, the forward and backward equations hold:
    \begin{align}\label{equation:finite_volume_process_forward_backward_equation_local_bounded_measurable}
        \Tl_t f(\eta) - f(\eta)
        = \int_0^t \, \Tl_s(\mathscr{L}_\Lambda f)(\eta) \, \d s
        = \int_0^t \, \mathscr{L}_\Lambda(\Tl_s f)(\eta) \, \d s, 
        \quad t \geq 0, \,\eta \in \Omega_\Lambda.
    \end{align}
    Additionally, we have the bound
    \begin{align}\label{equation:bound_generator_bounded_measurable}
        \vert \mathscr{L}_\Lambda f(\eta) \vert
        \leq 2 \Vert f\Vert_\infty \big(\Vert b\Vert_\infty \abs{\Lambda} + N_\Lambda(\eta)\big)
    \end{align}
\end{lemma}
\begin{proof}
    Extend \(f\) to \(\widetilde{f} \in \Lcal\) by \(\widetilde{f}(\eta) = f(\eta_\Lambda)\) for \(\eta \in \Omega\).
    Recall the auxiliary rate \(\widetilde{b}_\Lambda\), the associated formal generator \(\widetilde{\mathscr{L}}_\Lambda\) and the correspondence \(X_t^{(\Lambda, \eta)}
    = (X_t^{\Lambda, \eta})_\Lambda\), from \zcref{section:finite_volume_dynamics_construction}.
    Applying \zcref{lemma:infinite_volume_process_forward_equation_local_bounded_measurable} to the infinite-volume process yields continuous differentiability and 
    \begin{align*}
        \frac{\d}{\d t^+} \Tl_t f(\eta)
        = \frac{\d}{\d t^+} \mathbb{E}[\widetilde{f}(X_t^{\Lambda, \eta})]
        = \mathbb{E}[\widetilde{\mathscr{L}}_\Lambda\widetilde{f}(X_t^{\Lambda, \eta})]
        = \mathbb{E}[\mathscr{L}_\Lambda f(X_t^{(\Lambda, \eta)})]
        = \Tl_t( \mathscr{L}_\Lambda f)(\eta), \quad t \geq 0.
    \end{align*}
    Fix \(t \geq 0\). Since \(\Tl_t f\) is again bounded and measurable, the classical proof of the backward equation from the forward equation using the semigroup property (of \((\Tl_t)_{t \geq 0}\) on the bounded, measurable functions \(\Omega_\Lambda \to \R\)) applies. Indeed, by the computation above, 
    \begin{align*}
        \mathscr{L}_\Lambda (\Tl_t f)(\eta)
        = \lim_{h \downarrow 0} \frac{\Tl_h  (\Tl_t f)(\eta) -  (\Tl_t f)(\eta) }{h}
        = \lim_{h \downarrow 0} \frac{\Tl_{t+h} f(\eta) -  \Tl_t f(\eta) }{h}
        = \Tl_t(\mathscr{L}_\Lambda f)(\eta).
    \end{align*}

    The bound \zcref{equation:bound_generator_bounded_measurable} is just an application of the triangle inequality.
\end{proof}

Regarding the reversibility statement \zcref{proposition:finite_volume_reversibility}, let us first isolate the following infinitesimal reversibility and moment assertion.
\begin{lemma}\label{lemma:finite_volume_infinitesimal_reversibility}
    For all bounded, measurable \(f, g \colon \Omega_\Lambda \to \R\) it holds that
    \begin{align}\label{equation:finite_volume_reversibility}
        \nu_\Lambda[f (\mathscr{L}_\Lambda g) ]
        = \nu_\Lambda[(\mathscr{L}_\Lambda f) g].
    \end{align}
    The expectations are well-defined by \zcref{equation:bound_generator_bounded_measurable} and
    \begin{align}\label{equation:gibbs_marginal_first_moment}
        \nu_\Lambda[N_\Lambda] < \infty.
    \end{align}

    Furthermore, \eqref{equation:finite_volume_reversibility} holds also in the situation that \(f \in L^1(\nu_\Lambda)\) and \(g \in L^\infty(\nu_\Lambda)\) is such that \(g = 0\) on \(\{N_\Lambda > m\}\) for some \(m \in \N\). 
    Herein, \(\mathscr{L}_\Lambda f\) is \(\nu_\Lambda\)-a.e.\ well-defined for \(f \in L^1(\nu_\Lambda)\) and does not depend on the choice of the \(L^1(\nu_\Lambda)\)-representative.
\end{lemma}
\begin{proof}
    The GNZ equations imply
    \begin{align*}
        \nu_\Lambda[N_\Lambda] 
        = \int \nu(\d\eta) \int_{\R^d} \d x \, b(x, \eta) {\1}_{\Lambda}(x)
        \leq \Vert b\Vert_\infty \abs{\Lambda}.
    \end{align*}
    Now, denoting by \(\widetilde{f}\), \(\widetilde{g}\) the lifts of bounded, measurable \(f, g\) to \(\Omega\) with \(\widetilde{f}(\eta) = f(\eta_\Lambda)\) and \(\widetilde{g}(\eta) = g(\eta_\Lambda)\) for all \(\eta \in \Omega\), we get
    \begin{align*}
        \nu_\Lambda[f (\mathscr{L}_\Lambda g) ]
        = \nu[\widetilde{f} (\mathscr{L} \widetilde{g}) ]
        = \nu[(\mathscr{L} \widetilde{f}) \widetilde{g}]
        = \nu_\Lambda[(\mathscr{L}_\Lambda f) g]
    \end{align*} by the definition \(b_\Lambda^\nu(x, \eta_\Lambda) = \int \nu_{\Lambda \,\vert\, \Lambda^c}(\d\zeta_{\Lambda^c} \vert \eta_\Lambda)\, b(x, \eta_\Lambda \zeta_{\Lambda^c})\) and the infinite-volume infinitesimal reversibility established in the proof of \zcref{proposition:equivalence_of_GNZ_Gibbs_infinitesimally_reversible}.

    Let us now first see that \(\mathscr{L}_\Lambda f\) is \(\nu_\Lambda\)-a.e.\ well-defined for all \(f \in L^1(\nu_\Lambda)\) by ensuring that the sums and integrals in the definition of \(\mathscr{L}_\Lambda f(\eta)\) converge absolutely for \(\nu_\Lambda\)-a.a.\ \(\eta \in \Omega_\Lambda\), and that \(\mathscr{L}_\Lambda f\) does not depend on the representative of \(f\) in \(L^1(\nu_\Lambda)\). Indeed, for \(m \in \N\), we have 
    \begin{align}\label{equation:absolute_integrability_of_generator_expression_on_increasing_events}
        &\int \nu_\Lambda(\d\eta_\Lambda) {\1}_{\{N_\Lambda \leq m\}}(\eta_\Lambda) \Big\{ \int_{\Lambda} \d x \, b_\Lambda^\nu(x, \eta_\Lambda) \abs{f(\eta_\Lambda + \delta_x) - f(\eta_\Lambda)}  + \int_\Lambda \eta_\Lambda(\d x) \abs{f(\eta_\Lambda - \delta_x) - f(\eta_\Lambda)} \Big\} \\
        &\leq (2\Vert b\Vert_\infty \abs{\Lambda} + 2m + 1)\Vert f \Vert_{L^1(\nu_\Lambda)}, \nonumber
    \end{align} 
    because
     \begin{align*}
        &\int \nu_\Lambda(\d\eta_\Lambda) {\1}_{\{N_\Lambda \leq m\}}(\eta_\Lambda) \Big\{ \int_{\Lambda} \d x \, b_\Lambda^\nu(x, \eta_\Lambda) \abs{f(\eta_\Lambda + \delta_x)}  + \int_\Lambda \eta_\Lambda(\d x) \abs{f(\eta_\Lambda - \delta_x)} \Big\} \\
        &= \int \nu(\d\eta) \Big\{ \int_{\R^d} \d x \, b(x, \eta) \abs{f(\eta_\Lambda + \delta_x)} {\1}_\Lambda(x) {\1}_{\{N_\Lambda \leq m+1\}}(\eta_\Lambda + \delta_x)  \\
        &\qquad\qquad\qquad + \int \eta(\d x) \abs{f(\eta_\Lambda - \delta_x)}  {\1}_\Lambda(x) {\1}_{\{N_\Lambda \leq m -1 \}}(\eta_\Lambda - \delta_x)\Big\} \\
        &= \int \nu(\d\eta) \Big\{ \int \eta(\d x)  \abs{f(\eta_\Lambda)} {\1}_\Lambda(x) {\1}_{\{N_\Lambda \leq m+1\}}(\eta_\Lambda)  \\
        &\qquad\qquad\qquad + \int_{\R^d} \d x \, b(x, \eta) \abs{f(\eta_\Lambda)}  {\1}_\Lambda(x) {\1}_{\{N_\Lambda \leq m - 1\}}(\eta_\Lambda) \Big\} \\
        &\leq (m+1 +  \Vert b\Vert_\infty \abs{\Lambda})\Vert f\Vert_{L^1(\nu_\Lambda)} 
    \end{align*} by the definition of \(b_\Lambda^\nu(x, \eta_\Lambda) = \int \nu_{\Lambda \,\vert\, \Lambda^c}(\d\zeta_{\Lambda^c} \vert \eta_\Lambda)\, b(x, \eta_\Lambda \zeta_{\Lambda^c})\) and the GNZ equations. This shows that for \(\nu_\Lambda\)-integrable \(f\), \(\mathscr{L}_\Lambda f(\eta)\) converges absolutely for \(\nu_\Lambda\)-a.a.\ \(\eta \in \Omega_\Lambda\). Moreover, if \(\widetilde{f} =f \) \(\nu_\Lambda\)-a.s.\ it follows that \(\mathscr{L}_\Lambda \widetilde{f} = \mathscr{L}_\Lambda f\) \(\nu_\Lambda\)-a.s.\ by replacing \(f\) with \(f - \widetilde{f}\) in \eqref{equation:absolute_integrability_of_generator_expression_on_increasing_events}.

    Now suppose additionally to \(f \in L^1(\nu_\Lambda)\) that \(g \in L^\infty(\nu_\Lambda)\) is such that \(g = 0\) on \(\{N_\Lambda > m\}\) for fixed \(m \in \N\). Then, letting \(f_n := n \land (f \lor (-n))\), it holds
    \begin{align*}
        \nu_\Lambda[f (\mathscr{L}_\Lambda g)]
        = \lim_{n \to \infty} \nu_\Lambda[f_n (\mathscr{L}_\Lambda g)]
        = \lim_{n \to \infty} \nu_\Lambda[(\mathscr{L}_\Lambda f_n) g]
        = \nu_\Lambda[(\mathscr{L}_\Lambda f) g].
    \end{align*} The first equality is justified by the properties of \(g\) and \(L^1(\nu_\Lambda)\)-convergence of \(f_n \xrightarrow[n \to \infty]{} f\), the second equality is the infinitesimal reversibility for bounded functions from above, and the third equality follows from \(\abs{g} \leq \Vert g\Vert_\infty {\1}_{N_\Lambda \leq m}\) and \(L^1(\nu_\Lambda)\)-convergence \((\mathscr{L}_\Lambda f_n) {\1}_{N_\Lambda \leq m} \xrightarrow[n \to \infty]{} (\mathscr{L}_\Lambda f) {\1}_{N_\Lambda \leq m}\) by \eqref{equation:absolute_integrability_of_generator_expression_on_increasing_events}.
\end{proof}

We can now, with a classical argument, leverage this infinitesimal reversibility together with the backward equation from \zcref{lemma:finite_volume_process_forward_backward_equation_local_bounded_measurable} to give a
\begin{proof}[Proof of \zcref{proposition:finite_volume_reversibility}]
    By straightforward monotone approximation it is enough to consider bounded measurable \(f, g \colon \Omega_\Lambda \to \R\).
    Fix \(t > 0\) and define \(\phi \colon [0, t] \to \R\) by
    \begin{align*}
        \phi(s) 
        := \nu_\Lambda[(\Tl_s f)(\Tl_{t-s} g)].
    \end{align*}
    The function \(\phi\) has the property that \(\phi(0) = \nu_\Lambda[f(\Tl_{t} g)]\) and \(\phi(t) = \nu_\Lambda[(\Tl_{t} f) g]\), so that it suffices to show that \(\phi\) is continuous on \([0, t]\) and that it holds \(\phi'(s) = 0\) for \(s \in (0, t)\).
    Continuity of \(\phi\) follows from dominated convergence and the pointwise continuity of \(s \mapsto (\Tl_s f)(\Tl_{t-s} g)\) which follows e.g.\ from \zcref{lemma:finite_volume_process_forward_backward_equation_local_bounded_measurable}.
    Regarding the derivative of \(\phi\), we compute
    \begin{align}\label{equation:derivative_helper_forward_backward_function}
        \frac{\d}{\d s} (\Tl_s f)(\Tl_{t-s} g)
        = (\mathscr{L}_\Lambda \Tl_s f) (\Tl_{t-s} g) - (\Tl_s f)(\mathscr{L}_\Lambda \Tl_{t-s} g)
    \end{align}
    Applying the infinitesimal reversibility \zcref{lemma:finite_volume_infinitesimal_reversibility}, we get indeed
    \begin{align*}
        \phi'(s) = \nu_\Lambda[(\mathscr{L}_\Lambda \Tl_s f) (\Tl_{t-s} g) - (\Tl_s f)(\mathscr{L}_\Lambda \Tl_{t-s} g)]
        = 0,
    \end{align*} if we can justify the exchange of derivative and expectation.
    This is not so hard, as the dominated convergence theorem is applicable. Indeed, we get by \zcref{equation:derivative_helper_forward_backward_function} and \zcref{equation:bound_generator_bounded_measurable}, that
    \begin{align*}
        \Big\vert \frac{\d}{\d s} (\Tl_s f)(\Tl_{t-s} g) \Big\vert
        \leq 4 \Vert f\Vert_\infty \Vert g\Vert_\infty \big(\Vert b\Vert_\infty \abs{\Lambda} + N_\Lambda \big)
    \end{align*} and \(N_\Lambda\) is \(\nu_\Lambda\)-integrable by \zcref{equation:gibbs_marginal_first_moment}.
\end{proof}

Now with reversibility of the finite-volume dynamics w.r.t.\ \(\nu_\Lambda\) and the pointwise forward equation for bounded measurable functions at hand, we give a
\begin{proof}[Proof of \zcref{proposition:finite_volume_dynamics_semigroup}\textit{(ii)}]
    \emph{Extension to a \(C_0\)-semigroup on \(L^1(\nu_\Lambda)\):}
    By reversibility of \(\nu_\Lambda\) established in \zcref{proposition:finite_volume_reversibility}, we have
    \begin{align*}
        \Vert \Tl_t f\Vert_{L^1(\nu_\Lambda)}
        =  \nu_\Lambda[\abs{\Tl_t f}]
        \leq \nu_\Lambda[(\Tl_t \abs{f}) 1]
        = \nu_\Lambda[\abs{f} (\Tl_t  1)]
        = \Vert f\Vert_{L^1(\nu_\Lambda)},
    \end{align*} a priori for bounded measurable \(f\), so that the operators \((\Tl_t)_{t \geq 0}\) extend uniquely to \(L^1(\nu_\Lambda)\)-contractions. By \zcref{lemma:finite_volume_process_forward_backward_equation_local_bounded_measurable} and dominated convergence, 
    \begin{align*}
        \Vert \Tl_t f - f\Vert_{L^1(\nu_\Lambda)}
        \xrightarrow[t \downarrow 0]{} 0
    \end{align*} for bounded measurable \(f\). This convergence extends to general \(f \in L^1(\nu_\Lambda)\) by density of bounded measurable \(f\) in \(L^1(\nu_\Lambda)\) and the contraction property of the semigroup. Consequently, \((\Tl_t)_{t \geq 0}\) is strongly continuous on \(L^1(\nu_\Lambda)\).

    \emph{Properties of the generator:}

    Now denote the generator of \((\Tl_t)_{t \geq 0}\) on \(L^1(\nu_\Lambda)\) by \(\mathfrak{L}_\Lambda\). Fix bounded measurable \(f\). Then \(\mathscr{L}_\Lambda f \in L^1(\nu_\Lambda)\) by \eqref{equation:bound_generator_bounded_measurable} and \eqref{equation:gibbs_marginal_first_moment}, and the forward equation \eqref{equation:finite_volume_process_forward_backward_equation_local_bounded_measurable} together with Minkowski's integral inequality yield
    \begin{align*}
        &\Big\Vert \frac{\Tl_t f - f}{t} - \mathscr{L}_\Lambda f \Big\Vert_{L^1(\nu_\Lambda)}
        = \Big\Vert \frac{1}{t} \int_0^t \d s \, [\Tl_s \mathscr{L}_\Lambda f - \mathscr{L}_\Lambda f] \Big\Vert_{L^1(\nu_\Lambda)} \\
        &\leq \frac{1}{t} \int_0^t \d s \, \Vert \Tl_s \mathscr{L}_\Lambda f - \mathscr{L}_\Lambda f\Vert_{L^1(\nu_\Lambda)}
        \xrightarrow[t \downarrow 0]{} 0.
    \end{align*}
    Hence, \(f\) is an element of the domain of \(\mathfrak{L}_\Lambda\) and \(\mathfrak{L}_\Lambda f = \mathscr{L}_\Lambda f\) \(\nu_\Lambda\)-a.s.\ 

    It remains to identify the action of \(\mathfrak{L}_\Lambda\) on its whole domain as given by \(\mathscr{L}_\Lambda\). For that, fix \(f \in L^1(\nu_\Lambda)\) from the domain of \(\mathfrak{L}_\Lambda\).
    We use the existing identification on the bounded measurable functions and reversibility and infinitesimal reversibility from \zcref{proposition:finite_volume_reversibility}, \zcref{lemma:finite_volume_infinitesimal_reversibility} to show that
    \begin{align}\label{equation:identification_of_generator}
        &\nu_\Lambda[(\mathfrak{L}_\Lambda f) g]
        = \lim_{t \downarrow 0} \nu_\Lambda\Big[\Big( \frac{\Tl_t f - f}{t}\Big) g\Big]
        = \lim_{t \downarrow 0} \nu_\Lambda\Big[f \Big(\frac{\Tl_t g - g}{t} \Big)\Big]\\
        &= \nu_\Lambda\Big[f (\mathscr{L}_\Lambda g)\Big]
        = \nu_\Lambda\Big[(\mathscr{L}_\Lambda f) g \Big] \nonumber
    \end{align} for all bounded and measurable \(g\) with \(g = 0\) on \(\{N_\Lambda > m\}\) for some \(m \in \N\). Hence, \(\mathfrak{L}_\Lambda f\) and \(\mathscr{L}_\Lambda f\) agree \(\nu_\Lambda\)-a.s.\  on \(\{N_\Lambda \leq m\}\) for all \(m \in \N\), and as these events exhaust \(\Omega_\Lambda\), it holds  \(\mathfrak{L}_\Lambda f = \mathscr{L}_\Lambda f\)  \(\nu_\Lambda\)-a.s.\ 
    The only thing left to justify is the third equality in \zcref{equation:identification_of_generator}, and it follows by dominated convergence and the bound
    \begin{align*}
        \Big\Vert \frac{\Tl_t g - g}{t} \Big\Vert_\infty
        = \Big\Vert \frac{1}{t} \int_0^t \d s \, \Tl_s \mathscr{L}_\Lambda g \Big\Vert_\infty
        \leq \frac{1}{t} \int_0^t \d s \, \Vert \Tl_s \abs{\mathscr{L}_\Lambda g} \Vert_\infty
        \leq 2\Vert g\Vert_\infty(\Vert b\Vert_\infty \abs{\Lambda} + (m+1)),
    \end{align*} by
    \begin{align*}
        \abs{\mathscr{L}_\Lambda g}
        \leq 2\Vert g\Vert_\infty(\Vert b\Vert_\infty \abs{\Lambda} + (m+1)),
    \end{align*}
    for bounded and measurable \(g\) with \(g = 0\) on \(\{N_\Lambda > m\}\).
\end{proof}

It is important for our application to have the forward equation also for a larger class than bounded, measurable functions, namely the densities \(f = \d \mu_\Lambda / \d \nu_\Lambda\) of probabilities \(\mu_\Lambda\) having finite entropy w.r.t.\ \(\nu_\Lambda\). 
\begin{lemma}\label{lemma:L1_differentiability_from_finite_local_entropy}
    Let \(\Lambda \Subset \R^d\) and \(f \colon \Omega_\Lambda \to [0, \infty]\) measurable with
    \begin{align}
        \nu_\Lambda[f \log f] < \infty,
    \end{align} and write \(f_t := \Tl_t f\) for \(t \geq 0\).
    Then, 
    \begin{enumerate}[(i)]
        \item \(\mathscr{L}_\Lambda f_t \in L^1(\nu_\Lambda)\) for all \(t \geq 0\),
        
        \item the map \([0, \infty) \ni t \mapsto f_t\) is continuously differentiable in \(L^1(\nu_\Lambda)\) with derivative \([0, \infty) \ni t \mapsto \mathscr{L}_\Lambda f_t\).

        \item For any continuously differentiable \(\psi \colon [0, \infty) \to \R\) with bounded derivative \(\psi'\) it holds
        \begin{align}
            \nu_\Lambda[\psi(f_t)] - \nu_\Lambda[\psi(f)] 
            = \int_0^{t} \d s \, \nu_\Lambda[\psi'(f_s) \mathscr{L}_\Lambda f_s], \quad t \geq 0.
        \end{align}
    \end{enumerate}
\end{lemma}

\begin{proof}[Proof of \zcref{lemma:L1_differentiability_from_finite_local_entropy}]
    \textit{(i) and (ii)}: The statements follow from the fact that \((\Tl_t)_{t \geq 0}\) is a \(C_0\)-semigroup on \(L^1(\nu_\Lambda)\) shown in \zcref{proposition:finite_volume_dynamics_semigroup} and noticing that \(f\) is contained in the domain of the generator \(\mathfrak{L}_\Lambda\) of  \((\Tl_t)_{t \geq 0}\) on \(L^1(\nu_\Lambda)\) and the identification of \(\mathfrak{L}_\Lambda\) and \(\mathscr{L}_\Lambda\) given there. Let us validate the assertion that \(f\) is contained in the domain, which in turn follows from \zcref{lemma:in_domain_by_integrability_against_counting_variable} if it holds \(\nu_\Lambda[f N_\Lambda] < \infty\). For this, we use the entropy inequality (e.g.\ Fenchel--Young)
    \begin{align*}
        \nu_\Lambda[f N_\Lambda] 
        \leq \nu_\Lambda[f \log f] + \nu_\Lambda[\e^{N_\Lambda}] 
        < \infty
    \end{align*} 
    and, writing here \((N_\Lambda)_n\) for the falling factorial of \(N_\Lambda\) and repeatedly using GNZ as in the proof of \zcref{equation:gibbs_marginal_first_moment},
    \begin{align}\label{equation:gibbs_marginal_exponential_moment}
        \nu_\Lambda[\e^{N_\Lambda}]
        &= \nu\Big[\sum_{n = 0}^\infty \frac{(\e -1)^n}{n!} (N_\Lambda)_n\Big]
        = \sum_{n = 0}^\infty \frac{(\e -1)^n}{n!} \nu[(N_\Lambda)_n]
        \leq  \sum_{n = 0}^\infty \frac{(\e -1)^n}{n!} (\Vert b\Vert_\infty\abs{\Lambda})^n \\
        &= \exp((\e - 1) \Vert b\Vert_\infty\abs{\Lambda})
        < \infty. \nonumber
    \end{align}

    \textit{(iii)}: This is a slight variation of a well-known chain rule, we give the short proof anyways. It is enough to show that \(t \mapsto \nu_\Lambda[\psi(f_t)]\) is continuously differentiable with derivative \(t \mapsto \nu_\Lambda[\psi'(f_t) \mathscr{L}_\Lambda f_t]\). For that, consider the pointwise difference quotient
    \begin{align*}
        \frac{\psi(f_{t+h}) - \psi(f_t)}{h} 
        = \Big(\frac{\psi(f_{t+h}) - \psi(f_t)}{f_{t+h} - f_t} {\1}_{f_{t+h} - f_t \neq 0} + \psi'(f_t) {\1}_{f_{t+h} - f_t = 0} \Big) \frac{f_{t+h} - f_t}{h} 
    \end{align*} The second factor in the product of the right-hand side converges to \(\mathscr{L}_\Lambda f_t\) in \(L^1(\nu_\Lambda)\). The first factor converges in \(\nu_\Lambda\)-probability to \(\psi'(f_t)\) by the mean-value theorem, continuity of \(\psi'\), \(f_{t+h} \xrightarrow[h \downarrow 0]{} f_t\) in \(L^1(\nu_\Lambda)\). Additionally, the first factor is bounded in absolute value by \(\Vert \psi'\Vert_\infty\), again by the mean-value theorem. It follows that the product converges to \(\psi'(f_t) \mathscr{L}_\Lambda f_t\) in \(L^1(\nu_\Lambda)\), i.e.\
    \begin{align*}
        \frac{\nu_\Lambda[\psi(f_{t+h})] - \nu_\Lambda[\psi(f_t)]}{h} 
        = \nu_\Lambda\Big[\frac{\psi(f_{t+h}) - \psi(f_t)}{h} \Big]
        \xrightarrow[h \downarrow 0]{} \nu_\Lambda[\psi'(f_t) \mathscr{L}_\Lambda f_t].
    \end{align*}
    Continuity of the derivative \(t \mapsto \nu_\Lambda[\psi'(f_t) \mathscr{L}_\Lambda f_t]\) then follows with essentially the same argumentation.
\end{proof}

We needed the following
\begin{lemma}\label{lemma:in_domain_by_integrability_against_counting_variable}
    Let \(f \colon \Omega_\Lambda \to \R\) be measurable and such that 
    \begin{align*}
        \nu_\Lambda[N_\Lambda \abs{f}] < \infty.
    \end{align*}
    Then, \(f\) is contained in the domain of the generator \(\mathfrak{L}_\Lambda\) of  \((\Tl_t)_{t \geq 0}\) on \(L^1(\nu_\Lambda)\) and \(\mathfrak{L}_\Lambda f = \mathscr{L}_\Lambda f\) \(\nu_\Lambda\)-a.s.\
\end{lemma}
\begin{proof}
    Notice first that \(f \in L^1(\nu_\Lambda)\) by
    \begin{align*}
        \nu_\Lambda[\abs{f}] 
        \leq \nu_\Lambda[(N_\Lambda + {\1}_\emptyset)\abs{f}] 
        \leq \nu_\Lambda[N_\Lambda \abs{f}] + \nu_\Lambda(\{\emptyset\}) \abs{f}(\emptyset) < \infty.
    \end{align*}
    The identification \(\mathfrak{L}_\Lambda f = \mathscr{L}_\Lambda f\) for all \(f\) in the domain of the \(L^1(\nu_\Lambda)\) generator \(\mathfrak{L}_\Lambda\) was already established in \zcref{proposition:finite_volume_dynamics_semigroup}.
    By closedness of the generator of a \(C_0\)-semigroup it also suffices to show that for \(f_n := n \land (f \lor (-n))\), \(n \in \N\), we have \(\mathscr{L}_\Lambda f_n \xrightarrow[n \to \infty]{} \mathscr{L}_\Lambda f\) in \(L^1(\nu_\Lambda)\).
    But \(\mathscr{L}_\Lambda f \in L^1(\nu_\Lambda)\) and \(\mathscr{L}_\Lambda f_n \xrightarrow[n \to \infty]{} \mathscr{L}_\Lambda f\) in \(L^1(\nu_\Lambda)\) follow from 
    \begin{align*}
        \nu_\Lambda[\abs{\mathscr{L}_\Lambda g}]
        \leq 2 \nu_\Lambda[(\Vert b\Vert_\infty \abs{\Lambda} + N_\Lambda) \abs{g}]
    \end{align*} for \(g \in L^1(\nu_\Lambda)\), which follows from the GNZ equations exactly as in \zcref{equation:absolute_integrability_of_generator_expression_on_increasing_events}.
\end{proof}

\bibliographystyle{alpha}
\bibliography{references}

\end{document}